%% file: har-rev.tex
\documentclass{amsart}%
\usepackage{hyperref}
\usepackage{amssymb}
\usepackage{amsmath}
\usepackage{amsthm}
\usepackage{amsfonts}
\usepackage{graphicx}%
\allowdisplaybreaks
\makeatletter
\@ifundefined{pdfoutput}
{
\DeclareGraphicsRule{.wmf}{bmp}{}{}
\DeclareGraphicsRule{.jpg}{bmp}{}{}
\DeclareGraphicsRule{.png}{bmp}{}{}
\DeclareGraphicsRule{.cdr}{bmp}{}{}
\DeclareGraphicsRule{.gif}{bmp}{}{}
}
{
\ifnum\pdfoutput=0\relax
\DeclareGraphicsRule{.wmf}{bmp}{}{}
\DeclareGraphicsRule{.jpg}{bmp}{}{}
\DeclareGraphicsRule{.png}{bmp}{}{}
\DeclareGraphicsRule{.cdr}{bmp}{}{}
\DeclareGraphicsRule{.gif}{bmp}{}{}
\fi
\ifnum\pdfoutput=1\relax
\DeclareGraphicsRule{.eps}{pdf}{}{}
\DeclareGraphicsRule{.wmf}{jpg}{}{}
\fi
}
\makeatother
\newtheorem{thm}{Theorem}[section]
\newtheorem{cor}[thm]{Corollary}
\newtheorem{ex}[thm]{Example}

\newtheorem{assumption}{Assumption}
\newtheorem{lem}[thm]{Lemma}
\newtheorem{nota}[thm]{Notation}
\newtheorem{prop}[thm]{Proposition}

\theoremstyle{definition}
\newtheorem*{acknowledgement}{Acknowledgement}
\newtheorem{df}[thm]{Definition}

\newtheorem{example}{Example}[section]
\newtheorem{rem}[thm]{Remark}
\newtheorem{rems}[thm]{Remarks}
\newtheorem{notation}[thm]{Notation}
\numberwithin{equation}{section}
\begin{document}

\input{har-rev-title.tex}

\section{Introduction \label{s.1}}

\subsection{Basic setup\label{s.1.1}}

Let $(M,g)$ be a connected complete Riemannian manifold, $d:M\times
M\rightarrow\lbrack0,\infty)$ be the Riemannian distance function, $dV$ be the
Riemannian volume measure on $M,$ $\Delta$ be the Laplace--Beltrami operator
acting on the space of smooth differential forms, $\Omega\left(  M\right)  ,$
over $M,$ and $\Delta_{0}:=\Delta|_{\Omega_{c}^{0}\left(  M\right)  },$ where
$\Omega_{c}^{0}\left(  M\right)  :=C_{c}^{\infty}\left(  M\right)  $ is the
space of compactly supported smooth functions on $M.$ From Gaffney
\cite{Gaffney51}, Roelcke \cite{Roelcke60}, Chernoff \cite{Chernoff73} and
Strichartz \cite{Str1983}, we know that the $L^{2}\left(  M,dV\right)
$--closure, $\bar{\Delta}_{0},$ of $\Delta_{0}$ is a non-positive self-adjoint
operator on $L^{2}\left(  M,dV\right)  .$ Moreover, there exists an associated
smooth heat kernel, $\left(  0,\infty\right)  \times M\times M\ni\left(
t,x,y\right)  \rightarrow p_{t}\left(  x,y\right)  \in\left(  0,\infty\right)
,$ such that $p_{t}\left(  x,y\right)  =p_{t}\left(  y,x\right)  ,$%
\begin{equation}
\int_{M}p_{t}(x,y)dV(y)\leq1\text{ for all }x\in M,\text{ and} \label{e.1.1}%
\end{equation}%
\begin{equation}
\left(  e^{t\bar{\Delta}_{0}/2}f\right)  (x)=\int_{M}p_{t}(x,y)f(y)dV(y)\text{
for all }f\in L^{2}(M). \label{e.1.2}%
\end{equation}
For the bulk of this paper we will be considering the special case where $M=G$
is a Lie group equipped with a left invariant Riemannian metric as we now describe.

Let $G$ be a connected finite dimensional uni-modular Lie group,
$\mathfrak{g}=\mathrm{Lie}\left(  G\right)  $ be its Lie algebra, and suppose
that $\mathfrak{g}$ is equipped with an inner product, $\left(  \cdot
,\cdot\right)  =\left(  \cdot,\cdot\right)  _{\mathfrak{g}}.$ Let $\left\vert
A\right\vert _{\mathfrak{g}}:=\sqrt{\left(  A,A\right)  }$ for all
$A\in\mathfrak{g}.$ We endow $G$ with the unique left invariant Riemannian
metric which agrees with $\left(  \cdot,\cdot\right)  _{\mathfrak{g}}$ at
$e\in G,$ i.e. the unique metric on $G\ $such that $L_{g\ast}:\mathfrak{g}%
\rightarrow T_{g}G$ is isometric for all $g\in G.$ The Riemannian distance
between $x,y\in G$ will be denoted by $d\left(  x,y\right)  .$

For $A\in\mathfrak{g}$ let $\tilde{A}$ denote the unique left invariant vector
field on $G\ $such that $\tilde{A}\left(  e\right)  =A\in\mathfrak{g}$ and let
$L=\sum_{i=1}^{\dim\mathfrak{g}}\tilde{A}_{i}^{2}$ where $\left\{
A_{i}\right\}  _{i=1}^{\dim\mathfrak{g}}$ is an orthonormal basis for
$\mathfrak{g}.$ As is well-known, since $G$ is uni-modular, $L$ is the
Laplace-Beltrami operator (for example, see \cite[Remark 2.2]{Driver1997c} and
Lemma \ref{l.6.1} below) restricted to $C^{\infty}\left(  G\right)  .$ Since
$L_{g}:G\rightarrow G$ is an isometry for all $g\in G,$ if $p_{t}\left(
x,y\right)  $ is the heat kernel on $G,$ then $p_{t}\left(  gx,gy\right)
=p_{t}\left(  x,y\right)  $ for all $x,y,g\in G.$ Taking $g=x^{-1}$ then
implies that $p_{t}\left(  x,y\right)  =p_{t}\left(  e,x^{-1}y\right)  .$
Similarly, $d\left(  gx,gy\right)  =d\left(  x,y\right)  $ for all $x,y,g\in
G$ and therefore $d\left(  x,y\right)  =d\left(  e,x^{-1}y\right)  .$

\begin{nota}
\label{n.1.1}By a slight abuse of notation, let $p_{t}\left(  x\right)
:=p_{t}\left(  e,x\right)  $ for $x\in G.$ We will refer to $p_{t}\left(
\cdot\right)  $ as the \textbf{convolution} \textbf{heat kernel} on $G$ and to
the probability measure, $d\nu_{t}\left(  x\right)  :=p_{t}\left(  x\right)
dx,$ as the \textbf{heat kernel measure on }$G.$ We also write $dx$ for
$dV\left(  x\right)  $ and $\left\vert x\right\vert $ for $d\left(
e,x\right)  .$
\end{nota}

The following lemma is an immediate consequence of the comments above and the
basic properties of $p_{t}\left(  x,y\right)  .$

\begin{lem}
\label{l.1.2}For all $x,y\in G$

\begin{enumerate}
\item $d\left(  x,y\right)  =\left\vert x^{-1}y\right\vert ,$

\item $\left\vert x^{-1}\right\vert =\left\vert x\right\vert $

\item $p_{t}\left(  x^{-1}\right)  =p_{t}\left(  x\right)  $

\item $p_{t}\left(  x,y\right)  =p_{t}\left(  x^{-1}y\right)  =p_{t}\left(
y^{-1}x\right)  ,$

\item $dV$ is a bi-invariant Haar measure on $G,$

\item for $f\in L^{2}\left(  G,dV\right)  ,$%
\begin{align*}
\left(  e^{t\bar{\Delta}_{0}/2}f\right)  (x)  &  =\int_{G}p_{t}\left(
x^{-1}y\right)  f\left(  y\right)  dy\\
&  =\int_{G}p_{t}\left(  y^{-1}x\right)  f\left(  y\right)  dy\\
&  =\int_{G}p_{t}\left(  yx\right)  f\left(  y^{-1}\right)  dy.
\end{align*}

\end{enumerate}
\end{lem}

\subsection{The main theorems\label{s.1.2}}

\begin{df}
\label{d.1.3}For $A\in\mathfrak{g}$ and $T>0,$ let
\[
W_{A}^{T}\left(  x\right)  :=-\left(  \tilde{A}\ln p_{T}\right)  \left(
x\right)  =-\frac{\left(  \tilde{A}p_{T}\right)  \left(  x\right)  }%
{p_{T}\left(  x\right)  }.
\]

\end{df}

The significance of $W_{A}^{T}$ in the above definition stems from the
following integration by parts identity;%
\begin{equation}
\int_{G}\tilde{A}f\left(  x\right)  p_{T}\left(  x\right)  dx=\int_{G}f\left(
x\right)  W_{A}^{T}\left(  x\right)  p_{T}\left(  x\right)  dx~\forall~f\in
C_{c}^{\infty}\left(  G\right)  . \label{e.1.3}%
\end{equation}

We may now state the main theorems of this paper.

\begin{thm}
\label{t.1.4}If $T>0$ and $A\in\mathfrak{g},$ then%
\begin{equation}
\int_{G}e^{W_{A}^{T}\left(  x\right)  }p_{T}\left(  x\right)  dx\leq
\exp\left(  \frac{c\left(  kT\right)  }{2T}\left\vert A\right\vert
_{\mathfrak{g}}^{2}\right)  , \label{e.1.4}%
\end{equation}
where $c\left(  \cdot\right)  $ is as in Eq. (\ref{e.1.7}).
\end{thm}

The proof of this theorem relies on martingale inequalities applied to the
probabilistic representation for $\tilde{A}\ln p_{T}\left(  x\right)  $ in
Theorem \ref{t.6.4}. We also have another related integral bound on $W_{A}%
^{T}.$

\begin{thm}
\label{t.1.5}Continuing the notation in Theorem \ref{t.1.6} and in particular
let $c\left(  \cdot\right)  $ be as in Eq. (\ref{e.1.7}). Then for any
$p\in\left(  1,\infty\right)  $ there is a constant, $C_{p}<\infty$ such that%
\begin{equation}
\left\Vert W_{A}^{T}\right\Vert _{L^{p}\left(  \nu_{T}\right)  }\leq
C_{p}\sqrt{\frac{c\left(  kT\right)  }{T}}\left\vert A\right\vert \text{ for
all }A\in\mathfrak{g}. \label{e.1.5}%
\end{equation}

\end{thm}

These theorems will be proved in Sections \ref{s.5} and \ref{s.6} below. Also
see \cite[Theorem 5.11]{Driver2001b} for a version of this theorem valid on a
general compact Riemannian manifold and Proposition \ref{p.E.1} in Appendix
\ref{s.E} where we use a Hamilton type inequality to show that an inequality
similar to that in Eq. (\ref{e.1.4}) holds on any complete Riemannian
manifolds whose Ricci curvature is bounded from below. However, see Remark
\ref{r.E.2} where it is noted that, in general, we can not choose the
constants appearing in Proposition \ref{p.E.1} to be independent of dimension.

The following theorem is a corollary of Theorem \ref{t.1.4} above and Theorem
\ref{t.2.5} below. The details will be given in Section \ref{s.3} below.

\begin{thm}
\label{t.1.6}Let $T>0$ be given and let $k\in\mathbb{R}$ be a lower bound on
the Ricci curvature, $\operatorname{Ric}\geq kI.$ Then for every $y\in G$ and
$p\in\lbrack1,\infty),$%
\begin{equation}
\left(  \int_{G}\left[  \frac{p_{T}\left(  xy^{-1}\right)  }{p_{T}\left(
x\right)  }\right]  ^{p}p_{T}\left(  x\right)  dx\right)  ^{1/p}\leq
\exp\left(  \frac{c\left(  kT\right)  \left(  p-1\right)  }{2T}\left\vert
y\right\vert ^{2}\right)  \label{e.1.6}%
\end{equation}
where%
\begin{equation}
c\left(  t\right)  =\frac{t}{e^{t}-1}~\text{ for all }~t\in\mathbb{R}
\label{e.1.7}%
\end{equation}
with the convention that $c\left(  0\right)  =1.$
\end{thm}

From Theorem \ref{t.1.6} and Lemma \ref{l.1.2} we have,%
\begin{align}
\left(  \int_{G}\left[  \frac{p_{T}\left(  y,x\right)  }{p_{T}\left(
z,x\right)  }\right]  ^{p}p_{T}\left(  z,x\right)  dx\right)  ^{1/p} &
=\left(  \int_{G}\left[  \frac{p_{T}\left(  y^{-1}x\right)  }{p_{T}\left(
z^{-1}x\right)  }\right]  ^{p}p_{T}\left(  z^{-1}x\right)  dx\right)
^{1/p}\nonumber\\
&  =\left(  \int_{G}\left[  \frac{p_{T}\left(  y^{-1}zx\right)  }{p_{T}\left(
x\right)  }\right]  ^{p}p_{T}\left(  x\right)  dx\right)  ^{1/p}\nonumber\\
&  \leq\exp\left(  \frac{c\left(  kT\right)  \left(  p-1\right)  }%
{2T}\left\vert y^{-1}z\right\vert ^{2}\right)  \nonumber\\
&  =\exp\left(  \frac{c\left(  kT\right)  \left(  p-1\right)  }{2T}%
d^{2}\left(  y,z\right)  \right)  \label{e.1.8}%
\end{align}
for all $y,z\in G.$ This form of the integrated Harnack inequality makes sense
on any Riemannian manifold. We will show in Corollary \ref{c.D.3} of Appendix
\ref{s.D} below that Eq. (\ref{e.1.8}) does indeed hold when $G$ is replaced
by a complete connected Riemannian manifold with $\operatorname*{Ric}\geq kI$
for some $k\in\mathbb{R}.$ The key point is that the estimate in Eq.
(\ref{e.1.8}) is equivalent to Wang's dimension free Harnack inequality, see
\cite{Wang97a,Wang2004} and Theorem \ref{t.D.2} below. We are grateful to
Michael R\"{o}ckner for pointing out the relationship between Wang's
inequality and the integrated Harnack inequality in Eq. (\ref{e.1.8}).

\begin{rems}
\label{r.1.7}Some of the key features of Theorem \ref{t.1.6} are:

\begin{enumerate}
\item As seen in Example \ref{ex.1.1}) below, the estimate in Eq.
(\ref{e.1.6}) is sharp when $G=\mathbb{R}^{n}.$

\item For $T$ near zero, $c\left(  kT\right)  /T\cong1/T$ and for $T$ large,
$c\left(  kT\right)  /T\cong\max\left(  0,-k\right)  .$

\item The estimate in Eq. (\ref{e.1.6}) is dimension independent and therefore
has applications to infinite dimensional settings, see Section \ref{s.7} below.
\end{enumerate}
\end{rems}

Let $R_{y}:G\rightarrow G$ ($L_{y}:G\rightarrow G)$ be the operation of right
(left) multiplication by $y\in G,$ $\nu_{T}\circ R_{y}^{-1}$ $\left(  \nu
_{T}\circ L_{y}^{-1}\right)  $ be $\nu_{T}$ pushed forward by $R_{y}$ $\left(
L_{y}\right)  ,$ and $d\left(  \nu_{T}\circ R_{y}^{-1}\right)  /d\nu_{T}$
denote the Radon-Nikodym derivative of $\nu_{T}\circ R_{y}^{-1}$ with respect
to $\nu_{T}.$ For the infinite dimensional applications of Section \ref{s.7},
it is convenient to rewrite Eq. (\ref{e.1.6}) as%
\begin{equation}
\left\Vert \frac{d\left(  \nu_{T}\circ R_{y}^{-1}\right)  }{d\nu_{T}%
}\right\Vert _{L^{p}\left(  G,\nu_{T}\right)  }\leq\exp\left(  \frac{c\left(
kT\right)  \left(  p-1\right)  }{2T}d^{2}\left(  e,y\right)  \right)  .
\label{e.1.9}%
\end{equation}
By Lemma \ref{l.1.2}, Eq. (\ref{e.1.6}) may be also be expressed as%
\begin{equation}
\left(  \int_{G}\left[  \frac{p_{T}\left(  xy\right)  }{p_{T}\left(  x\right)
}\right]  ^{p}p_{T}\left(  x\right)  dx\right)  ^{1/p}\leq\exp\left(
\frac{c\left(  kT\right)  \left(  p-1\right)  }{2T}\left\vert y\right\vert
^{2}\right)  \label{e.1.10}%
\end{equation}
or as%
\begin{equation}
\left(  \int_{G}\left[  \frac{p_{T}\left(  y^{-1}x\right)  }{p_{T}\left(
x\right)  }\right]  ^{p}p_{T}\left(  x\right)  dx\right)  ^{1/p}\leq
\exp\left(  \frac{c\left(  kT\right)  \left(  p-1\right)  }{2T}\left\vert
y\right\vert ^{2}\right)  . \label{e.1.11}%
\end{equation}
This last equality is equivalent to the left translation analogue of Eq.
(\ref{e.1.9}), namely%
\begin{equation}
\left\Vert \frac{d\left(  \nu_{T}\circ L_{y}^{-1}\right)  }{d\nu_{T}\left(
\cdot\right)  }\right\Vert _{L^{p}\left(  G,\nu_{T}\right)  }\leq\exp\left(
\frac{c\left(  kT\right)  \left(  p-1\right)  }{2T}\left\vert y\right\vert
^{2}\right)  . \label{e.1.12}%
\end{equation}

\subsection{Examples and applications\label{s.1.3}}

\begin{example}
\label{ex.1.1}Suppose $G=\mathbb{R}^{n}$ so that $\mathfrak{g}\cong
\mathbb{R}^{n}$ which we assume has been equipped with the standard inner
product. In this case%
\[
p_{T}\left(  x\right)  =\left(  \frac{1}{2\pi T}\right)  ^{n/2}\exp\left(
-\frac{\left\vert x\right\vert ^{2}}{2T}\right)  ,
\]
where $\left\vert x\right\vert ^{2}:=\sum_{i=1}^{n}x_{i}^{2}.$ For
$A\in\mathfrak{g}$ and $f\in C_{c}^{1}\left(  \mathbb{R}^{n}\right)  $ we have
$\tilde{A}=\partial_{A}$ and
\[
\int_{\mathbb{R}^{n}}\tilde{A}f\left(  x\right)  p_{T}\left(  x\right)
dx=-\int_{\mathbb{R}^{n}}f\left(  x\right)  \partial_{A}p_{T}\left(  x\right)
dx=\int_{\mathbb{R}^{n}}f\left(  x\right)  \frac{x\cdot A}{T}p_{T}\left(
x\right)  dx
\]
from which it follows that $W_{A}^{T}\left(  x\right)  =\frac{x\cdot A}{T}.$
By simple Gaussian integrations,%
\[
\int_{\mathbb{R}^{n}}e^{W_{A}^{T}\left(  x\right)  }p_{T}\left(  x\right)
dx=\exp\left(  \frac{\left\vert A\right\vert _{\mathfrak{g}}^{2}}{2T}\right)
,
\]%
\begin{align*}
\left(  \int_{\mathbb{R}^{n}}\left[  \frac{p_{T}\left(  x-y\right)  }%
{p_{T}\left(  x\right)  }\right]  ^{p}p_{T}\left(  x\right)  dx\right)
^{1/p}  &  =\left(  \int_{\mathbb{R}^{n}}\left[  e^{-\frac{1}{2T}\left\vert
y\right\vert ^{2}+\frac{1}{T}x\cdot y}\right]  ^{p}p_{T}\left(  x\right)
dx\right)  ^{1/p}\\
&  =e^{-\frac{\left(  p-1\right)  }{2T}\left\vert y\right\vert ^{2}}%
=\exp\left(  c\left(  0\right)  \frac{\left(  p-1\right)  }{2T}\left\vert
y\right\vert ^{2}\right)  ,
\end{align*}
and%
\begin{equation}
\int_{\mathbb{R}^{n}}\left\vert W_{A}^{T}\left(  x\right)  \right\vert
^{p}p_{T}\left(  x\right)  dx=\int_{\mathbb{R}^{n}}\left\vert \frac{x\cdot
A}{T}\right\vert ^{p}p_{T}\left(  x\right)  dx=T^{p/2}\left\vert A\right\vert
^{p}\tilde{C}_{p}^{p}, \label{e.1.13}%
\end{equation}
where
\[
\tilde{C}_{p}^{p}:=\int_{\mathbb{R}^{n}}\left\vert x\right\vert ^{p}%
p_{1}\left(  x\right)  dx.
\]
The first two results show the estimates in Eqs. (\ref{e.1.4}) and
(\ref{e.1.6}) are sharp. The identity in Eq. (\ref{e.1.13}) shows the form of
Eq. (\ref{e.1.5}) is sharp. We do not know if, in general, the constant
$C_{p}$ appearing in Eq. \ref{e.1.5} can be taken to be $\tilde{C}_{p}$
defined above.
\end{example}

Our main interest in Theorem \ref{t.1.6} is in its application to proving that
certain \textquotedblleft heat kernel measures\textquotedblright\ on infinite
dimensional Lie groups, $G,$ are quasi-invariant under left and right
translations by elements of a certain subgroup, $G_{0}.$ We will postpone our
discussion of this application to Section \ref{s.7}. For now let us give a
couple of finite dimensional applications of Theorems \ref{t.1.6} and
\ref{t.1.5}.

\begin{prop}
\label{p.1.8}Suppose that $T>0,$ $p>1,$ and $f\in L^{p}\left(  \nu_{T}\right)
$ is a harmonic function, i.e. $\Delta f=0.$ Then
\begin{equation}
\int_{G}p_{T}\left(  y,x\right)  f\left(  x\right)  dx=f\left(  y\right)
\text{ }\text{ for all }~y\in G. \label{e.1.14}%
\end{equation}

\end{prop}

At an informal level we expect
\[
\int_{G}p_{t}\left(  y,x\right)  f\left(  x\right)  dx=\left(  e^{t\bar
{\Delta}_{0}/2}f\right)  \left(  y\right)
\]
and hence%
\[
\frac{d}{dt}\int_{G}p_{t}\left(  y,x\right)  f\left(  x\right)  dx=\frac
{d}{dt}\left(  e^{t\bar{\Delta}_{0}/2}f\right)  \left(  y\right)  =\left(
e^{t\bar{\Delta}_{0}/2}\frac{\bar{\Delta}_{0}}{2}f\right)  \left(  y\right)
=0.
\]
Therefore it is reasonable to conclude that%
\[
\int_{G}p_{T}\left(  y,x\right)  f\left(  x\right)  dx=\left(  e^{T\bar
{\Delta}_{0}/2}f\right)  \left(  y\right)  =\left(  e^{0\bar{\Delta}_{0}%
/2}f\right)  \left(  y\right)  =f\left(  y\right)  .
\]
However, this argument is not rigorous as $f$ is only square--integrable
relative to the rapidly decaying measure, $\nu_{T},$ rather than to Haar
measure on $G.$ The rigorous proof of Proposition \ref{p.1.8} will be given in
Section \ref{s.7}.

The following corollary is a simple consequence of Proposition \ref{p.1.8},
Eq. (\ref{e.7.4}) in the proof of this proposition, and Theorem \ref{t.1.6} in
the form of Eq. (\ref{e.1.11}).

\begin{cor}
\label{c.1.9}Suppose that $p\in\left(  1,\infty\right)  .$ Under the
hypothesis of Theorem \ref{t.1.6}, if $f\in L^{p}\left(  \nu_{T}\right)  $ and
$f$ is harmonic (i.e. $\Delta f=0),$ then
\begin{equation}
\left\vert f\left(  y\right)  \right\vert \leq\left\Vert f\right\Vert
_{L^{p}\left(  \nu_{T}\right)  }\exp\left(  \frac{c\left(  kT\right)
}{2T\left(  p-1\right)  }\left\vert y\right\vert ^{2}\right)  . \label{e.1.15}%
\end{equation}
In particular, if $G$ is further assumed to be a complex Lie group and $f\in
L^{p}\left(  \nu_{T}\right)  $ is assumed to be holomorphic, then the
pointwise bound in Eq. (\ref{e.1.15}) is still valid.
\end{cor}

\begin{rem}
\label{r.1.10}When $f$ is holomorphic, $p=2,$ $T=1/2,$ and $G=\mathbb{C}^{d},$
the inequality in Eq. (\ref{e.1.15}) is Bargmann's pointwise bound in
\cite[(Eq. (1.7)]{Bargmann61} except that the constant in the exponent is off
by a factor of two. More generally, when $G$ is a general complex Lie group
and $f$ is holomorphic, it has been shown in \cite[Corollary 5.4]{Driver1997c}
that
\[
\left\vert f(y)\right\vert \leq\left\Vert f\right\Vert _{L^{2}\left(
\nu_{t/2}\right)  }e^{|y|^{2}/2t}\text{ for all }y\in G.
\]
The reason for the discrepancy in the coefficients in the exponents between
these inequalities is that $p_{t/2}\left(  x,y\right)  $ is not the
reproducing kernel for the holomorphic functions in $L^{2}\left(  \nu
_{t/2}\right)  $ in that $y\rightarrow p_{t/2}\left(  x,y\right)  $ is not
holomorphic. The coefficient in the exponent of Eq. (\ref{e.1.15}) is also not
sharp since $y\rightarrow p_{T}\left(  x,y\right)  $ is not harmonic.
\end{rem}

\section{$L^{p}$ -- Jacobian estimates\label{s.2}}

Let $M$ be a finite dimensional manifold, $\mu$ be a probability measure on
$M$ with a smooth, strictly positive density in each coordinate chart. For
$r>0,$ let $\left\Vert f\right\Vert _{r}:=\left(  \int_{M}\left\vert
f\right\vert ^{r}d\mu\right)  ^{1/r}$ denote the $L^{r}\left(  \mu\right)  $
-- norm of $f:M\rightarrow\mathbb{C}.$

Let $X_{t}$ be a time dependent vector field and let $S_{t}$ denote its flow,
i.e. $S_{t}\left(  m\right)  $ solves,
\begin{equation}
\frac{d}{dt}S_{t}\left(  m\right)  =X_{t}\circ S_{t}\left(  m\right)  \text{
with }S_{0}\left(  m\right)  =m\text{ for all }m\in M. \label{e.2.1}%
\end{equation}
We will assume that $X_{t}$ is forward complete, i.e. $S_{t}\left(  m\right)
$ exists for all $t\geq0$ and $m\in M.$ Define%
\[
\mu_{t}=\left(  S_{t}\right)  _{\ast}\mu=\mu\circ S_{t}^{-1}.
\]
Since $\mu_{t}$ also has a strictly positive density in each coordinate chart
the Radon-Nikodym derivative
\[
J_{t}=d\mu_{t}/d\mu
\]
exists for all $t\geq0.$ Our goal of this section is to prove Theorem
\ref{t.2.5} below which gives an upper bound on $\left\Vert J_{t}\right\Vert
_{p}$ for $p\in\left(  1,\infty\right)  .$ This result is a slight extension
of the part of Theorem 2.14 in Galaz-Fontes, Gross, and Sontz \cite{GGS01} to
the setting of time dependent vector fields, $X_{t}.$ For the readers
convenience we will sketch the method introduced in \cite[Theorem 2.14]%
{GGS01}. In what follows, $0\ln0$ is to always be interpreted to be $0$.

\begin{lem}
\label{l.2.1}Suppose that $\left(  t,m\right)  \in\left(  0,T\right)  \times
M\rightarrow h_{t}\left(  m\right)  \in\lbrack0,\infty)$ is a smooth bounded
function and $r:\left(  0,T\right)  \rightarrow(1,\infty)$ is a $C^{1}$ --
function. Then%
\begin{equation}
\frac{d}{dt}\ln\left\Vert h_{t}\right\Vert _{r\left(  t\right)  }=\frac
{\dot{r}\left(  t\right)  }{r\left(  t\right)  }\int_{M}\frac{h_{t}^{r\left(
t\right)  }}{\left\Vert h_{t}\right\Vert _{r\left(  t\right)  }^{r\left(
t\right)  }}\left(  \ln\frac{h_{t}}{\left\Vert h_{t}\right\Vert _{r\left(
t\right)  }}\right)  ~d\mu+\frac{1}{r\left(  t\right)  }\int_{M}\frac{\frac
{d}{ds}|_{s=t}h_{s}^{r\left(  t\right)  }}{\left\Vert h_{t}\right\Vert
_{r\left(  t\right)  }^{r\left(  t\right)  }}d\mu. \label{e.2.2}%
\end{equation}

\end{lem}

\begin{proof}
For the reader's convenience we will give a formal derivation of this identity
and refer the reader to Gross \cite[Lemma 1.1]{Gross75} for the technical
details. For $r>0$ and any bounded measurable function, $g:M\rightarrow
\mathbb{R},$ a straight forward calculation shows%
\[
\frac{d}{dr}\ln\left\Vert g\right\Vert _{r}=\frac{1}{r}\int_{M}\frac
{\left\vert g\right\vert ^{r}}{\left\Vert g\right\Vert _{r}^{r}}\left(
\ln\frac{\left\vert g\right\vert }{\left\Vert g\right\Vert _{r}}\right)
~d\mu.
\]
If we further assume that $r>1$ and $v:M\rightarrow\mathbb{R}$ is another
bounded measurable function, then%
\begin{align*}
\partial_{v}\ln\left\Vert g\right\Vert _{r}  &  =\partial_{v}\left[  \frac
{1}{r}\ln\left(  \int_{M}\left\vert g\right\vert ^{r}d\mu\right)  \right]
=\frac{1}{r}\frac{\int_{M}\partial_{v}\left\vert g\right\vert ^{r}d\mu}%
{\int_{M}\left\vert g\right\vert ^{r}d\mu}\\
&  =\frac{1}{r}\int_{M}\frac{\partial_{v}\left\vert g\right\vert ^{r}%
}{\left\Vert g\right\Vert _{r}^{r}}d\mu=\int_{M}\frac{\left\vert g\right\vert
^{r-1}\mathrm{sgn}(g)}{\left\Vert g\right\Vert _{r}^{r}}v~d\mu.
\end{align*}
These two identities along with the chain rule,%
\[
\frac{d}{dt}\ln\left\Vert h_{t}\right\Vert _{r\left(  t\right)  }=\frac{d}%
{ds}|_{s=t}\left[  \ln\left\Vert h_{t}\right\Vert _{r\left(  s\right)  }%
+\ln\left\Vert h_{s}\right\Vert _{r\left(  t\right)  }\right]  ,
\]
easily give Eq. (\ref{e.2.2}).
\end{proof}

\begin{lem}
\label{l.2.2}Let $W\in L^{1}\left(  \mu\right)  $ and $f\geq0$ be a bounded
measurable function. Then, for all $s>0,$
\begin{equation}
\int_{M}Wfd\mu\leq s\int_{M}f\ln\frac{f}{\mu\left(  f\right)  }d\mu
+s\mathcal{B}\left(  W/s\right)  \int_{M}fd\mu\label{e.2.3}%
\end{equation}
where
\[
\mathcal{B}\left(  W\right)  :=\ln\left(  \mu\left(  e^{W}\right)  \right)
=\ln\left(  \int_{M}e^{W}d\mu\right)  .
\]

\end{lem}

\begin{proof}
Recall that Young's inequality states, $xy\leq e^{x}+y\ln y-y$ for
$x\in\mathbb{R}$ and $y\geq0,$ where $0\ln0:=0.$ Applying Young's inequality
with $x=W$ and $y=f$ and then integrating the result gives%
\[
\int_{M}Wfd\mu\leq\int_{M}e^{W}d\mu+\int_{M}\left[  f\ln f-f\right]  d\mu.
\]
Replacing $f$ by $\lambda f$ with $\lambda>0$ in this inequality then shows%
\begin{align*}
\int_{M}Wfd\mu &  \leq\lambda^{-1}\left[  \int_{M}e^{W}d\mu+\int_{M}\left[
\lambda f\ln\left(  \lambda f\right)  -\lambda f\right]  d\mu\right] \\
&  =\lambda^{-1}\int_{M}e^{W}d\mu+\ln\lambda\int_{M}fd\mu+\int_{M}\left[  f\ln
f-f\right]  d\mu.
\end{align*}
The minimizer of the right side of this inequality occurs at $\lambda=\left(
\int_{M}e^{W}d\mu\right)  \cdot\left(  \int_{M}fd\mu\right)  ^{-1}$ and using
this value for $\lambda$ gives%
\begin{equation}
\int_{M}Wfd\mu\leq\int_{M}f\ln\frac{f}{\mu\left(  f\right)  }d\mu
+\mathcal{B}\left(  W\right)  \int_{M}fd\mu. \label{e.2.4}%
\end{equation}
(The proof of Eq. (\ref{e.2.4}) was predicated on the assumption that
$\mathcal{B}\left(  W\right)  <\infty$ but clearly Eq. (\ref{e.2.4}) remains
valid when $\mathcal{B}\left(  W\right)  =\infty.)$ The estimate in Eq.
(\ref{e.2.3}) follows directly from this by replacing $W$ by $W/s.$
\end{proof}

\begin{df}
\label{d.2.3}The $\mu$\textit{--divergence} of a smooth vector field, $X,$ on
$M$ is the function $W=W_{X}^{\mu}$ defined by
\[
\int_{M}X\varphi d\mu=\int_{M}\varphi Wd\mu,\ \text{for all}\ \varphi\in
C_{c}^{1}(M).
\]

\end{df}

\begin{prop}
\label{p.2.4}Let $X_{t}$ and $S_{t}$ be as in Eq. (\ref{e.2.1}),
$W_{t}:=W_{X_{t}}$ be the $\mu$\textit{--divergence} of $X_{t},$ $h\in
C^{1}\left(  M,[0,\infty)\right)  ,$ $h_{t}:=h\circ S_{t}^{-1},$ and $r\in
C^{1}\left(  \left(  0,\tau\right)  ,(1,\infty)\right)  .$ Then for any $s>0$
we have%
\begin{equation}
\frac{d}{dt}\ln\left\Vert h_{t}\right\Vert _{r\left(  t\right)  }\geq\left(
\frac{\dot{r}}{r}-s\right)  \int_{M}\frac{h_{t}^{r}}{\left\Vert h_{t}%
\right\Vert _{r}^{r}}\left(  \ln\frac{h_{t}}{\left\Vert h_{t}\right\Vert _{r}%
}\right)  ~d\mu-\frac{s}{r}\mathcal{B}\left(  s^{-1}W_{t}\right)  .
\label{e.2.5}%
\end{equation}

\end{prop}

\begin{proof}
Differentiating the identity $S_{t}\circ S_{t}^{-1}\left(  m\right)  =m$ and
making use of the flow Eq. (\ref{e.2.1}) implies%
\[
X_{t}\left(  m\right)  +\left(  S_{t}\right)  _{\ast}\frac{d}{dt}S_{t}%
^{-1}\left(  m\right)  =0.
\]
Therefore,%
\[
\frac{d}{dt}S_{t}^{-1}\left(  m\right)  =-\left(  S_{t}^{-1}\right)  _{\ast
}X_{t}\left(  m\right)
\]
or equivalently,%
\[
\frac{d}{dt}f\left(  S_{t}^{-1}\left(  m\right)  \right)
=-X_{t}\left( f\circ S_{t}^{-1}\right)  \left(  m\right)  \text{ for
all }f\in C^{1}\left( M\right).
\]
Using this identity along with Eq. (\ref{e.2.2}) shows%
\begin{equation}
\frac{d}{dt}\ln\left\Vert h_{t}\right\Vert _{r\left(  t\right)  }=\frac
{\dot{r}}{r}\int_{M}\frac{h_{t}^{r}}{\left\Vert h_{t}\right\Vert _{r}^{r}%
}\left(  \ln\frac{h_{t}}{\left\Vert h_{t}\right\Vert _{r}}\right)  ~d\mu
-\frac{1}{r}\int_{M}\frac{X_{t}h_{t}^{r}}{\left\Vert h_{t}\right\Vert _{r}%
^{r}}d\mu\label{e.2.6}%
\end{equation}
where $r=r\left(  t\right)  $ and $\dot{r}=\dot{r}\left(  t\right)  .$
Combining this identity with the definition of $W_{t}$ and the estimate in Eq.
(\ref{e.2.3}) with $W=W_{t}$ and $f=\frac{h_{t}^{r}}{\left\Vert h_{t}%
\right\Vert _{r}^{r}}$ then implies,%
\begin{align*}
\frac{d}{dt}\ln\left\Vert h_{t}\right\Vert _{r\left(  t\right)  }=  &
\frac{\dot{r}}{r}\int_{M}\frac{h_{t}^{r}}{\left\Vert h_{t}\right\Vert _{r}%
^{r}}\left(  \ln\frac{h_{t}}{\left\Vert h_{t}\right\Vert _{r}}\right)
~d\mu-\frac{1}{r}\int_{M}W_{t}\frac{h_{t}^{r}}{\left\Vert h_{t}\right\Vert
_{r}^{r}}d\mu\\
\geq &  \frac{\dot{r}}{r}\int_{M}\frac{h_{t}^{r}}{\left\Vert h_{t}\right\Vert
_{r}^{r}}\left(  \ln\frac{h_{t}}{\left\Vert h_{t}\right\Vert _{r}}\right)
~d\mu\\
&  \qquad-\frac{s}{r}\left[  \int_{M}\frac{h_{t}^{r}}{\left\Vert
h_{t}\right\Vert _{r}^{r}}\ln\frac{h_{t}^{r}}{\left\Vert h_{t}\right\Vert
_{r}^{r}}d\mu+\mathcal{B}\left(  W_{t}/s\right)  \right]
\end{align*}
which is the same as Eq. (\ref{e.2.5}).
\end{proof}

The following theorem is the extension of Galaz-Fontes, Gross, and Sontz
\cite[Theorem 2.14]{GGS01} from time--independent vector fields to
time--dependent vector fields. These results generalize the fundamental
results of Cruzerio \cite{Cruzeiro83a} -- also see
\cite{Bell85,BogMayer1999,Cruzeiro83b,Driver1997b,Peters95b,Peters96b} for
other related results.

\begin{thm}
[Jacobian Estimate]\label{t.2.5}Let $p>1$ and $r\in C\left(  \left[
0,\tau\right]  ,[1,\infty)\right)  \cap C^{1}\left(  \left(  0,\tau\right)
,\left(  1,\infty\right)  \right)  $ such that $r\left(  0\right)  =1,$
$r\left(  \tau\right)  =p$ and $\dot{r}\left(  t\right)  >0$ for $0<t<\tau,$
then%
\begin{equation}
\left\Vert J_{\tau}\right\Vert _{p^{\prime}}\leq e^{\Lambda\left(  r\right)
}, \label{e.2.7}%
\end{equation}
where $p^{\prime}:=p/\left(  p-1\right)  $ is the conjugate exponent to $p$
and
\begin{equation}
\Lambda\left(  r\right)  =\Lambda_{X}\left(  r\right)  :=\int_{0}^{\tau}%
\frac{\dot{r}\left(  t\right)  }{r^{2}\left(  t\right)  }\mathcal{B}\left(
\frac{r\left(  t\right)  }{\dot{r}\left(  t\right)  }W_{t}\right)  dt.
\label{e.2.8}%
\end{equation}

\end{thm}

\begin{proof}
Taking $s=\dot{r}/r$ in Eq. (\ref{e.2.5}) gives%
\[
\frac{d}{dt}\ln\left\Vert h_{t}\right\Vert _{r\left(  t\right)  }\geq
-\frac{\dot{r}}{r^{2}}\mathcal{B}\left(  \frac{r}{\dot{r}}W_{t}\right)
\]
which integrates to
\[
\left\Vert h\circ S_{\tau}^{-1}\right\Vert _{p}=\left\Vert h_{\tau}\right\Vert
_{p}\geq\left\Vert h\right\Vert _{1}\exp\left(  -\int_{0}^{\tau}\frac{\dot
{r}\left(  t\right)  }{r^{2}\left(  t\right)  }\mathcal{B}\left(
\frac{r\left(  t\right)  }{\dot{r}\left(  t\right)  }W_{t}\right)  dt\right)
.
\]
Replacing $h$ by $h\circ S_{\tau}$ in this inequality implies%
\begin{equation}
\int_{M}hJ_{\tau}d\mu=\left\Vert h\circ S_{\tau}\right\Vert _{1}\leq\left\Vert
h\right\Vert _{p}e^{\Lambda\left(  r\right)  }. \label{e.2.9}%
\end{equation}
Let $L^{p}\left(  \mu\right)  ^{+}$ denote the almost everywhere non-negative
functions in $L^{p}\left(  \mu\right)  .$ Since Eq. (\ref{e.2.9}) is valid for
all $h\in C^{1}\left(  M,[0,\infty)\right)  $ and the latter functions are
dense in $L^{p}\left(  \mu\right)  ^{+}$ (see the proof of Lemma 2.8 in
\cite{GGS01}), it follow that Eq. (\ref{e.2.9}) is valid for all $h\in
L^{p}\left(  \mu\right)  ^{+}.$ Equation \ref{e.2.7}) now follows by the
converse to H\"{o}lder's inequality. Indeed, let $K\subset M$ be a compact set
and take $h=J_{\tau}^{p^{\prime}-1}1_{K}=J_{\tau}^{1/\left(  p-1\right)
}1_{K}$ in Eq. (\ref{e.2.9}) to find%
\[
\int_{M}J_{\tau}^{p^{\prime}}1_{K}d\mu\leq\left\Vert J_{\tau}^{1/\left(
p-1\right)  }1_{K}\right\Vert _{p}e^{\Lambda\left(  r\right)  }=\left(
\int_{M}J_{\tau}^{p^{\prime}}1_{K}d\mu\right)  ^{1/p}e^{\Lambda\left(
r\right)  }.
\]
This inequality is equivalent to%
\[
\left\Vert J_{\tau}1_{K}\right\Vert _{p^{\prime}}=\left(  \int_{M}J_{\tau
}^{p^{\prime}}1_{K}d\mu\right)  ^{1-1/p}\leq e^{\Lambda\left(  r\right)  }.
\]
Now replacing $K$ by $K_{n}$ with $K_{n}$ compact and $K_{n}\uparrow M$ and
passing to the limit as $n\rightarrow\infty$ in the previous inequality gives
the estimate in Eq. (\ref{e.2.7}).
\end{proof}

\section{Proof of Theorem \ref{t.1.6}\label{s.3}}

In this section we will give a proof of Theorem \ref{t.1.6} assuming that
Theorem \ref{t.1.4} holds.

\begin{proof}
(Proof of Theorem \ref{t.1.6}.) In order to abbreviate the notation, let
$c:=c\left(  kT\right)  /T.$ Let $g\in C^{1}\left(  \left[  0,1\right]
,G\right)  $ be such that $g\left(  0\right)  =e\in G$ and $g\left(  1\right)
=y\in G$ and define $A_{t}:=L_{g\left(  t\right)  \ast}^{-1}\dot{g}\left(
t\right)  \in\mathfrak{g.}$ If we now let $X_{t}:=\tilde{A}_{t}\in
\Gamma\left(  TG\right)  ,$ then the flow, $S_{t},$ of $X_{t}$ satisfies,
$S_{t}\left(  x\right)  =xg\left(  t\right)  .$ Indeed, because $X_{t}$ is
left invariant,%
\[
\frac{d}{dt}xg\left(  t\right)  =L_{x\ast}\dot{g}\left(  t\right)  =L_{x\ast
}L_{g\left(  t\right)  \ast}A_{t}=L_{xg\left(  t\right)  \ast}A_{t}%
=X_{t}\left(  xg\left(  t\right)  \right)  .
\]

In order to apply the Jacobian estimate in Theorem \ref{t.2.5}, let
$d\mu\left(  x\right)  =d\nu_{T}\left(  x\right)  :=p_{T}\left(  x\right)  dx$
and observe that
\begin{align*}
\int_{G}h\left(  S_{1}\left(  x\right)  \right)  d\mu\left(  x\right)   &
=\int_{G}h\left(  xy\right)  d\mu\left(  x\right)  =\int_{G}h\left(
xy\right)  p_{T}\left(  x\right)  dx\\
&  =\int_{G}h\left(  x\right)  p_{T}\left(  xy^{-1}\right)  dx=\int
_{G}h\left(  x\right)  \frac{p_{T}\left(  xy^{-1}\right)  }{p_{T}\left(
x\right)  }d\mu\left(  x\right)
\end{align*}
from which it follows that
\begin{equation}
J_{1}\left(  x\right)  :=\frac{d\left(  S_{1}\right)  _{\ast}\mu}{d\mu}\left(
x\right)  =\frac{p_{T}\left(  xy^{-1}\right)  }{p_{T}\left(  x\right)  }.
\label{e.3.1}%
\end{equation}
Moreover, if $W_{t}=W_{X_{t}}^{\nu_{T}}$ is the $\mu=\nu_{T}$ -- divergence of
$X_{t},$ by Theorem \ref{t.1.4},
\begin{equation}
\mathcal{B}\left(  \lambda W_{t}\right)  =\ln\left(  \int_{G}e^{\lambda W_{t}%
}d\mu\right)  \leq\frac{c\left(  kT\right)  }{T}\lambda^{2}\left\vert
A_{t}\right\vert _{\mathfrak{g}}^{2}. \label{e.3.2}%
\end{equation}
Hence it follows from Theorem \ref{t.2.5} that%
\begin{equation}
\left[  \int_{G}\left(  \frac{p_{T}\left(  xy^{-1}\right)  }{p_{T}\left(
x\right)  }\right)  ^{p^{\prime}}p_{T}\left(  x\right)  dx\right]
^{1/p^{\prime}}=\left\Vert J_{1}\right\Vert _{p^{\prime}}\leq e^{\Lambda
\left(  r\right)  }, \label{e.3.3}%
\end{equation}
where%
\begin{align*}
\Lambda\left(  r\right)   &  =\int_{0}^{1}\frac{\dot{r}\left(  t\right)
}{r^{2}\left(  t\right)  }\mathcal{B}\left(  \frac{r\left(  t\right)  }%
{\dot{r}\left(  t\right)  }W_{t}\right)  dt\\
&  \leq c\int_{0}^{1}\frac{\dot{r}\left(  t\right)  }{r^{2}\left(  t\right)
}\frac{r^{2}\left(  t\right)  }{\dot{r}^{2}\left(  t\right)  }\left\vert
A_{t}\right\vert _{\mathfrak{g}}^{2}dt=c\int_{0}^{1}\frac{\left\vert
A_{t}\right\vert _{\mathfrak{g}}^{2}}{\dot{r}\left(  t\right)  }dt,
\end{align*}
and $r\in C\left(  \left[  0,1\right]  ,[1,\infty)\right)  \cap C^{1}\left(
\left(  0,1\right)  ,\left(  1,\infty\right)  \right)  $ such that $r\left(
0\right)  =1,$ $r\left(  1\right)  =p$ and $\dot{r}\left(  t\right)  >0$ for
$0<t<1.$

We now want to choose $r\left(  t\right)  $ so as to minimize $\Lambda\left(
r\right)  $ subject to the constraints $\dot{r}\left(  t\right)  >0,$
$r\left(  0\right)  =1$ and $r\left(  1\right)  =p.$ To see how to choose $r,$
let us differentiate $\Lambda\left(  r\right)  $ in a direction $v$ such that
$v\left(  0\right)  =0=v\left(  1\right)  $ and then require%
\[
0\overset{\text{set}}{=}\left(  \partial_{v}\Lambda\right)  \left(  r\right)
=-\frac{c}{2}\int_{0}^{1}\frac{\left\vert A_{t}\right\vert _{\mathfrak{g}}%
^{2}}{\dot{r}^{2}\left(  t\right)  }\dot{v}\left(  t\right)  dt=-\frac{c}%
{2}\int_{0}^{1}v\left(  t\right)  \frac{d}{dt}\left(  \frac{\left\vert
A_{t}\right\vert _{\mathfrak{g}}^{2}}{\dot{r}^{2}\left(  t\right)  }\right)
dt.
\]
Since $v\left(  t\right)  $ is arbitrary, we should require $\frac{\left\vert
A_{t}\right\vert _{\mathfrak{g}}^{2}}{\dot{r}^{2}\left(  t\right)  }%
=\kappa^{-2},$ where $\kappa>0$ is a constant, i.e. $\dot{r}\left(  t\right)
=\kappa\left\vert A_{t}\right\vert _{\mathfrak{g}}.$ Hence we take
\[
r\left(  t\right)  =1+\kappa\int_{0}^{t}\left\vert A_{\tau}\right\vert
_{\mathfrak{g}}d\tau,
\]
where%
\[
\kappa:=\left(  p-1\right)  \left(  \int_{0}^{1}\left\vert A_{\tau}\right\vert
_{\mathfrak{g}}d\tau\right)  ^{-1}%
\]
has been chosen so that $r\left(  1\right)  =p.$ With this choice of $r,$%
\[
\Lambda\left(  r\right)  :=\frac{c}{2}\int_{0}^{1}\frac{\left\vert
A_{t}\right\vert _{\mathfrak{g}}^{2}}{\kappa\left\vert A_{t}\right\vert
_{\mathfrak{g}}}dt=\frac{c}{2\kappa}\int_{0}^{1}\left\vert A_{t}\right\vert
_{\mathfrak{g}}dt=\frac{c}{2\left(  p-1\right)  }\left(  \int_{0}%
^{1}\left\vert A_{t}\right\vert _{\mathfrak{g}}dt\right)  ^{2}%
\]
and using this value for $\Lambda\left(  r\right)  $ in Eq. (\ref{e.3.3})
along with the identity, $\left(  p-1\right)  ^{-1}=p^{\prime}-1$ implies%
\[
\left(  \int_{G}\left[  \frac{p_{T}\left(  xy^{-1}\right)  }{p_{T}\left(
x\right)  }\right]  ^{p^{\prime}}p_{T}\left(  x\right)  dx\right)
^{1/p^{\prime}}=\left\Vert J_{1}\right\Vert _{p^{\prime}}\leq\exp\left(
\frac{c\left(  p^{\prime}-1\right)  }{2}\left(  \int_{0}^{1}\left\vert
A_{t}\right\vert _{\mathfrak{g}}dt\right)  ^{2}\right)  .
\]
Upon noting that $p^{\prime}:=p\left(  p-1\right)  ^{-1}$ ranges over $\left(
1,\infty\right)  $ as $p$ ranges over $\left(  1,\infty\right)  ,$ the proof
of Theorem \ref{t.1.6} is complete.
\end{proof}

\section{Properties of the Hodge -- de Rham semigroups\label{s.4}}

This section gathers a number of technical functional analytic results needed
to establish the representation formula in Theorem \ref{t.5.4} below. Let
$\left(  M,g\right)  $ be a complete Riemannian manifold, $dV$\ denote the
volume measure on $M$ associated to $g,$ $\nabla$ denote the Levi-Civita
covariant derivative, $\Lambda^{k}=\Lambda^{k}\left(  T^{\ast}M\right)  ,$
$\Lambda=\oplus_{k=0}^{\dim M}\Lambda^{k},$ $\Omega^{k}\left(  M\right)  $
($\Omega_{c}^{k}\left(  M\right)  $) denote the space of (compactly supported)
smooth $k$ -- forms over $M,$ and $\Omega\left(  M\right)  =\oplus_{k=0}^{\dim
M}\Omega^{k}\left(  M\right)  $ be the space of all smooth forms over $M.$ If
$\alpha$ and $\beta$ are measurable $k$ -- forms, let%
\[
\left\langle \alpha,\beta\right\rangle _{m}:=\sum_{j_{1},\dots,j_{k}=1}%
^{d}\alpha\left(  e_{j_{1}},\dots,e_{j_{k}}\right)  \beta\left(  e_{j_{1}%
},\dots,e_{j_{k}}\right)  ,
\]
where $\left\{  e_{j}\right\}  _{j=1}^{d}$ is any orthonormal frame for
$T_{m}M.$ When $m\rightarrow\left\langle \alpha,\beta\right\rangle _{m}$ is
integrable, let
\[
\left(  \alpha,\beta\right)  :=\int_{M}\left\langle \alpha,\beta\right\rangle
dV
\]
and let $L^{2}\left(  \Lambda^{k}\right)  $ denote the measurable $k$ --
forms, $\alpha,$ such that $\left(  \alpha,\alpha\right)  <\infty.$ Further
let
\[
L^{2}\left(  \Lambda\right)  :=\oplus_{k=0}^{\dim M}L^{2}\left(  \Lambda
^{k}\right)  .
\]
Two measurable $k$ -- forms, $\alpha$ and $\beta,$ are take to be equivalent
if $\alpha=\beta$ a.e.

Let $d:\Omega\left(  M\right)  \rightarrow\Omega\left(  M\right)  $ be the
differential operator taking $k$ -- forms to $k+1$ -- forms, $\delta$ be the
formal $L^{2}$ -- adjoint of $-d$,
\[
\Delta:=-\left(  \delta d+d\delta\right)  =-\left(  d+\delta\right)  ^{2}%
\]
be the Hodge-de Rham Laplacian on $\Omega\left(  M\right)  ,$ and $\square$ be
the Bochner (i.e. flat) Laplacian on $\Omega\left(  M\right)  .$ More
precisely if $\alpha$ is a $k$ -- form, $\delta\alpha$ is the $k-1$ form
defined by%
\begin{equation}
\left(  \delta\alpha\right)  _{m}:=\sum_{j=1}^{d}\left(  \nabla_{e_{j}}%
\alpha\right)  \left(  e_{j},\text{--}\right)  \label{e.4.1}%
\end{equation}
and%
\[
\left(  \square\alpha\right)  _{m}:=\sum_{j=1}^{d}\nabla_{e_{j}\otimes e_{j}%
}^{2}\alpha:=\sum_{j=1}^{d}\left(  \nabla_{E_{j}}^{2}\alpha-\nabla
_{\nabla_{E_{j}}E_{j}}\alpha\right)  _{m}%
\]
where $\left\{  E_{j}\right\}  _{j=1}^{\dim M}$ is an local orthonormal frame
for $TM$ defined in a neighborhood of $m.$ The next two theorems summarize the
properties about these operators that will be needed in this paper.

\begin{thm}
\label{t.4.1}The operators, $d_{k}:=d|_{\Omega_{c}^{k}\left(  M\right)
}:\Omega_{c}^{k}\left(  M\right)  \rightarrow\Omega_{c}^{k+1}\left(  M\right)
$ for $k=0,1,2\dots,\dim M-1$ are $L^{2}\left(  \Lambda^{k}\right)  $ --
closable with closure denoted by $\bar{d}_{k}.$ Let us now further assume that
$\left(  M,g\right)  $ is complete. Then:

\begin{enumerate}
\item Each of the operators, $\Delta_{k}:=\Delta|_{\Omega_{c}^{k}\left(
M\right)  }$ for $k=0,1,2\dots,\dim M$ thought of as unbounded operators on
$L^{2}\left(  \Lambda^{k}\right)  ,$ are essentially self-adjoint operators.
Let $\bar{\Delta}_{k}$ denote the (self-adjoint) closure of $\Delta_{k}.$

\item Each operator, $\bar{\Delta}_{k},$ is non-negative. Let $e^{t\bar
{\Delta}_{k}}$ denotes the contraction semi-group on $L^{2}\left(  \Lambda
^{k}\right)  $ associated to $\bar{\Delta}_{k}.$

\item For $k\in\left\{  0,1,\dots,\dim M-1\right\}  $ and $t>0,$ $\bar{d}%
_{k}e^{t\bar{\Delta}_{k}}=e^{t\bar{\Delta}_{k+1}}\bar{d}_{k}$ on the domain of
$\bar{d}_{k}.$

\item $\delta_{k}e^{t\bar{\Delta}_{k}}\omega=e^{t\bar{\Delta}_{k-1}}\delta
_{k}\omega$ for all $\omega\in\Omega_{c}^{k}\left(  M\right)  $ with
$k=1,2,\dots,\dim M.$
\end{enumerate}
\end{thm}

\begin{proof}
Let $\delta_{k}:=\delta|_{\Omega_{c}^{k}\left(  M\right)  }:\Omega_{c}%
^{k}\left(  M\right)  \rightarrow\Omega_{c}^{k-1}\left(  M\right)  .$ As
$-\delta_{k+1}\subset d_{k}^{\ast},$ $d_{k}^{\ast}$ is densely defined and
hence $d_{k}$ is closable. For items 1. and 2., see Gaffney \cite{Gaffney51},
Roelcke \cite{Roelcke60}, Chernoff \cite{Chernoff73}, \cite{Yau78}, and
Strichartz \cite{Str1983}.

Item 3. is a simple application of Theorem \ref{t.A.2} of Appendix \ref{s.A}
below. In applying this theorem, take $W=L^{2}\left(  \Lambda^{k-1}\right)  ,$
$X=L^{2}\left(  \Lambda^{k}\right)  ,$ $Y=L^{2}\left(  \Lambda^{k+1}\right)  $
and $Z=L^{2}\left(  \Lambda^{k+2}\right)  $ with $A=\bar{d}_{k-1},$ $B=\bar
{d}_{k},$ and $C:=\bar{d}_{k+1}.$ By convention $\Omega^{-1}\left(  M\right)
=\left\{  0\right\}  =\Omega^{\dim M+1}\left(  M\right)  $ and $d_{-1}%
=0=d_{\dim M}.$ With these assignments, the self-adjoint operators, $L$ and
$S,$ in Theorem \ref{t.A.2} become%
\begin{equation}
L=\bar{d}_{k-1}d_{k-1}^{\ast}+d_{k}^{\ast}\bar{d}_{k}\text{ and }S=\bar{d}%
_{k}d_{k}^{\ast}+d_{k+1}^{\ast}\bar{d}_{k+1}. \label{e.4.2}%
\end{equation}
As $\Delta_{k}|_{\Omega_{c}^{k}\left(  M\right)  }\subset-L$ and $-L$ is
self-adjoint (see Theorem \ref{t.A.1} below), it follows that $\bar{\Delta
}_{k}=-L$ and similarly, $\bar{\Delta}_{k+1}=-S.$

For item 4., let $\omega\in\Omega_{c}^{k}\left(  M\right)  $ and $\varphi
\in\Omega_{c}^{k-1}\left(  M\right)  .$ Then%
\begin{align*}
\left(  \delta e^{t\bar{\Delta}_{k}}\omega,\varphi\right)   &  =-\left(
e^{t\bar{\Delta}_{k}}\omega,d\varphi\right)  =-\left(  \omega,e^{t\bar{\Delta
}_{k}}\bar{d}\varphi\right)  =-\left(  \omega,\bar{d}e^{t\bar{\Delta}_{k-1}%
}\varphi\right) \\
&  =\left(  \delta\omega,e^{t\bar{\Delta}_{k-1}}\varphi\right)  =\left(
e^{t\bar{\Delta}_{k-1}}\delta\omega,\varphi\right)  .
\end{align*}

\end{proof}

\begin{rem}
\label{r.4.2}With a little more work it is possible to show that $\bar{d}%
_{k}=-\delta_{k+1}^{\ast}$ and that $\bar{\delta}_{k}e^{t\bar{\Delta}_{k}%
}=e^{t\bar{\Delta}_{k-1}}\bar{\delta}_{k}$ on the domain of $\bar{\delta}%
_{k}.$ We will omit the proof of these results as they are not needed for this paper.
\end{rem}

We are primarily concerned with zero and one forms. A key ingredient in the
sequel is the Bochner identity,
\begin{equation}
\Delta\alpha=\square\alpha-\alpha\circ\operatorname{Ric}\text{ for all }%
\alpha\in\Omega^{1}\left(  M\right)  . \label{e.4.3}%
\end{equation}

\begin{assumption}
\label{ass.1}For the rest of this paper we will assume that $\left(
M,g\right)  $ is a complete Riemannian manifold such that $\operatorname{Ric}%
\geq k$ for some $k\in\mathbb{R},$ i.e. $\operatorname{Ric}_{m}\geq
kI_{T_{m}M}$ for all $m\in M.$
\end{assumption}

\begin{thm}
[Semi-group domination]\label{t.4.3}Suppose that $\left(  M,g\right)  $ is a
complete Riemannian manifold such that $\operatorname{Ric}\geq k$ for some
$k\in\mathbb{R}.$ Then for all $f\in L^{2}\left(  \Lambda^{0}\right)  $ and
$\alpha\in L^{2}\left(  \Lambda^{1}\right)  ,$%
\begin{equation}
\left\vert e^{t\bar{\Delta}_{0}}f\right\vert \leq e^{t\bar{\Delta}_{0}%
}\left\vert f\right\vert \leq\left\Vert f\right\Vert _{\infty}\text{ a.e.}
\label{e.4.4}%
\end{equation}
and
\begin{equation}
\left\vert e^{t\bar{\Delta}_{1}}\alpha\right\vert \leq e^{-kt}e^{t\bar{\Delta
}_{0}}\left\vert \alpha\right\vert \leq e^{-kt}\left\Vert \alpha\right\Vert
_{\infty}\text{ a.e.} \label{e.4.5}%
\end{equation}
where $\left\Vert f\right\Vert _{\infty}$ and $\left\Vert \alpha\right\Vert
_{\infty}$ denote the essential supremums of the functions, $\left\vert
f\right\vert $ and $m\rightarrow\left\vert \alpha_{m}\right\vert $ respectively.
\end{thm}

\begin{proof}
The inequality in Eq. (\ref{e.4.4}) is an immediate consequence Eqs.
(\ref{e.1.2}), (\ref{e.1.1}) and the positivity of the heat kernel,
$p_{t}\left(  x,y\right)  .$ This inequality may also be proved using the
semi-group domination ideas that will be used below to prove Eq. (\ref{e.4.5}).

The proof of Eq. (\ref{e.4.5}) will be an application of the results in Simon
\cite{Simon1977,Simon1979} and Hess, Schrader, and Uhlenbrock \cite{HSU1977}
along with a Kato \cite{Kato1972} type inequality. The general Kato inequality
we need is given in Theorem \ref{t.B.2} of Appendix \ref{s.B}. We apply
Theorem \ref{t.B.2} with $E=\Lambda^{1}\left(  T^{\ast}M\right)  $ to
conclude,%
\begin{equation}
\left(  \square\alpha,\varphi\,\mathrm{sgn}_{e}\left(  \alpha\right)  \right)
\leq\left(  \left\vert \alpha\right\vert ,\Delta\varphi\right)  \label{e.4.6}%
\end{equation}
for all $\alpha\in\Omega_{c}^{1}\left(  M\right)  $ and $\varphi\in C^{\infty
}\left(  M\right)  _{+}:=C^{\infty}\left(  M\rightarrow\lbrack0,\infty
)\right)  .$ In Eq. (\ref{e.4.6}),%
\[
\mathrm{sgn}_{e}\left(  \alpha\right)  :=1_{\alpha\neq0}\frac{\alpha
}{\left\vert \alpha\right\vert }+1_{\alpha=0}e,
\]
where $e$ is any measurable section of $E$ such that $\left\langle
\square\alpha,e\right\rangle =0$ on $M.$ This inequality and the Bochner
identity in Eq. (\ref{e.4.3}) shows%
\begin{align}
\left(  \Delta_{1}\alpha,\varphi\,\mathrm{sgn}_{e}\left(  \alpha\right)
\right)   &  =\left(  \square\alpha,\varphi\,\mathrm{sgn}_{e}\left(
\alpha\right)  \right)  -\left(  \alpha\circ\operatorname{Ric},\varphi
\,\mathrm{sgn}_{e}\left(  \alpha\right)  \right) \nonumber\\
&  \leq\left(  \left\vert \alpha\right\vert ,\Delta\varphi\right)  -\left(
\alpha\circ\operatorname{Ric},\varphi\,\mathrm{sgn}_{e}\left(  \alpha\right)
\right)  . \label{e.4.7}%
\end{align}
To evaluate the last term, let $Y$ be the vector field on $M$ such that
$\alpha=\left\langle Y,\cdot\right\rangle .$ Then $\alpha\circ
\operatorname{Ric}=\left\langle \operatorname{Ric}Y,\cdot\right\rangle $ and%
\begin{align*}
\left\langle \alpha\circ\operatorname{Ric},\mathrm{sgn}_{e}\left(
\alpha\right)  \right\rangle  &  =1_{\alpha\neq0}\frac{1}{\left\vert
\alpha\right\vert }\left\langle \alpha\circ\operatorname{Ric},\alpha
\right\rangle =1_{\alpha\neq0}\frac{1}{\left\vert \alpha\right\vert
}\left\langle \operatorname{Ric}Y,Y\right\rangle \\
&  \geq k1_{\alpha\neq0}\frac{1}{\left\vert \alpha\right\vert }\left\langle
Y,Y\right\rangle =k1_{\alpha\neq0}\frac{1}{\left\vert \alpha\right\vert
}\left\vert \alpha\right\vert ^{2}=k\left\vert \alpha\right\vert .
\end{align*}
Therefore,%
\[
\left(  \alpha\circ\operatorname{Ric},\varphi\,\mathrm{sgn}_{e}\left(
\alpha\right)  \right)  =\int_{M}\left\langle \alpha\circ\operatorname{Ric}%
,\mathrm{sgn}_{e}\left(  \alpha\right)  \right\rangle \varphi dV\geq k\left(
\left\vert \alpha\right\vert ,\varphi\right)
\]
which combined with Eq. (\ref{e.4.7}) implies%
\begin{equation}
\left(  \Delta_{1}\alpha,\varphi\,\mathrm{sgn}_{e}\left(  \alpha\right)
\right)  \leq\left(  \left\vert \alpha\right\vert ,\Delta\varphi\right)
-k\left(  \left\vert \alpha\right\vert ,\varphi\right)  \label{e.4.8}%
\end{equation}
or equivalently,
\[
\left(  H_{0}\alpha,\varphi\,\mathrm{sgn}_{e}\left(  \alpha\right)  \right)
\geq\left(  \left\vert \alpha\right\vert ,-\Delta\varphi\right)
\]
where $H_{0}:=-\left(  \Delta+k\right)  |_{\Omega_{c}^{1}\left(  M\right)  }.$
In particular if $g\in C_{c}^{\infty}\left(  M\right)  _{+},$ $\lambda>0,$
$\varphi=\left(  -\bar{\Delta}_{0}+\lambda\right)  ^{-1}g,$ and $\alpha_{1}%
\in\Omega_{c}^{1}\left(  M\right)  $ and we define $\alpha_{2}:=\varphi
\,\mathrm{sgn}_{e}\left(  \alpha_{1}\right)  \in L^{2}\left(  \Lambda
^{1}\right)  ,$ then $\left(  \alpha_{1},\alpha_{2}\right)  _{L^{2}\left(
\Lambda_{1}\right)  }=\left(  \left\vert \alpha_{1}\right\vert ,\left\vert
\alpha_{2}\right\vert \right)  _{L^{2}\left(  \Lambda_{0}\right)  },$
$\left\vert \alpha_{2}\right\vert =\varphi,$ and
\[
\left(  H_{0}\alpha_{1},\alpha_{2}\right)  _{L^{2}\left(  \Lambda_{1}\right)
}\geq\left(  \left\vert \alpha_{1}\right\vert ,-\bar{\Delta}_{0}%
\varphi\right)  _{L^{2}\left(  \Lambda_{0}\right)  }.
\]
Hence we have verified the hypothesis of Proposition 2.14 and Theorem 2.15 in
\cite{HSU1977} and as a consequence,%
\begin{equation}
\left\vert e^{-t\bar{H}_{0}}\alpha\right\vert \leq e^{-t\left(  -\bar{\Delta
}_{0}\right)  }\left\vert \alpha\right\vert \text{ a.e. }\text{ for all
}~\alpha\in L^{2}\left(  \Lambda^{1}\right)  . \label{e.4.9}%
\end{equation}
As $\bar{H}_{0}=-\bar{\Delta}_{1}-k$ and hence, $e^{-t\bar{H}_{0}}%
=e^{t\bar{\Delta}_{1}}e^{tk},$ Eq. (\ref{e.4.9}) is equivalent to the first
inequality in Eq. (\ref{e.4.5}).
\end{proof}

\section{A path integral derivative formula\label{s.5}}

\subsection{Brownian motion and the divergence formula\label{s.5.1}}

Let $\bigl(\Omega,\mathcal{F},\{\mathcal{F}_{t}\}_{t\geq0},\mathbb{P}\bigr)$
be a filtered probability space satisfying the usual hypothesis, and for each
$x\in M$ let $\left\{  \Sigma_{t}^{x}:t<\zeta(x)\right\}  $ be an $M$ --
valued Brownian motion on $\bigl(\Omega,\mathcal{F},\{\mathcal{F}_{t}%
\}_{t\geq0},\mathbb{P}\bigr)$, starting from $x$, with possibly finite
lifetime $\zeta(x)$. Recall $\Sigma_{t}^{x}$ is said to be an $M$--valued
Brownian motion provided it is a Markov diffusion process starting at $x$ with
transition semi-group determined by the heat kernel, $p_{t}\left(  \cdot
,\cdot\right)  .$ Because of our standing assumption, $\operatorname{Ric}\geq
k,$ it is well--known that $\int_{M}p_{t}(x,y)dy=1$ for all $x\in M$ and
consequently that $\zeta\left(  x\right)  =\infty,$ see
\cite{Az,Gaffney59,Yau78,Dodziuk83,Li84,Grig1,Grig2,Grig3} and the books
\cite[Theorem 8.62]{Stroock00}, \cite[Chapter 4.]{Hsu2002} and \cite[Theorem
5.2.6]{Davies89}. For our purposes it will be convenient to construct
$\Sigma_{t}^{x}$ as a solution to a stochastic differential equation which we
will describe shortly.

\begin{nota}
\label{n.5.1}Given two isometric isomorphic real finite--dimensional inner
product spaces, $V$ and $W,$ let $O\left(  V,W\right)  $ denote the set of
linear isometries from $V$ to $W.$
\end{nota}

Let $//_{t}\left(  \sigma\right)  $ denote parallel translation along a curve
$\sigma$ in $TM$ and all associated bundles. We also introduce the horizontal
vector fields on the orthogonal frame bundle over $M$ as
\[
B_{v}\left(  u\right)  =\frac{d}{dt}|_{0}//_{t}\left(  \sigma\right)  u\text{
for }v\in\mathbb{R}^{d}\text{ and }u\in O\left(  \mathbb{R}^{d},T_{x}M\right)
,
\]
where $\sigma\left(  t\right)  $ is a curve in $M$ such that $\dot{\sigma
}\left(  0\right)  =uv.$

\begin{nota}
\label{n.5.2}Given a semi-martingale, $Y_{t},$ we will denote its It\^{o}
differential by $dY_{t}$ and its Fisk-Stratonovich differential by $\circ
dY_{t}.$
\end{nota}

Let $b_{t}$ denote a $\mathbb{R}^{d}$ -- valued Brownian motion, $x\in M,$ and
$u_{0}\in O\left(  \mathbb{R}^{d},T_{x}M\right)  ,$ then $\Sigma_{t}^{x}$ may
be defined as the unique solution to the stochastic differential equation,
\begin{align*}
\circ d\Sigma_{t}^{x}  &  =u_{t}\circ db_{t}\text{ with }\Sigma_{0}^{x}=x,\\
\circ du_{t}  &  =B_{\circ db_{t}}\left(  u_{t}\right)  \text{ with }u_{0}.
\end{align*}
The \emph{stochastic parallel translation along }$\Sigma_{t}^{x}$ up to time
$t$ is taken to be, $//_{t}:=u_{t}u_{0}^{-1}\in O\left(  T_{x}M,T_{\Sigma
_{t}^{x}}M\right)  .$ Suppose that $f\left(  t,m\right)  $ ($\alpha\left(
t,m\right)  )$ is a smooth time dependent function (one form), then the
It\^{o} differentials of $f\left(  t,\Sigma_{t}^{x}\right)  $ and
$\alpha\left(  t,\Sigma_{t}^{x}\right)  //_{t}$ are%
\begin{equation}
d\left[  f\left(  t,\Sigma_{t}^{x}\right)  \right]  =\left(  \frac{\partial
}{\partial t}f\left(  t,\Sigma_{t}^{x}\right)  +\frac{1}{2}\Delta_{0}f\left(
t,\Sigma_{t}^{x}\right)  \right)  dt+\left\langle \operatorname*{grad}f\left(
t,\cdot\right)  ,//_{t}db_{t}\right\rangle \label{e.5.1}%
\end{equation}
and
\begin{equation}
d\left[  \alpha\left(  t,\Sigma_{t}^{x}\right)  //_{t}\right]  =\left(
\frac{\partial}{\partial t}\alpha\left(  t,\Sigma_{t}^{x}\right)  +\frac{1}%
{2}\square\alpha\left(  t,\Sigma_{t}^{x}\right)  \right)  dt+\left[
\nabla_{//_{t}db_{t}}\alpha\left(  t,\cdot\right)  \right]  //_{t}.
\label{e.5.2}%
\end{equation}
See (for example) \cite{Emery89,Mal97,Stroock00,Hsu2002,Driver2004} for more
on the general background used in this section.

\subsection{The divergence formula\label{s.5.2}}

Let $Q_{t}$ denote the $\operatorname*{End}\left(  T_{x}M\right)  $ -- valued
process satisfying the ordinary differential equation,%
\begin{equation}
\frac{d}{dt}Q_{t}=-\frac{1}{2}\operatorname{Ric}^{//_{t}}Q_{t}\quad\text{with
}Q_{0}=id_{T_{x}M}. \label{e.5.3}%
\end{equation}
where
\begin{equation}
\operatorname{Ric}^{//_{t}}:=//_{t}^{-1}\operatorname{Ric}_{\Sigma_{t}^{x}%
}//_{t}. \label{e.5.4}%
\end{equation}

\begin{lem}
\label{l.5.3}If $\operatorname{Ric}\geq k$ for some $k\in\mathbb{R}$ and
$\left\Vert \cdot\right\Vert _{op}$ denotes the operator norm on $T_{x}M,$
then
\begin{equation}
\left\Vert Q_{t}\right\Vert _{op}\leq e^{-kt/2}. \label{e.5.5}%
\end{equation}
Similarly if $\operatorname{Ric}\leq K$ for some $K\in\mathbb{R},$ then%
\begin{equation}
\left\Vert Q_{t}^{-1}\right\Vert _{op}\leq e^{Kt/2}. \label{e.5.6}%
\end{equation}

\end{lem}

\begin{proof}
For any $v\in T_{x}M,$ we have%
\[
\frac{d}{dt}\left\vert Q_{t}v\right\vert ^{2}=\left\langle -\operatorname{Ric}%
^{//_{t}}Q_{t}v,Q_{t}v\right\rangle \leq-k\left\vert Q_{t}v\right\vert ^{2}%
\]
from which Eq. (\ref{e.5.5}) easily follows. To prove Eq. (\ref{e.5.6}), let
$R_{t}:=\left(  Q_{t}^{-1}\right)  ^{\ast}$ and observe that
\[
\frac{d}{dt}R_{t}=-\left(  Q_{t}^{-1}\dot{Q}_{t}Q_{t}^{-1}\right)  ^{\ast
}=\frac{1}{2}\left(  Q_{t}^{-1}\operatorname{Ric}^{//_{t}}Q_{t}Q_{t}%
^{-1}\right)  ^{\ast}=\frac{1}{2}\operatorname{Ric}^{//_{t}}R_{t}.
\]
Hence reasoning as above we may conclude that
\[
\left\Vert Q_{t}^{-1}\right\Vert _{op}=\left\Vert \left(  Q_{t}^{-1}\right)
^{\ast}\right\Vert _{op}=\left\Vert R_{t}\right\Vert _{op}\leq e^{Kt/2}.
\]

\end{proof}

When $M$ is compact, the following result is Theorem 5.10 of Driver and
Thalmaier \cite{Driver2001b}.

\begin{thm}
[A divergence formula]\label{t.5.4}Assume the Ricci curvature,
$\operatorname{Ric},$ on $M$ satisfies, $k\leq\operatorname{Ric}\leq K$ for
some $-\infty<k\leq K<\infty.$ Let $T>0$ and $\tilde{\ell}$ be a $C^{1}$ --
adapted real-valued process such that $\tilde{\ell}_{0}=0,$ $\tilde{\ell}%
_{T}=1,$ and
\begin{equation}
\int_{0}^{T}\left\vert \frac{d}{d\tau}\tilde{\ell}_{\tau}\right\vert d\tau\leq
C, \label{e.5.7}%
\end{equation}
where $C<\infty$ is a non-random constant. Then for every $C^{2}$ -- vector
field, $Y,$ on $M$ with compact support the following identity holds
\begin{equation}
\mathbb{E}\left[  \nabla\cdot Y\left(  \Sigma_{T}^{x}\right)  \right]
=\mathbb{E}\left[  \left\langle Y(\Sigma_{T}^{x}),//_{T}Q_{T}\int_{0}%
^{T}\tilde{\ell}_{t}^{\prime}Q_{t}^{-1}db_{t}\right\rangle \right]  ,
\label{e.5.8}%
\end{equation}
where $\nabla\cdot Y$ is the divergence of $Y$ and $\tilde{\ell}_{t}^{\prime
}:=\frac{d}{dt}\tilde{\ell}_{t}.$
\end{thm}

\begin{proof}
The proof will consist of adding some technical details to the proof of
Theorem 5.10 in \cite{Driver2001b}. Suppose $a$ is a smooth one form on $M$
with compact support,
\begin{equation}
a_{t}:=e^{\left(  T-t\right)  \bar{\Delta}_{1}/2}a, \label{e.5.9}%
\end{equation}
$\tilde{\ell}_{\tau}$ is an adapted continuously differentiable real--valued
process, and $\ell_{0}$ is a fixed vector in $T_{x}M.$ Then as shown in
\cite[Theorem 3.4]{Driver2001b} (and repeated below in Lemma \ref{l.C.1} for
the readers convenience) the process,%
\begin{equation}
Z_{t}:=\left(  a_{t}\left(  \Sigma_{t}^{x}\right)  \circ//_{t}\right)
Q_{t}\left[  \int_{0}^{t}Q_{\tau}^{-1}\left(  \frac{d}{d\tau}\tilde{\ell
}_{\tau}\right)  db_{\tau}+\ell_{0}\right]  -\left(  \delta a_{t}\right)
\left(  \Sigma_{t}^{x}\right)  \tilde{\ell}_{t} \label{e.5.10}%
\end{equation}
is a local martingale.

From Theorems \ref{t.4.1} and \ref{t.4.3} we have%
\[
\left\vert a_{t}\right\vert \leq e^{-\left(  T-t\right)  k/2}\left\Vert
a\right\Vert _{\infty}\leq e^{T\left\vert k\right\vert /2}\left\Vert
a\right\Vert _{\infty}%
\]
and%
\[
\left\vert \delta a_{t}\right\vert =\left\vert e^{\left(  T-t\right)
\bar{\Delta}_{0}/2}\delta a\right\vert \leq\left\Vert \delta a\right\Vert
_{\infty}.
\]
Making use of these estimates along with Lemma \ref{l.5.3} and Eq.
(\ref{e.5.7}) shows that $Z_{t}$ is a bounded local martingale and hence, by a
localization argument, a martingale. In particular, it follows that
$t\rightarrow\mathbb{E}Z_{t}$ is constant for $0\leq t\leq T$ and hence%
\begin{align*}
\left(  e^{T\Delta/2}a\right)  \left(  \Sigma_{0}^{x}\right)  \ell_{0}  &
-\delta\left(  e^{T\Delta/2}a\right)  \left(  \Sigma_{0}^{x}\right)
\tilde{\ell}_{0}=Z_{0}=\mathbb{E}Z_{T}\\
&  =\mathbb{E}\left[  \left(  a\left(  \Sigma_{T}^{x}\right)  \circ
//_{T}\right)  Q_{T}\left[  \int_{0}^{T}Q_{\tau}^{-1}\left(  \frac{d}{d\tau
}\tilde{\ell}_{\tau}\right)  db_{\tau}+\ell_{0}\right]  -\delta a\left(
\Sigma_{T}^{x}\right)  \tilde{\ell}_{T}\right]  .
\end{align*}
If we now suppose that $\ell_{0}=0,$ $\tilde{\ell}_{0}=0,$ and $\tilde{\ell
}_{T}=1,$ the above formula reduces to%
\[
0=\mathbb{E}\left[  \left(  a\left(  \Sigma_{T}^{x}\right)  \circ
//_{T}\right)  Q_{T}\int_{0}^{T}Q_{\tau}^{-1}\left(  \frac{d}{d\tau}%
\tilde{\ell}_{\tau}\right)  db_{\tau}-\delta a\left(  \Sigma_{T}^{x}\right)
\right]  .
\]
This identity is equivalent to the identity in Eq. (\ref{e.5.8}) as is seen by
taking $a\left(  x\right)  v:=\left\langle Y\left(  x\right)  ,v\right\rangle
$ for all $x\in M$ and $v\in T_{x}M$ and recalling that
\[
\delta a=\sum_{i=1}^{d}i_{e_{i}}\nabla_{e_{i}}\left\langle Y,\cdot
\right\rangle =\sum_{i=1}^{d}i_{e_{i}}\left\langle \nabla_{e_{i}}%
Y,\cdot\right\rangle =\nabla\cdot Y.
\]

\end{proof}

\begin{ex}
\label{ex.5.5}Taking $\tilde{\ell}_{t}=t/T$ in Eq. (\ref{e.5.8}) shows%
\begin{equation}
\mathbb{E}\left[  \nabla\cdot Y\left(  \Sigma_{T}^{x}\right)  \right]
=\frac{1}{T}\mathbb{E}\left[  \left\langle Y(\Sigma_{T}^{x}),//_{T}Q_{T}%
\int_{0}^{T}Q_{t}^{-1}db_{t}\right\rangle \right]  . \label{e.5.11}%
\end{equation}

\end{ex}

\section{Exponential integrability of $W_{A}^{T}$\label{s.6}}

In this section and for the remainder of the paper we will again go back to
the setting where $M=G$ is a connected uni-modular Lie group equipped with a
left - invariant Riemannian metric as described in the introduction. We are
now going to use Theorem \ref{t.5.4} to estimate $W_{A}:=W_{A}^{T}$ in
Definition \ref{d.1.3}. In order to do this we will use Eq. (\ref{e.5.8}) to
find a useful path integral expression for $W_{A},$ see Theorem \ref{t.6.4} below.

For $A,B\in\mathfrak{g},$ let $D_{A}B:=\nabla_{A}\tilde{B}\in\mathfrak{g}$
where $\nabla$ is the Levi-Civita covariant derivative on $TG.$ Observe that
$\nabla_{\tilde{A}}\tilde{B}$ is a left invariant vector field and $\left(
\nabla_{\tilde{A}}\tilde{B}\right)  \left(  e\right)  =\nabla_{A}\tilde
{B}=D_{A}B.$ Hence we have the identity, $\nabla_{\tilde{A}}\tilde
{B}=\widetilde{D_{A}B}.$

\begin{lem}
\label{l.6.1}Suppose that $\left\{  A_{i}\right\}  _{i=1}^{\dim\mathfrak{g}}$
is an orthonormal basis for $\mathfrak{g}$ and $G$ is uni-modular. Then

\begin{enumerate}
\item $\sum_{i=1}^{\dim\mathfrak{g}}D_{A_{i}}A_{i}=0$ or equivalently
$\sum_{i=1}^{\dim\mathfrak{g}}\nabla_{\tilde{A}_{i}}\tilde{A}_{i}=0.$

\item The divergence of $\tilde{B},$ $\nabla\cdot\tilde{B}$, is zero for all
$B\in\mathfrak{g}.$

\item $\Delta_{0}=\sum_{i=1}^{\dim\mathfrak{g}}\tilde{A}_{i}^{2}$ is the
Laplace Beltrami operator on $G.$
\end{enumerate}
\end{lem}

\begin{proof}
\begin{enumerate}
\item The formula for $D_{A}B$ is%
\[
D_{A}B=\frac{1}{2}\left(  ad_{A}B-ad_{A}^{\ast}B-ad_{B}^{\ast}A\right)
\]
and hence $D_{A}A=-ad_{A}^{\ast}A$ and for any $B\in\mathfrak{g}$ we find%
\begin{align*}
\left(  \sum_{i=1}^{\dim\mathfrak{g}}D_{A_{i}}A_{i},B\right)  _{\mathfrak{g}}
&  =-\sum_{i=1}^{\dim\mathfrak{g}}\left(  A_{i},ad_{A_{i}}B\right)
_{\mathfrak{g}}\\
&  =-\sum_{i=1}^{\dim\mathfrak{g}}\left(  A_{i},ad_{B}A_{i}\right)
_{\mathfrak{g}}=-\operatorname{tr}\left(  ad_{B}\right)  .
\end{align*}
Since $G$ is uni-modular, $\det\left(  Ad_{e^{tB}}\right)  =0$ for all $t$ and
therefore $\operatorname{tr}\left(  ad_{B}\right)  =0.$

\item The following simple computation shows $\nabla\cdot\tilde{B}=0$%
\begin{align*}
\nabla\cdot\tilde{B}  &  =\sum_{i=1}^{\dim\mathfrak{g}}\left(  \nabla
_{\tilde{A}_{i}}\tilde{B},\tilde{A}_{i}\right)  _{TG}=\sum_{i=1}%
^{\dim\mathfrak{g}}\left(  D_{A_{i}}B,A_{i}\right)  _{\mathfrak{g}}\\
&  =-\sum_{i=1}^{\dim\mathfrak{g}}\left(  B,D_{A_{i}}A_{i}\right)
_{\mathfrak{g}}=0.
\end{align*}

\item Observe that $\left\{  \tilde{A}_{i}\right\}  _{i=1}^{\dim\mathfrak{g}}$
is a globally defined orthonormal frame for $TG$ and that
\[
\Delta_{0}=\sum_{i=1}^{\dim\mathfrak{g}}\left[  \tilde{A}_{i}^{2}%
-\nabla_{\tilde{A}_{i}}\tilde{A}_{i}\right]  =\sum_{i=1}^{\dim\mathfrak{g}%
}\tilde{A}_{i}^{2}.
\]

\end{enumerate}

\end{proof}

In Theorem \ref{t.6.4} below, we will specialize Theorem \ref{t.5.4} in order
to find a probabilistic representation for $W_{A}$ of Definition \ref{d.1.3}.
This representation will then be used to estimate $\int_{G}e^{W_{A}}d\nu_{T}$
for all $A\in\mathfrak{g}.$ Let $\left\{  \Sigma_{t}\right\}  _{t\geq0}$ be a
Brownian motion on $G$ such that $\Sigma_{0}=e,$ $b_{t}$ be the $\mathfrak{g}$
-- valued Brownian motion defined by,
\[
b_{t}:=\int_{0}^{t}//_{\tau}\left(  \Sigma\right)  ^{-1}\circ d\Sigma_{\tau},
\]
and $\beta_{t}$ be the $\mathfrak{g}$ -- valued semi-martingale defined by%
\[
\beta_{t}:=\int_{0}^{t}\theta\left(  \circ d\Sigma_{\tau}\right)  =\int
_{0}^{t}L_{\Sigma_{\tau}^{-1}\ast}\circ d\Sigma_{\tau},
\]
where $\theta\left(  v_{g}\right)  :=L_{g^{-1}\ast}v_{g}$ for all $v_{g}\in
T_{g}G.$ As a reflection of the fact that $\sum_{i=1}^{\dim\mathfrak{g}}%
\tilde{A}_{i}^{2}$ is the Laplace--Beltrami operator, $\beta_{t}$ is another
$\mathfrak{g}$--valued Brownian motion. This will also be evident from the
following proposition.

\begin{prop}
\label{p.6.2}Fix $T>0$ and let $U_{t}\in O\left(  \mathfrak{g}\right)  $ be
the unique solution to the stochastic differential equation%
\begin{equation}
dU_{t}+D_{\circ d\beta_{t}}U_{t}=0\text{ with }U_{0}=I. \label{e.6.1}%
\end{equation}
Further define $Y_{t}:=U_{t}Q_{t},$ and $V_{t}:=Y_{T}Y_{t}^{-1}.$ Then
\begin{equation}
//_{t}:=L_{\Sigma_{t}\ast}U_{t} \label{e.6.2}%
\end{equation}
and
\begin{equation}
\int_{0}^{t}U_{\tau}^{-1}\circ d\beta_{\tau}=\int_{0}^{t}U_{\tau}^{-1}%
d\beta_{\tau}=\int_{0}^{t}//_{\tau}{}^{-1}\circ d\Sigma_{\tau}=b_{t}.
\label{e.6.3}%
\end{equation}

\end{prop}

\begin{proof}
The fact that $//_{t}:=L_{\Sigma_{t}\ast}U_{t}$ is explained in \cite[Theorem
6.6]{Driver1997a} and hence%
\[
b_{t}=\int_{0}^{t}U_{t}^{-1}L_{\Sigma_{t}\ast}^{-1}\circ d\Sigma_{\tau}%
=\int_{0}^{t}U_{t}^{-1}\theta\left(  \circ d\Sigma_{\tau}\right)  =\int
_{0}^{t}U_{\tau}^{-1}\circ d\beta_{\tau},
\]
i.e. $d\beta_{t}=U_{t}\circ db_{t}.$ Letting $\left\{  A_{i}\right\}
_{i=1}^{\dim\mathfrak{g}}$ be an orthonormal basis for $\mathfrak{g,}$ it
follows from Lemma \ref{l.6.1} and the fact that $\left\{  U_{t}A_{i}\right\}
_{i=1}^{\dim\mathfrak{g}}$ is also an orthonormal basis for $\mathfrak{g}$
that
\begin{align*}
dU_{t}db_{t}  &  =-\frac{1}{2}D_{d\beta_{t}}U_{t}db_{t}=-\frac{1}{2}%
D_{U_{t}db_{t}}U_{t}db_{t}\\
&  =-\frac{1}{2}\sum_{i=1}^{\dim\mathfrak{g}}D_{U_{t}A_{i}}U_{t}A_{i}\,dt=0.
\end{align*}
This allows us to conclude that $d\beta_{t}=U_{t}\circ db_{t}=U_{t}db_{t}$
which completes the proof of the proposition.
\end{proof}

\begin{prop}
\label{p.6.3}Let $Y_{t}:=U_{t}Q_{t}$ and for fixed $T>0$ let $V_{t}%
:=Y_{T}Y_{t}^{-1}$ and $\mathcal{G}_{t}:=\overline{\sigma\left(  \beta_{\tau
}-\beta_{s}:t\leq s,\tau\leq T\right)  }$ -- the completion of the $\sigma$ --
algebra generated by $\left\{  \beta_{\tau}-\beta_{s}:t\leq s,\tau\leq
T\right\}  .$ Then

\begin{enumerate}
\item $V_{t}$ is $\mathcal{G}_{t}$ -- measurable, and

\item $V_{t}$ is the unique solution to the backwards stochastic differential
equation,%
\[
dV_{t}=V_{t}\left(  D_{\circ d\beta_{t}}+\frac{1}{2}\operatorname{Ric}%
_{e}\,dt\right)  \text{ with }V_{T}=I.
\]

\end{enumerate}
\end{prop}

\begin{proof}
Because $L_{\Sigma_{t}\ast}$ is an isometry of $G,$ it follows that%
\begin{equation}
\operatorname{Ric}^{//_{t}}=//_{t}^{-1}\operatorname{Ric}_{\Sigma_{t}}%
//_{t}=U_{t}^{-1}L_{\Sigma_{t}\ast}^{-1}\operatorname{Ric}_{\Sigma_{t}%
}L_{\Sigma_{t}\ast}U_{t}=U_{t}^{-1}\operatorname{Ric}_{e}U_{t}. \label{e.6.4}%
\end{equation}
Using this identity and the definition of $Y_{t}$ we find, $Y_{0}=Id$ and
\begin{align}
dY_{t}  &  =-D_{\circ d\beta_{t}}U_{t}Q_{t}-\frac{1}{2}U_{t}\operatorname{Ric}%
^{//_{t}}Q_{t}dt\nonumber\\
&  =-D_{\circ d\beta_{t}}Y_{t}-\frac{1}{2}U_{t}\operatorname{Ric}^{//_{t}%
}U_{t}^{-1}Y_{t}dt\label{e.6.5}\\
&  =-D_{\circ d\beta_{t}}Y_{t}-\frac{1}{2}\operatorname{Ric}_{e}\,Y_{t}dt.
\label{e.6.6}%
\end{align}
Since $dY_{t}^{-1}=-Y_{t}^{-1}\left(  \circ dY_{t}\right)  Y_{t}^{-1},$ it
follows that $Y_{t}^{-1}$ satisfies,%
\begin{equation}
dY_{t}^{-1}=Y_{t}^{-1}D_{\circ d\beta_{t}}+\frac{1}{2}Y_{t}^{-1}%
\operatorname{Ric}_{e}\,dt\text{ with }Y_{0}^{-1}=Id. \label{e.6.7}%
\end{equation}

For $T\geq t\geq0,$ let $Y_{T,t}$ solve,%
\[
d_{T}Y_{T,t}=-D_{\circ d\beta_{T}}Y_{T,t}-\frac{1}{2}\operatorname{Ric}%
_{e}\,Y_{T,t}dT\text{ with }Y_{t,t}=Id,
\]
and observe that $Y_{T,t}$ is $\overline{\sigma\left(  \beta_{\tau}-\beta
_{s}:t\leq s,\tau\leq T\right)  }$ -- measurable. By the uniqueness of
solutions to linear stochastic differential equations we may conclude%
\[
Y_{T}=Y_{T,t}Y_{t}\text{ a.s. for all }0\leq t\leq T
\]
and hence it follows that $V_{t}=Y_{T}Y_{t}^{-1}\overset{\text{a.s.}}%
{=}Y_{T,t}$ is also $\overline{\sigma\left(  \beta_{\tau}-\beta_{s}:t\leq
s,\tau\leq T\right)  }$ -- measurable. Moreover we have,
\begin{align*}
dV_{t}  &  =Y_{T}d\left(  Y_{t}^{-1}\right)  =-Y_{T}Y_{t}^{-1}\left(  \circ
dY_{t}\right)  Y_{t}^{-1}\\
&  =-V_{t}\left(  -D_{\circ d\beta_{t}}-\frac{1}{2}\operatorname{Ric}%
_{e}\,dt\right) \\
&  =V_{t}\left(  D_{\circ d\beta_{t}}+\frac{1}{2}\operatorname{Ric}%
_{e}\,dt\right)  \text{ with }V_{T}=Id.
\end{align*}
See \cite[Section 4.1]{Driver1997a} for more on the backwards stochastic
integral interpretation of this equation.
\end{proof}

\begin{thm}
\label{t.6.4}If $A\in\mathfrak{g}$ and $\ell\in C^{1}\left(  \left[
0,T\right]  ,\mathbb{R}\right)  $ with $\ell\left(  0\right)  =0$ and
$\ell\left(  T\right)  =1,$ then%
\begin{equation}
W_{A}\left(  x\right)  =\mathbb{E}\left[  \left.  \left(  A,\int_{0}^{T}%
\dot{\ell}\left(  \tau\right)  V_{\tau}d\overleftarrow{\beta}_{\tau}\right)
\right\vert \Sigma_{T}=x\right]  , \label{e.6.8}%
\end{equation}
where $\int_{0}^{T}\dot{\ell}\left(  \tau\right)  V_{\tau}d\overleftarrow
{\beta}_{\tau}$ is a backwards It\^{o} integral and $V_{t}$ satisfies the
(backwards) stochastic differential equation,%
\[
dV_{t}=\frac{1}{2}V_{t}\operatorname{Ric}_{e}\,dt+V_{t}D_{\circ d\beta_{t}%
}\text{ with }V_{T}=Id.
\]

\end{thm}

\begin{proof}
Let $f\in C_{c}^{\infty}\left(  G\right)  $ and
\[
Y\left(  x\right)  :=f\left(  x\right)  \tilde{A}\left(  x\right)  =f\left(
x\right)  L_{x\ast}A.
\]
As shown in Lemma \ref{l.6.1}, $\nabla\cdot\tilde{A}=0$ from which it follows
that%
\[
\nabla\cdot Y=\left(  \operatorname*{grad}f,\tilde{A}\right)  _{TG}=\tilde
{A}f.
\]
Therefore an application of Theorem \ref{t.5.4} (with $\tilde{\ell}_{t}$ now
being denoted by $\ell\left(  t\right)  )$ shows,%
\begin{align}
\mathbb{E}\left[  \left(  \tilde{A}f\right)  \left(  \Sigma_{T}\right)
\right]   &  =\mathbb{E}\left[  f\left(  \Sigma_{T}\right)  \left\langle
\tilde{A}\left(  \Sigma_{T}\right)  ,//_{T}Q_{T}\int_{0}^{T}\dot{\ell}\left(
\tau\right)  Q_{\tau}^{-1}db_{\tau}\right\rangle \right] \nonumber\\
&  =\mathbb{E}\left[  f\left(  \Sigma_{T}\right)  \left\langle A,L_{\Sigma
_{T}^{-1}\ast}//_{T}Q_{T}\int_{0}^{T}\dot{\ell}\left(  \tau\right)  Q_{\tau
}^{-1}db_{\tau}\right\rangle \right]  . \label{e.6.9}%
\end{align}
From Eq. (\ref{e.6.3})%
\begin{align}
\left\langle A,L_{\Sigma_{T}^{-1}\ast}//_{T}Q_{T}\int_{0}^{T}\dot{\ell}\left(
\tau\right)  Q_{\tau}^{-1}db_{\tau}\right\rangle  &  =\left\langle
A,U_{T}Q_{T}\int_{0}^{T}\dot{\ell}\left(  \tau\right)  Q_{\tau}^{-1}U_{\tau
}^{-1}d\beta_{\tau}\right\rangle \nonumber\\
&  =\left\langle A,Y_{T}\int_{0}^{T}\dot{\ell}\left(  \tau\right)  Y_{\tau
}^{-1}d\beta_{\tau}\right\rangle \label{e.6.10}\\
&  =\left\langle A,\int_{0}^{T}\dot{\ell}\left(  \tau\right)  V_{\tau}%
d\beta_{\tau}\right\rangle . \label{e.6.11}%
\end{align}
Moreover, we may write the last expression as a backwards It\^{o} integral,
since%
\[
dV_{\tau}d\beta_{\tau}=V_{\tau}D_{d\beta_{\tau}}d\beta_{\tau}=V_{\tau}%
\sum_{A\in ONB\left(  \mathfrak{g}\right)  }D_{A}A\cdot dt=0
\]
wherein we have used Lemma \ref{l.6.1} again for the last equality. Hence we
now have%
\[
\left\langle A,L_{\Sigma_{T}^{-1}\ast}//_{T}Q_{T}\int_{0}^{T}\dot{\ell}\left(
\tau\right)  Q_{\tau}^{-1}db_{\tau}\right\rangle =\left\langle A,\int_{0}%
^{T}\dot{\ell}\left(  \tau\right)  V_{\tau}d\overleftarrow{\beta}_{\tau
}\right\rangle .
\]
These computations may be justified by the same methods introduced in
\cite{Driver1997a}. This completes the proof because,%
\[
\mathbb{E}\left[  W_{A}\left(  \Sigma_{T}\right)  f\left(  \Sigma_{T}\right)
\right]  =\mathbb{E}\left[  \left(  \tilde{A}f\right)  \left(  \Sigma
_{T}\right)  \right]  =\mathbb{E}\left[  f\left(  \Sigma_{T}\right)
\left\langle A,\int_{0}^{T}\dot{\ell}\left(  \tau\right)  V_{\tau
}d\overleftarrow{\beta}_{\tau}\right\rangle \right]
\]
for all $f\in C_{c}^{\infty}\left(  G\right)  .$
\end{proof}

Our next goal is to bound $\int_{G}e^{W_{A}}d\nu_{T}$ for all $A\in
\mathfrak{g}.$ In order to do this it will be necessary to estimate the size
of the process $V_{t}.$

\begin{lem}
\label{l.6.5}Suppose $k\in\mathbb{R}$ is chosen so that $\operatorname{Ric}%
\geq kI,$ then
\begin{equation}
\left\vert V_{t}^{\ast}A\right\vert ^{2}\leq\left\vert A\right\vert
^{2}e^{-k\left(  T-t\right)  }\text{ for all }A\in\mathfrak{g}. \label{e.6.12}%
\end{equation}

\end{lem}

\begin{proof}
Since%
\[
dV_{t}=\frac{1}{2}V_{t}\operatorname{Ric}_{e}dt+V_{t}D_{\circ d\beta_{t}},
\]
we have%
\[
dV_{t}^{\ast}=\frac{1}{2}\operatorname{Ric}_{e}V_{t}^{\ast}dt-D_{\circ
d\beta_{t}}V_{t}^{\ast}%
\]
wherein we have used the fact that $D_{A}:\mathfrak{g}\rightarrow\mathfrak{g}$
is antisymmetric. In particular it now follows that
\begin{align*}
d\left\vert V_{t}^{\ast}A\right\vert ^{2}  &  =2\left(  \circ dV_{t}^{\ast
}A,V_{t}^{\ast}A\right)  =2\left(  \frac{1}{2}\operatorname{Ric}_{e}%
V_{t}^{\ast}Adt-D_{\circ d\beta_{t}}V_{t}^{\ast}A,V_{t}^{\ast}A\right) \\
&  =\left(  \operatorname{Ric}_{e}V_{t}^{\ast}A,V_{t}^{\ast}A\right)  dt\geq
k\left\vert V_{t}^{\ast}A\right\vert ^{2}dt\text{ with }\left\vert V_{T}%
^{\ast}A\right\vert ^{2}=\left\vert A\right\vert ^{2}.
\end{align*}
We may write this inequality as%
\[
\frac{d}{dt}\ln\left\vert V_{t}^{\ast}A\right\vert ^{2}\geq k\text{ with
}\left\vert V_{T}^{\ast}A\right\vert ^{2}=\left\vert A\right\vert ^{2}%
\]
which upon integration gives,%
\[
\ln\left\vert A\right\vert ^{2}-\ln\left\vert V_{t}^{\ast}A\right\vert
^{2}=\ln\left\vert V_{T}^{\ast}A\right\vert ^{2}-\ln\left\vert V_{t}^{\ast
}A\right\vert ^{2}\geq k\left(  T-t\right)  .
\]
Hence $\left\vert A\right\vert ^{2}/\left\vert V_{t}^{\ast}A\right\vert
^{2}\geq e^{k\left(  T-t\right)  }$ which is equivalent to Eq. (\ref{e.6.12}).
\end{proof}

\begin{lem}
\label{l.6.6}Let $k\in\mathbb{R}$ and $T>0,$ then%
\begin{equation}
\inf\left\{  \int_{0}^{T}\dot{\ell}^{2}\left(  \tau\right)  e^{-k\left(
T-\tau\right)  }d\tau\right\}  \leq\frac{k}{e^{kT}-1} \label{e.6.13}%
\end{equation}
where the infimum is taken over all $\ell\in C^{1}\left(  \left[  0,T\right]
,\mathbb{R}\right)  $ such that $\ell\left(  0\right)  =0$ and $\ell\left(
T\right)  =1.$
\end{lem}

\begin{proof}
By a simple calculus of variation argument, $\ell\in C^{1}\left(  \left[
0,T\right]  ,\mathbb{R}\right)  $ with $\ell\left(  0\right)  =0$ and
$\ell\left(  T\right)  =1$ is a critical point for the function,
\begin{equation}
K\left(  \ell\right)  :=\int_{0}^{T}\dot{\ell}^{2}\left(  \tau\right)
e^{-k\left(  T-\tau\right)  }d\tau, \label{e.6.14}%
\end{equation}
iff $\dot{\ell}\left(  \tau\right)  e^{k\tau}$ is constant in $\tau.$ This
constraint and the boundary conditions imply that $K$ has a unique critical
point at
\[
\ell_{c}\left(  \tau\right)  =\frac{e^{-k\tau}-1}{e^{-kT}-1}.
\]
Plugging this value of $\ell_{c}$ into $K$ then shows $K\left(  \ell
_{c}\right)  =k\left(  1-e^{-kT}\right)  ^{-1}$ from which Eq. (\ref{e.6.13}) follows.
\end{proof}

\subsection{Proof of Theorems \ref{t.1.4} and \ref{t.1.5}}

With the above results as preparation, we are now in position to complete the
proofs of Theorem \ref{t.1.4} and \ref{t.1.5}.

\begin{proof}
Proof of Theorem \ref{t.1.4}. Let $\ell\in C^{1}\left(  \left[  0,T\right]
,\mathbb{R}\right)  $ such that $\ell\left(  0\right)  =0$ and $\ell\left(
T\right)  =1.$ From Theorem \ref{t.6.4}, Lemma \ref{l.6.5}, Jensen's
inequality for conditional expectations, and a standard martingale argument
(see the proof of Lemma 7.6 and especially Eq. 7.17 in \cite{Driver1997b}) we
have%
\begin{align*}
\int_{G}e^{W_{A}}d\nu_{T}  &  =\mathbb{E}\left[  \exp\left(  \mathbb{E}\left[
\left.  \left\langle A,\int_{0}^{T}\dot{\ell}\left(  \tau\right)  V_{\tau
}d\overleftarrow{\beta}_{\tau}\right\rangle \right\vert \sigma\left(
\Sigma_{T}\right)  \right]  \right)  \right] \\
&  \leq\mathbb{E}\left[  \mathbb{E}\left[  \exp\left.  \left(  \left\langle
A,\int_{0}^{T}\dot{\ell}\left(  \tau\right)  V_{\tau}d\overleftarrow{\beta
}_{\tau}\right\rangle \right)  \right\vert \sigma\left(  \Sigma_{T}\right)
\right]  \right] \\
&  =\mathbb{E}\left[  \exp\left(  \left\langle A,\int_{0}^{T}\dot{\ell}\left(
\tau\right)  V_{\tau}d\overleftarrow{\beta}_{\tau}\right\rangle \right)
\right] \\
&  \leq\exp\left(  \frac{1}{2}\left\Vert \int_{0}^{T}\dot{\ell}^{2}\left(
\tau\right)  \left\vert V_{\tau}^{\ast}A\right\vert ^{2}d\tau\right\Vert
_{L^{\infty}\left(  P\right)  }\right) \\
&  \leq\exp\left(  \frac{\left\vert A\right\vert ^{2}}{2}\int_{0}^{T}\dot
{\ell}^{2}\left(  \tau\right)  e^{-k\left(  T-\tau\right)  }d\tau\right)  ,
\end{align*}
where $P$ is the underlying probability measure. Since $\ell$ was arbitrary,
it follows from Lemma \ref{l.6.6} that,%
\begin{align*}
\int_{G}e^{W_{A}}d\nu_{T}  &  \leq\inf_{\ell}\exp\left(  \frac{1}{2}\int
_{0}^{T}\dot{\ell}^{2}\left(  \tau\right)  \left\vert A\right\vert
^{2}e^{-k\left(  T-\tau\right)  }d\tau\right) \\
&  \leq\exp\left(  \frac{1}{2}\frac{k}{e^{kT}-1}\left\vert A\right\vert
^{2}\right)  =\exp\left(  \frac{1}{2T}c\left(  kT\right)  \left\vert
A\right\vert ^{2}\right)  .
\end{align*}

\end{proof}

\begin{proof}
(Proof of Theorem \ref{t.1.5}.) From Theorem \ref{t.6.4}, Lemma \ref{l.6.5},
Jensen's inequality for conditional expectations, and Burkholder-Davis-Gundy
inequality (see for example \cite[Corollary 6.3.1a on p.344]{StroockBook},
\cite[Appendix A.2]{Nualart95}, or \cite[p. 212]{Metivier82} and \cite[Theorem
17.7]{Kall} for the real case), there exists $C_{p}<\infty$ such that%
\begin{align*}
\int_{G}\left\vert W_{A}\right\vert ^{p}d\nu_{T}  &  =\mathbb{E}\left[
\left\vert \mathbb{E}\left[  \left.  \left\langle A,\int_{0}^{T}\dot{\ell
}\left(  \tau\right)  V_{\tau}d\overleftarrow{\beta}_{\tau}\right\rangle
\right\vert \sigma\left(  \Sigma_{T}\right)  \right]  \right\vert ^{p}\right]
\\
&  \leq\mathbb{E}\left[  \mathbb{E}\left[  \left.  \left\vert \left\langle
A,\int_{0}^{T}\dot{\ell}\left(  \tau\right)  V_{\tau}d\overleftarrow{\beta
}_{\tau}\right\rangle \right\vert ^{p}~\right\vert \sigma\left(  \Sigma
_{T}\right)  \right]  \right] \\
&  =\mathbb{E}\left[  \left\vert \left\langle A,\int_{0}^{T}\dot{\ell}\left(
\tau\right)  V_{\tau}d\overleftarrow{\beta}_{\tau}\right\rangle \right\vert
^{p}\right]  =\mathbb{E}\left[  \left\vert \int_{0}^{T}\dot{\ell}\left(
\tau\right)  \left\langle V_{\tau}^{\ast}A,d\overleftarrow{\beta}_{\tau
}\right\rangle \right\vert ^{p}\right] \\
&  \leq C_{p}^{p}\mathbb{E}\left[  \left\vert \int_{0}^{T}\dot{\ell}%
^{2}\left(  \tau\right)  \left\vert V_{\tau}^{\ast}A\right\vert ^{2}%
d\tau\right\vert ^{p/2}\right] \\
&  \leq C_{p}^{p}\left(  \left\vert A\right\vert ^{2}\int_{0}^{T}\dot{\ell
}^{2}\left(  \tau\right)  e^{-k\left(  T-\tau\right)  }d\tau\right)  ^{p/2}.
\end{align*}
Using Lemma \ref{l.6.6}, we may optimize this last estimate over the
admissible $\ell$ to find,%
\[
\int_{G}\left\vert W_{A}\right\vert ^{p}d\nu_{T}\leq C_{p}^{p}\left(
\left\vert A\right\vert ^{2}\frac{k}{e^{kT}-1}\right)  ^{p/2}=C_{p}^{p}\left(
\left\vert A\right\vert ^{2}\frac{c\left(  kT\right)  }{T}\right)  ^{p/2}%
\]
which is equivalent to Eq. (\ref{e.1.5}).
\end{proof}

\section{Applications\label{s.7}}

\begin{lem}
\label{l.7.1}Suppose that $T>0,$ $p>1,$ and $f\in L^{p}\left(  \nu_{T}\right)
\cap C^{2}\left(  G\right)  $ such that $\Delta f\in L^{p}\left(  \nu
_{T}\right)  .$ Then $f,\Delta f\in L^{p}\left(  \nu_{t}\right)  $ for
$0<t\leq T$ and
\begin{equation}
\frac{\partial}{\partial t}\int_{G}p_{t}\left(  x,y\right)  f\left(  y\right)
dy=\frac{1}{2}\int_{G}p_{t}\left(  x,y\right)  \Delta f\left(  y\right)
dy~\text{ for all }~0<t<T. \label{e.7.1}%
\end{equation}

\end{lem}

\begin{proof}
Since the Ricci curvature is left translation invariant, it is bounded on $G.$
Applying the Li -- Yau Harnack inequality (see Eq. (\ref{e.D.6} below), we
have for any $\gamma>1/2$ that there exists $K=K\left(  \gamma,T\right)
<\infty$ such that
\begin{equation}
p_{t}(x)\leq K\left(  \frac{T}{t}\right)  ^{d\gamma}p_{T}(x)\text{~}%
\forall~\left(  x,t\right)  \in G\times(0,T]. \label{e.7.2}%
\end{equation}
In particular it follows that
\begin{equation}
\left\Vert f\right\Vert _{L^{p}\left(  \nu_{t}\right)  }\leq K\left(  \frac
{T}{t}\right)  ^{d\gamma/p}\left\Vert f\right\Vert _{L^{p}\left(  \nu
_{T}\right)  }\text{~}\forall~\text{ }0<t\leq T. \label{e.7.3}%
\end{equation}

Using $p^{\prime}-1=\left(  p-1\right)  ^{-1}$ and Eq. (\ref{e.1.8}), it
follows that
\begin{align}
\int_{G}p_{t}\left(  y,x\right)  \left\vert f\left(  x\right)  \right\vert dx
&  =\int_{G}\frac{p_{t}\left(  y,x\right)  }{p_{t}\left(  x\right)
}\left\vert f\left(  x\right)  \right\vert d\nu_{t}\left(  x\right)
\nonumber\\
&  \leq\left\Vert \frac{p_{t}\left(  y,\cdot\right)  }{p_{t}\left(
\cdot\right)  }\right\Vert _{L^{p^{\prime}}\left(  \nu_{t}\right)  }%
\cdot\left\Vert f\right\Vert _{L^{p}\left(  \nu_{t}\right)  }\nonumber\\
&  \leq\left\Vert f\right\Vert _{L^{p}\left(  \nu_{t}\right)  }\exp\left(
\frac{c\left(  kt\right)  \left(  p^{\prime}-1\right)  }{2t}\left\vert
y\right\vert ^{2}\right) \nonumber\\
&  \leq\left\Vert f\right\Vert _{L^{p}\left(  \nu_{t}\right)  }\exp\left(
\frac{c\left(  kt\right)  }{2t\left(  p-1\right)  }\left\vert y\right\vert
^{2}\right)  . \label{e.7.4}%
\end{align}
Therefore the integrals in Eq. (\ref{e.7.1}) are well defined. Moreover,
\begin{align*}
\int_{G}p_{t}\left(  y,x\right)  f\left(  x\right)  dx  &  =\int_{G}%
p_{t}\left(  y^{-1}x\right)  f\left(  x\right)  dx=\int_{G}p_{t}\left(
x\right)  f\left(  yx\right)  dx\\
&  =\int_{G}f\circ L_{y}\left(  x\right)  p_{t}\left(  x\right)  dx
\end{align*}
and for any $q\in\left(  1,p\right)  ,$%
\begin{align*}
\left\Vert f\circ L_{y}\right\Vert _{L^{q}\left(  \nu_{t}\right)  }^{q}  &
=\int_{G}\left\vert f\left(  yx\right)  \right\vert ^{q}p_{t}\left(  x\right)
dx=\int_{G}\left\vert f\left(  x\right)  \right\vert ^{q}p_{t}\left(
y^{-1}x\right)  dx\\
&  =\int_{G}\left\vert f\left(  x\right)  \right\vert ^{q}\frac{p_{t}\left(
y^{-1}x\right)  }{p_{t}\left(  x\right)  }d\nu_{t}\left(  x\right) \\
&  \leq\left\Vert f\right\Vert _{L^{p}\left(  \nu_{t}\right)  }\exp\left(
\frac{c\left(  kt\right)  p\left(  p-q\right)  ^{-1}}{2t}\left\vert
y\right\vert ^{2}\right)
\end{align*}
wherein we have used H\"{o}lder's inequality and Eq. (\ref{e.1.11}) for the
last inequality. From these remarks and the fact that $\Delta\left(  f\circ
L_{y}\right)  =\left(  \Delta f\right)  \circ L_{y},$ it suffices to prove Eq.
(\ref{e.7.1}) in the special case where $y=e.$

From Eq. (\ref{e.7.2}) and the Dominated convergence theorem, the function,%
\[
F\left(  t\right)  =\int_{G}f\left(  x\right)  d\nu_{t}\left(  x\right)
~\text{ for all }~t\in(0,T],
\]
is continuous. Our goal now is to show $F$ is differentiable and that $\dot
{F}\left(  t\right)  =\frac{1}{2}\int_{G}\Delta f\left(  x\right)  d\nu
_{t}\left(  x\right)  $ for all $0<t<T.$ To prove this suppose that $h\in
C_{c}^{\infty}\left(  G\right)  $ and consider,%
\[
F_{h}\left(  t\right)  :=\int_{G}f\left(  x\right)  h\left(  x\right)
p_{t}\left(  x\right)  dx.
\]
To simplify notation in the computation below, let $\left\{  A_{i}\right\}
_{i=1}^{\dim\mathfrak{g}}$ be an orthonormal basis for $\mathfrak{g,}$ $\nabla
f=\left(  \tilde{A}_{i}f\right)  _{i=1}^{\dim\mathfrak{g}},$ and $\nabla\cdot
U=\sum\tilde{A}_{i}U_{i}$ where $U=\left(  U_{i}\right)  _{i=1}^{\dim
\mathfrak{g}}$ with $U_{i}\in C^{\infty}\left(  G\right)  .$ Using
$\frac{\partial}{\partial t}p_{t}\left(  x\right)  =\frac{1}{2}\Delta
p_{t}\left(  x\right)  ,$ and a few integration by parts we find%
\begin{align}
\dot{F}_{h}\left(  t\right)   &  =\frac{1}{2}\int_{G}f\left(  x\right)
h\left(  x\right)  \Delta p_{t}\left(  x\right)  dx\nonumber\\
&  =\frac{1}{2}\int_{G}\Delta\left(  fh\right)  ~p_{t}\,dV=\frac{1}{2}\int
_{G}\left(  f\Delta h+2\nabla f\cdot\nabla h+h\Delta f\right)  ~p_{t}%
\,dV\nonumber\\
&  =\frac{1}{2}\int_{G}\left(  f\Delta h+h\Delta f\right)  ~p_{t}\,dV-\int
_{G}f~\nabla\cdot\left[  \nabla h~p_{t}\right]  \,dV\nonumber\\
&  =\frac{1}{2}\int_{G}\left(  f\Delta h+h\Delta f\right)  ~p_{t}\,dV-\int
_{G}f\left[  \Delta h~p_{t}+\nabla h\cdot\nabla p_{t}\right]  \,dV\nonumber\\
&  =-\frac{1}{2}\int_{G}f\Delta h~d\nu_{t}-\int_{G}f~\nabla h\cdot\frac{\nabla
p_{t}}{p_{t}}d\nu_{t}+\frac{1}{2}\int_{G}h\Delta f~d\nu_{t}. \label{e.7.5}%
\end{align}
Therefore,%
\[
\dot{F}_{h}\left(  t\right)  -\frac{1}{2}\int_{G}\Delta f~d\nu_{t}=-\frac
{1}{2}R_{h}\left(  t\right)  -S_{h}\left(  t\right)  +\frac{1}{2}U_{h}\left(
t\right)
\]
where, making use of Eqs. (\ref{e.7.3}) and (\ref{e.1.5}), we have
\begin{align}
\left\vert R_{h}\left(  t\right)  \right\vert  &  \leq\int_{G}\left\vert
f\right\vert \left\vert \Delta h\right\vert ~d\nu_{t}\leq\left\Vert
f\right\Vert _{L^{p}\left(  \nu_{t}\right)  }\left\Vert \Delta h\right\Vert
_{L^{p^{\prime}}\left(  \nu_{t}\right)  }\nonumber\\
&  \leq K^{2}\left(  \frac{T}{t}\right)  ^{d\gamma}\left\Vert f\right\Vert
_{L^{p}\left(  \nu_{T}\right)  }\left\Vert \Delta h\right\Vert _{L^{p^{\prime
}}\left(  \nu_{T}\right)  }, \label{e.7.6}%
\end{align}%
\begin{align}
\left\vert S_{h}\left(  t\right)  \right\vert  &  =\sum_{i}\int_{G}\left\vert
f\right\vert \left\vert \tilde{A}_{i}h\right\vert \left\vert W_{A_{i}}%
^{t}\right\vert d\nu_{t}\leq\sum_{i}\left\Vert f\cdot\tilde{A}_{i}h\right\Vert
_{L^{p}\left(  \nu_{t}\right)  }\left\Vert W_{A_{i}}^{t}\right\Vert
_{L^{p^{\prime}}\left(  \nu_{t}\right)  }\nonumber\\
&  \leq C_{p}\sqrt{\frac{c\left(  kt\right)  }{t}}K\left(  \frac{T}{t}\right)
^{d\gamma/p}\sum_{i}\left\Vert f\cdot\tilde{A}_{i}h\right\Vert _{L^{p}\left(
\nu_{T}\right)  }. \label{e.7.7}%
\end{align}
and%
\begin{align}
\left\vert U_{h}\left(  t\right)  \right\vert  &  \leq\int_{G}\left\vert
\Delta f\right\vert \left\vert h-1\right\vert ~d\nu_{t}\leq\left\Vert \Delta
f\right\Vert _{L^{p}\left(  \nu_{t}\right)  }\left\Vert 1-h\right\Vert
_{L^{p^{\prime}}\left(  \nu_{t}\right)  }\nonumber\\
&  \leq K^{2}\left(  \frac{T}{t}\right)  ^{d\gamma}\left\Vert \Delta
f\right\Vert _{L^{p}\left(  \nu_{T}\right)  }\left\Vert 1-h\right\Vert
_{L^{p^{\prime}}\left(  \nu_{T}\right)  }. \label{e.7.8}%
\end{align}

From \cite[Lemma 3.6]{Driver1997c}, we may choose $\left\{  h_{n}\right\}
_{n=1}^{\infty}\subset C_{c}^{\infty}(G,\left[  0,1\right]  )$ such that
$h_{n}(x)=1$ whenever $\left\vert x\right\vert \leq n$ and $\sup_{n}\sup_{x\in
G}\left\vert \left(  \tilde{A}_{i_{1}}\dots\tilde{A}_{i_{k}}h_{n}\right)
(x)\right\vert <\infty$ for all $i_{1},\dots,i_{k}\in\left\{  1,2,\dots
,\dim\mathfrak{g}\right\}  $ and $k\in\mathbb{N}.$ It then follows from Eqs.
(\ref{e.7.3}), (\ref{e.7.5}), (\ref{e.7.6}), (\ref{e.7.7}), and (\ref{e.7.8})
and the dominated convergence theorem that%
\[
\left\vert \dot{F}_{h_{n}}\left(  t\right)  -\frac{1}{2}\int_{G}\Delta
f~d\nu_{t}\right\vert \leq\frac{1}{2}\left\vert R_{h_{n}}\left(  t\right)
\right\vert +\left\vert S_{h_{n}}\left(  t\right)  \right\vert +\frac{1}%
{2}\left\vert U_{h_{n}}\left(  t\right)  \right\vert \rightarrow0\text{ as
\thinspace}n\rightarrow\infty
\]
uniformly on compact subsets of $\left(  0,T\right)  .$ Moreover, by the
dominated convergence theorem, $F_{h_{n}}\left(  t\right)  \rightarrow
F\left(  t\right)  $ as $n\rightarrow\infty$ and therefore we may conclude
that $\dot{F}\left(  t\right)  =\frac{1}{2}\int_{G}\Delta f~d\nu_{t}$ for
$t\in\left(  0,T\right)  .$
\end{proof}

\subsection{The proof of Proposition \ref{p.1.8}}

\begin{proof}
Now suppose, as in Proposition \ref{p.1.8}, $T>0,$ $p>1,$ and $f\in
L^{p}\left(  \nu_{T}\right)  $ such that $\Delta f=0.$ As in the proof of
Lemma \ref{l.7.1}, we may reduce the proof to the case where $y=e.$ Let
$F\left(  t\right)  :=\int_{G}fd\nu_{t}.$ By Lemma \ref{l.7.1} and the mean
value theorem, $F\left(  T\right)  =F\left(  t\right)  $ for all $t\in\left(
0,T\right)  $ and in particular, $F\left(  T\right)  =\lim_{t\downarrow
0}F\left(  t\right)  .$ We are going to finish the proof by showing
$\lim_{t\downarrow0}F\left(  t\right)  =f\left(  e\right)  .$ To do this, let
$h\in C_{c}^{\infty}\left(  G,\left[  0,1\right]  \right)  $ be chosen so that
$h\left(  x\right)  =1$ if $\left\vert x\right\vert \leq1.$ Then%
\[
F\left(  t\right)  =\int_{G}f\left(  x\right)  h\left(  x\right)  p_{t}\left(
x\right)  dx+r\left(  t\right)
\]
where
\begin{align}
\left\vert r\left(  t\right)  \right\vert  &  \leq\int_{G}\left\vert f\left(
x\right)  \right\vert \left\vert 1-h\left(  x\right)  \right\vert p_{t}\left(
x\right)  dx\leq\int_{\left\vert x\right\vert \geq1}\left\vert f\left(
x\right)  \right\vert p_{t}\left(  x\right)  dx\nonumber\\
&  =\int_{\left\vert x\right\vert \geq1}\left\vert f\left(  x\right)
\right\vert \frac{p_{t}\left(  x\right)  }{p_{T}\left(  x\right)  }d\nu
_{T}\left(  x\right)  \leq\sup_{\left\vert x\right\vert \geq1}\frac
{p_{t}\left(  x\right)  }{p_{T}\left(  x\right)  }\left\Vert f\right\Vert
_{L^{1}\left(  \nu_{T}\right)  }. \label{e.7.9}%
\end{align}
Since $\lim_{t\downarrow0}\int_{G}f\left(  x\right)  h\left(  x\right)
p_{t}\left(  x\right)  dx=f\left(  e\right)  h\left(  e\right)  =f\left(
e\right)  ,$ it suffices to show $\lim_{t\downarrow0}\left\vert r\left(
t\right)  \right\vert =0.$

To estimate $r\left(  t\right)  $ we will make use of some crude upper and
lower bounds on the heat kernel, $p_{t}\left(  x\right)  ,$ for example see
\cite[Theorem V.4.4 or Theorem IX.1.2.]{VSC} for more precise bounds.
According to either of these theorems, there exists a constant $c>0$ such
that
\[
\frac{p_{t}\left(  x\right)  }{p_{T}\left(  x\right)  }\leq\frac{ct^{-d/2}%
\exp\left(  -c\left\vert x\right\vert ^{2}/t\right)  }{c^{-1}T^{-d/2}%
\exp\left(  -c^{-1}\left\vert x\right\vert ^{2}/T\right)  }=c^{2}\left(
\frac{T}{t}\right)  ^{d/2}\exp\left(  \left(  \frac{1}{cT}-\frac{c}{t}\right)
\left\vert x\right\vert ^{2}\right)  .
\]
From this estimate it follows that $\lim_{t\downarrow0}\sup_{\left\vert
x\right\vert \geq1}\left(  p_{t}\left(  x\right)  /p_{T}\left(  x\right)
\right)  =0$ which combined with Eq. (\ref{e.7.9}) shows $\lim_{t\downarrow
0}\left\vert r\left(  t\right)  \right\vert =0.$
\end{proof}

\subsection{Applications to infinite--dimensional groups\label{s.7.2}}

For this section, suppose that $G$ is a topological group, $\mathcal{B}$ is
the Borel $\sigma$ -- algebra over $G,$ and $G_{0}$ is a dense subgroup of $G$
which is endowed with the structure of an infinite--dimensional Hilbert Lie
group. Further assume that $\mathfrak{g}_{0}:=\operatorname*{Lie}\left(
G_{0}\right)  =T_{e}G_{0}$ is equipped with a Hilbertian inner product,
$\left\langle \cdot,\cdot\right\rangle _{\mathfrak{g}_{0}}.$ We will also
assume that $\left(  G,\mathcal{B}\right)  $ is also equipped with a
probability measure, $\nu,$ to be thought of as the \textquotedblleft heat
kernel\textquotedblright\ measure at some time $T>0$ associated to the given
inner product on $\mathfrak{g}_{0}.$ We will now give two theorems which
guarantee that $\nu$ is quasi-invariant under both left and right translations
by elements of $G_{0}.$ The two cases considered are where $G$ can be thought
of as either a projective or inductive limit of finite--dimensional Lie groups.

\begin{thm}
[Projective Limits]\label{t.7.2}Suppose that $T>0,$ $A$ is a directed set,
$\left\{  G_{\alpha}\right\}  _{\alpha\in A}$ is a collection of finite
dimensional uni-modular Lie groups, and $\left\{  \pi_{\alpha}:G\rightarrow
G_{\alpha}\right\}  _{\alpha\in A}$ is a collection of continuous group
homomorphisms satisfying the following properties.

\begin{enumerate}
\item $\mathcal{B}$ is equal to the $\sigma$ -- algebra generated by the
projections, $\left\{  \pi_{\alpha}\right\}  _{\alpha\in A}.$

\item $\pi_{\alpha}|_{G_{0}}:G_{0}\rightarrow G_{\alpha}$ is a smooth
surjection. Let $d\pi_{\alpha}:\mathfrak{g}_{0}\rightarrow\mathfrak{g}%
_{\alpha}$ be the differential of $\pi_{\alpha}$ at $e.$

\item $\nu_{\alpha}:=\left(  \pi_{\alpha}\right)  _{\ast}\nu=\nu\circ
\pi_{\alpha}^{-1}$ is the time $T$ heat kernel measure on $G_{\alpha}$
determined by the unique inner product. $\left(  \cdot,\cdot\right)  _{\alpha
}$ on $\mathfrak{g}_{\alpha}$ which makes%
\[
d\pi_{\alpha}|_{\operatorname*{Nul}\left(  \pi_{\alpha}\right)  }%
:\operatorname*{Nul}\left(  \pi_{\alpha}\right)  ^{\perp}\rightarrow
\mathfrak{g}_{\alpha}%
\]
an isometric isomorphism of inner product spaces.

\item There exists $k\in\mathbb{R}$ such that $\operatorname{Ric}_{\alpha}\geq
kg_{\alpha}$ for all $\alpha\in A,$ where $\operatorname{Ric}_{\alpha}$ is the
Ricci tensor on $G_{\alpha}$ equipped with the left invariant metric
determined by $\left\langle \cdot,\cdot\right\rangle _{\alpha}.$
\end{enumerate}

Under these assumptions, to each $h\in G_{0},$ $\nu\circ R_{h}^{-1}$ is
absolutely continuous relative to $\nu.$ Moreover, if $J_{h}:=d\left(
\nu\circ R_{h}^{-1}\right)  /d\nu$ is the Radon-Nikodym derivative of
$\nu\circ R_{h}^{-1}$ with respect to $\nu$ and $1\leq p<\infty,$ then%
\begin{equation}
\left\Vert J_{h}\right\Vert _{L^{p}\left(  \nu\right)  }\leq\exp\left(
\frac{c\left(  kT\right)  \left(  p-1\right)  }{2T}d_{G_{0}}^{2}\left(
e,h\right)  \right)  , \label{e.7.10}%
\end{equation}
where $d_{G_{0}}$ is the Riemannian distance function on $G_{0}.$
\end{thm}

\begin{proof}
Since the estimate in Eq. (\ref{e.7.10}) holds for $p=1,$ we may assume
without loss of generality that $1<p<\infty.$ Let $\mathbb{H}$ denote the
linear space of bounded measurable functions of the form $f=u\circ\pi_{\alpha
}$ where $\alpha\in A$ and $u:G_{\alpha}\rightarrow\mathbb{R}$ is a bounded
measurable function on $G_{\alpha}.$ Because of assumption 1., $\mathbb{H}$ is
dense in $L^{p}\left(  G,\nu\right)  .$ (An easy proof may be given using a
functional form of the monotone class theorem, see for example \cite[Theorem
A.1 on p. 309]{Janson97}.) By Theorem \ref{t.1.6} in the form of Eq.
(\ref{e.1.9}),%
\[
J_{\alpha}\left(  x\right)  :=\frac{\nu_{\alpha}\left(  dx\cdot\pi_{\alpha
}\left(  h^{-1}\right)  \right)  }{\nu_{\alpha}\left(  dx\right)  }\text{ for
}x\in G_{\alpha},
\]
satisfies%
\[
\left\Vert J_{\alpha}\right\Vert _{L^{p}\left(  G_{\alpha},\nu_{\alpha
}\right)  }\leq\exp\left(  \frac{c\left(  kT\right)  \left(  p-1\right)  }%
{2T}d_{G_{\alpha}}^{2}\left(  e,\pi_{\alpha}\left(  h\right)  \right)
\right)  \text{ for all }1<p<\infty.
\]
Using this result and assumption 3, if $f=u\circ\pi_{\alpha}\in\mathbb{H},$
then
\begin{align*}
\int_{G}\left\vert f\left(  xh\right)  \right\vert d\nu\left(  x\right)   &
=\int_{G}\left\vert u\circ\pi_{\alpha}\left(  xh\right)  \right\vert
d\nu\left(  x\right)  =\int_{G}\left\vert u\left(  \pi_{\alpha}\left(
x\right)  \pi_{\alpha}\left(  h\right)  \right)  \right\vert d\nu\left(
x\right) \\
&  =\int_{G_{\alpha}}\left\vert u\left(  y\cdot\pi_{\alpha}\left(  h\right)
\right)  \right\vert d\nu_{\alpha}\left(  y\right)  =\int_{G_{\alpha}%
}\left\vert u\left(  y\right)  \right\vert J_{\alpha}\left(  y\right)
d\nu_{\alpha}\left(  y\right)  .
\end{align*}
An application of H\"{o}lder's inequality then implies,%
\begin{align}
\int_{G}\left\vert f\left(  xh\right)  \right\vert d\nu\left(  x\right)   &
\leq\left\Vert u\right\Vert _{L^{p}\left(  G_{\alpha},\nu_{\alpha}\right)
}\cdot\left\Vert J_{\alpha}\right\Vert _{L^{p^{\prime}}\left(  G_{\alpha}%
,\nu_{\alpha}\right)  }\nonumber\\
&  \leq\left\Vert f\right\Vert _{L^{p}\left(  G,\nu\right)  }\exp\left(
\frac{c\left(  kT\right)  \left(  p^{\prime}-1\right)  }{2T}d_{G_{\alpha}}%
^{2}\left(  e,\pi_{\alpha}\left(  h\right)  \right)  \right)  . \label{e.7.11}%
\end{align}

Now suppose that $k\in C^{1}\left(  \left[  0,1\right]  ,G_{0}\right)  $ such
that $k\left(  0\right)  =e$ and $k\left(  1\right)  =h.$ Then the length of
$t\rightarrow\pi_{\alpha}\left(  k\left(  t\right)  \right)  \in G_{\alpha}$
is given by
\[
\ell_{G_{\alpha}}\left(  \pi_{\alpha}\circ k\right)  =\int_{0}^{1}\left\vert
L_{\pi_{\alpha}\left(  k\left(  t\right)  \right)  ^{-1}\ast}\pi_{\alpha
}\left(  \dot{k}\left(  t\right)  \right)  \right\vert _{\mathfrak{g}_{\alpha
}}dt.
\]
Since%
\begin{align*}
L_{\pi_{\alpha}\left(  k\left(  t\right)  \right)  ^{-1}\ast}\pi_{\alpha
}\left(  \dot{k}\left(  t\right)  \right)   &  =\frac{d}{ds}|_{0}\pi_{\alpha
}\left(  k\left(  t\right)  \right)  ^{-1}\pi_{\alpha}\left(  k\left(
t+s\right)  \right) \\
&  =\frac{d}{ds}|_{0}\pi_{\alpha}\left(  k\left(  t\right)  ^{-1}k\left(
t+s\right)  \right)  =d\pi_{\alpha}\left(  L_{k\left(  t\right)  ^{-1}\ast
}\dot{k}\left(  t\right)  \right)
\end{align*}
and%
\[
\left\vert L_{\pi_{\alpha}\left(  k\left(  t\right)  \right)  ^{-1}\ast}%
\pi_{\alpha}\left(  \dot{k}\left(  t\right)  \right)  \right\vert
_{\mathfrak{g}_{\alpha}}=\left\vert d\pi_{\alpha}\left(  L_{k\left(  t\right)
^{-1}\ast}\dot{k}\left(  t\right)  \right)  \right\vert _{\mathfrak{g}%
_{\alpha}}\leq\left\vert L_{k\left(  t\right)  ^{-1}\ast}\dot{k}\left(
t\right)  \right\vert _{\mathfrak{g}_{0}},
\]
it follows that
\[
d_{G_{\alpha}}\left(  e,\pi_{\alpha}\left(  h\right)  \right)  \leq
\ell_{G_{\alpha}}\left(  \pi_{\alpha}\circ k\right)  \leq\int_{0}%
^{1}\left\vert L_{k\left(  t\right)  ^{-1}\ast}\dot{k}\left(  t\right)
\right\vert _{\mathfrak{g}_{0}}dt=\ell_{G_{0}}\left(  k\right)  .
\]
Taking the infimum over all such $k$ implies%
\[
d_{G_{\alpha}}\left(  e,\pi_{\alpha}\left(  h\right)  \right)  \leq d_{G_{0}%
}\left(  e,h\right)  .
\]
Combining this inequality with Eq. (\ref{e.7.11}) gives the estimate,
\begin{equation}
\int_{G}\left\vert f\left(  xh\right)  \right\vert d\nu\left(  x\right)
\leq\left\Vert f\right\Vert _{L^{p}\left(  G,\nu\right)  }\exp\left(
\frac{c\left(  kT\right)  \left(  p^{\prime}-1\right)  }{2T}d_{G_{0}}%
^{2}\left(  e,h\right)  \right)  . \label{e.7.12}%
\end{equation}

The afore mentioned density of $\mathbb{H}$ in $L^{p}\left(  G,\nu\right)  $
along with Eq. (\ref{e.7.12}) shows the linear functional $\varphi
:\mathbb{H\rightarrow R},$ defined by%
\[
\varphi_{h}\left(  f\right)  :=\int_{G}f\left(  xh\right)  d\nu\left(
x\right)  ,
\]
extends uniquely to a continuous linear functional, $\bar{\varphi}_{h},$ on
$L^{p}\left(  G,\nu\right)  $ satisfying%
\[
\left\vert \bar{\varphi}_{h}\left(  f\right)  \right\vert \leq\left\Vert
f\right\Vert _{L^{p}\left(  G,\nu\right)  }\exp\left(  \frac{c\left(
kT\right)  \left(  p^{\prime}-1\right)  }{2T}d_{G_{0}}^{2}\left(  e,h\right)
\right)  \text{ for all }f\in L^{p}\left(  G,\nu\right)  .
\]
Since $L^{p}\left(  G,\nu\right)  ^{\ast}\cong L^{p^{\prime}}\left(
G,\nu\right)  ,$ there exists $J_{h}\in L^{p^{\prime}}\left(  G,\nu\right)  $
such that%
\[
\left\Vert J_{h}\right\Vert _{L^{p^{\prime}}\left(  G,\nu\right)  }\leq
\exp\left(  \frac{c\left(  kT\right)  \left(  p^{\prime}-1\right)  }%
{2T}d_{G_{0}}^{2}\left(  e,h\right)  \right)
\]
and
\[
\bar{\varphi}_{h}\left(  f\right)  =\int_{G}f\left(  x\right)  J_{h}\left(
x\right)  d\nu\left(  x\right)  \text{~}\text{ for all }~f\in L^{p}\left(
G,\nu\right)  .
\]
Restricting this formula $\mathbb{H}$ shows,%
\begin{equation}
\int_{G}f\left(  x\right)  \nu\left(  dxh^{-1}\right)  =\int_{G}f\left(
xh\right)  d\nu\left(  x\right)  =\bar{\varphi}_{h}\left(  f\right)  =\int
_{G}f\left(  x\right)  J_{h}\left(  x\right)  d\nu\left(  x\right)  ~\text{
for all }~f\in\mathbb{H}. \label{e.7.13}%
\end{equation}
Another monotone class argument (again use \cite[Theorem A.1 on p.
309]{Janson97})) shows that Eq. (\ref{e.7.13}) remains valid for all bounded
measurable functions, $f:G\rightarrow\mathbb{R}.$ Therefore, we have shown
that $J_{h}:=d\nu\circ R_{h}^{-1}/d\nu$ exists and satisfies the bound in Eq.
(\ref{e.7.10}).
\end{proof}

We now turn to the inductive limit quasi-invariance theorem. The following
result is an abstraction of the quasi-invariance result in \cite{Driver1997b}.
For related results of this type see, Fang \cite{Fang99} and Airault and
Malliavin \cite{AirMall2006}.

\begin{thm}
[Inductive Limits]\label{t.7.3}Again, let $T>0,$ $G_{0}\subset G,$ and
$\left(  G,\mathcal{B},\nu\right)  $ be as described at the start of this
section. Further assume there exists, $\left\{  G_{\alpha}\right\}
_{\alpha\in A},$ where $A$ is a directed set and for each $\alpha\in A,$
$G_{\alpha}$ is a finite dimensional uni-modular Lie subgroup of $G_{0}$ such
that $G_{\alpha}\subset G_{\beta}$ if $\alpha<\beta.$ Let $i_{\alpha
}:G_{\alpha}\rightarrow G_{0}\ $denote the smooth injection map. The following
properties are assumed to hold.

\begin{enumerate}
\item $\cup_{\alpha\in A}G_{\alpha}$ is a dense subgroup of $G_{0}.$

\item For all $f\in BC\left(  G,\mathbb{R}\right)  $ (the bounded continuous
maps from $G$ to $\mathbb{R)},$
\[
\int_{G}fd\nu=\lim_{\alpha\rightarrow\infty}\int_{G_{\alpha}}\left(  f\circ
i_{\alpha}\right)  d\nu_{\alpha},
\]
where $\nu_{\alpha}$ is the time, $T,$ heat kernel measure on $G_{\alpha}$
associated to inner product, $\left(  \cdot,\cdot\right)  _{\mathfrak{g}%
_{\alpha}},$ defined to be the restriction of $\left(  \cdot,\cdot\right)
_{\mathfrak{g}_{0}}$ to $\mathfrak{g}_{\alpha}\times\mathfrak{g}_{\alpha}.$

\item There exists $k\in\mathbb{R}$ such that $\operatorname{Ric}_{\alpha}\geq
kg_{\alpha}$ for all $\alpha\in A,$ where $\operatorname{Ric}_{\alpha}$ and
$g_{\alpha}$ are the left invariant Ricci and the metric tensors on
$G_{\alpha}$ induced by $\left(  \cdot,\cdot\right)  _{\mathfrak{g}_{\alpha}%
}.$

\item For each $\alpha\in A,$ there exits a smooth section, $s_{\alpha}%
:G_{0}\rightarrow G_{\alpha}$ (i.e. $s_{\alpha}\circ i_{\alpha}=id_{G_{\alpha
}})$ satisfying the following property. Given $\alpha_{0}\in A,$ and $k\in
C^{1}\left(  \left[  0,1\right]  ,G_{0}\right)  $ with $k\left(  0\right)
=e,$ there exists an increasing sequence, $\left\{  \alpha_{n}\right\}
_{n=1}^{\infty}\subset A$ (i.e. $\alpha_{0}<\alpha_{1}<\alpha_{2}<\dots),$
such that
\begin{equation}
\ell_{G_{0}}\left(  k\left(  \cdot\right)  \right)  =\lim_{n\rightarrow\infty
}\ell_{G_{\alpha_{n}}}\left(  s_{\alpha_{n}}\circ k\right)  . \label{e.7.14}%
\end{equation}
(We do \textbf{not} assume that $s_{\alpha}:G_{0}\rightarrow G_{\alpha}$ is a homomorphism.)
\end{enumerate}

Under these assumptions, to each $h\in G_{0},$ $\nu\circ R_{h}^{-1}$ is
absolutely continuous relative to $\nu$ and the Moreover, the Radon-Nikodym
derivative, $J_{h}:=d\left(  \nu\circ R_{h}^{-1}\right)  /d\nu,$ again
satisfies the bounds in Eq. (\ref{e.7.10}).
\end{thm}

\begin{proof}
As in the proof of Theorem \ref{t.7.2} it suffices to assume $p\in\left(
1,\infty\right)  $ throughout the proof. Let $\alpha_{0}\in A,$ $h\in
G_{\alpha_{0}},$ and $\alpha_{0}<\alpha_{1}<\alpha_{2}<\dots<\alpha_{n}<\dots$
be as in item 4. above. By Theorem \ref{t.1.6} in the form of Eq.
(\ref{e.1.9}), the Radon-Nikodym derivative, $J_{\alpha_{n}}\left(  x\right)
,$ of $\nu_{\alpha_{n}}\left(  dx\cdot s_{\alpha_{n}}\left(  h\right)
^{-1}\right)  =\nu_{\alpha_{n}}\left(  dx\cdot h^{-1}\right)  $ relative to
$\nu_{\alpha_{n}}\left(  dx\right)  $ satisfies the estimate,%
\begin{align*}
\left\Vert J_{\alpha_{n}}\right\Vert _{L^{p^{\prime}}\left(  G_{\alpha_{n}%
},\nu_{\alpha_{n}}\right)  }  &  \leq\exp\left(  \frac{c\left(  kT\right)
\left(  p^{\prime}-1\right)  }{2T}d_{G_{\alpha_{n}}}^{2}\left(  e,h^{-1}%
\right)  \right) \\
&  =\exp\left(  \frac{c\left(  kT\right)  \left(  p^{\prime}-1\right)  }%
{2T}d_{G_{\alpha_{n}}}^{2}\left(  e,h\right)  \right) \\
&  \leq\exp\left(  \frac{c\left(  kT\right)  \left(  p^{\prime}-1\right)
}{2T}\ell_{G_{\alpha_{n}}}^{2}\left(  s_{\alpha_{n}}\circ\sigma\right)
\right)  ,
\end{align*}
where $\sigma$ is any path in $C^{1}\left(  \left[  0,1\right]  ,G_{0}\right)
$ such that $\sigma\left(  0\right)  =e$ and $\sigma\left(  1\right)  =h.$
Assuming the $f\in BC\left(  G\right)  ,$ by the definition of $J_{\alpha_{n}%
}$ and H\"{o}lder's inequality,
\begin{align*}
\int_{G_{\alpha_{n}}}\left\vert f\left(  x\cdot h\right)  \right\vert
d\nu_{\alpha_{n}}\left(  x\right)   &  =\int_{G_{\alpha_{n}}}J_{\alpha_{n}%
}\left(  x\right)  \left\vert f\left(  x\right)  \right\vert d\nu_{\alpha_{n}%
}\left(  x\right) \\
&  \leq\left\Vert f\right\Vert _{L^{p}\left(  G_{\alpha_{n}},\nu_{\alpha_{n}%
}\right)  }\cdot\exp\left(  \frac{c\left(  kT\right)  \left(  p^{\prime
}-1\right)  }{2T}\ell_{G_{\alpha_{n}}}^{2}\left(  s_{\alpha_{n}}\circ
\sigma\right)  \right)  .
\end{align*}
Using the assumptions in items 2. and 4. of the theorem, we may pass to the
limit $\left(  n\rightarrow\infty\right)  $ in this inequality to find,%
\begin{equation}
\int_{G}\left\vert f\left(  x\cdot h\right)  \right\vert d\nu\left(  x\right)
\leq\left\Vert f\right\Vert _{L^{p}\left(  G,\nu\right)  }\cdot\exp\left(
\frac{c\left(  kT\right)  \left(  p^{\prime}-1\right)  }{2T}\ell_{G_{0}}%
^{2}\left(  \sigma\right)  \right)  . \label{e.7.15}%
\end{equation}

Optimizing this inequality over $\sigma\in C^{1}\left(  \left[  0,1\right]
,G_{0}\right)  $ joining $e$ to $h$ gives%
\begin{equation}
\int_{G}\left\vert f\left(  x\cdot h\right)  \right\vert d\nu\left(  x\right)
\leq\left\Vert f\right\Vert _{L^{p}\left(  G,\nu\right)  }\cdot\exp\left(
\frac{c\left(  kT\right)  \left(  p^{\prime}-1\right)  }{2T}d_{G_{0}}%
^{2}\left(  e,h\right)  \right)  . \label{e.7.16}%
\end{equation}

Up to now we have verified Eq. (\ref{e.7.16}) for any $h\in\cup_{\alpha\in
A}G_{\alpha}.$ As the latter set is dense in $G_{0},$ the dominated
convergence theorem along with the continuity of $d_{G_{0}}^{2}\left(
e,h\right)  $ in $h$ allows us to conclude that the estimate in Eq.
(\ref{e.7.16}) is valid for all $h\in G_{0}.$ Since $BC\left(  G,\mathbb{R}%
\right)  $ is dense in $L^{p}\left(  G,\nu\right)  $ (again use \cite[Theorem
A.1 on p. 309]{Janson97}) and because of Eq. (\ref{e.7.16}), the linear
functional, $\varphi_{h}:BC\left(  G\right)  \rightarrow\mathbb{R}$ defined by%
\begin{equation}
\varphi_{h}\left(  f\right)  =\int_{G}f\left(  xh\right)  d\nu\left(
x\right)  , \label{e.7.17}%
\end{equation}
has a unique extension to an element, $\bar{\varphi}_{h},$ of $L^{p}\left(
G,\nu\right)  ^{\ast}$ satisfying%
\begin{equation}
\left\vert \bar{\varphi}_{h}\left(  f\right)  \right\vert \leq\left\Vert
f\right\Vert _{L^{p}\left(  G,\nu\right)  }\cdot\exp\left(  \frac{c\left(
kT\right)  \left(  p^{\prime}-1\right)  }{2T}d_{G_{0}}^{2}\left(  e,h\right)
\right)  ~\text{ for all }~f\in L^{p}\left(  G,\nu\right)  . \label{e.7.18}%
\end{equation}
As in the latter part of the proof of Theorem \ref{t.7.2}, the estimate in Eq.
(\ref{e.7.18}) implies the existence of a function, $J_{h}\in L^{p^{\prime}%
}\left(  G,\nu\right)  ,$ such that%
\begin{equation}
\bar{\varphi}_{h}\left(  f\right)  =\int_{G}f\left(  x\right)  J_{h}\left(
x\right)  d\nu\left(  x\right)  \label{e.7.19}%
\end{equation}
and%
\[
\left\Vert J_{h}\right\Vert _{L^{p^{\prime}}\left(  G,\nu\right)  }\leq
\exp\left(  \frac{c\left(  kT\right)  \left(  p^{\prime}-1\right)  }%
{2T}d_{G_{0}}^{2}\left(  e,h\right)  \right)  .
\]
Furthermore, from Eqs. (\ref{e.7.17}) and (\ref{e.7.19}) it follows that%
\begin{equation}
\int_{G}f\left(  xh\right)  d\nu\left(  x\right)  =\int_{G}f\left(  x\right)
J_{h}\left(  x\right)  d\nu\left(  x\right)  \text{ for all }f\in BC\left(
G\right)  . \label{e.7.20}%
\end{equation}
Another monotone class argument \cite[Theorem A.1 on p. 309]{Janson97} then
shows Eq. (\ref{e.7.20}) is valid for all bounded measurable functions,
$f:G\rightarrow\mathbb{R}.$ Hence $\nu\left(  dxh^{-1}\right)  =J_{h}\left(
x\right)  \nu\left(  dx\right)  $ and $J_{h}\left(  x\right)  $ satisfies the
estimate in Eq. (\ref{e.7.10}).
\end{proof}

\begin{cor}
\label{c.7.4}Under the hypothesis of either Theorem \ref{t.7.2} or
\ref{t.7.3}, the heat kernel measure, $\nu,$ is quasi-invariant under left
translations by elements of $h\in G_{0}.$ Moreover, the Radon-Nikodym
derivative, $J_{h}^{\text{l}}:=d\left(  \nu\circ L_{h}^{-1}\right)  /d\nu$
satisfies the same bound as $d\left(  \nu\circ R_{h}^{-1}\right)  /d\nu$ which
is given in Eq. (\ref{e.7.10}).
\end{cor}

\begin{proof}
Since the heat kernel measures $\left\{  \nu_{\alpha}\right\}  _{\alpha\in A}$
on the Lie groups, $\left\{  G_{\alpha}\right\}  _{\alpha\in A},$ are
invariant under inversion, $x\rightarrow x^{-1},$ it follows that $\nu$ also
inherits this property. Hence if $f:G\rightarrow\mathbb{R}$ is a bounded
measurable function, then
\begin{align*}
\int_{G}f\left(  hx\right)  d\nu\left(  x\right)   &  =\int_{G}f\left(
hx^{-1}\right)  d\nu\left(  x\right)  =\int_{G}f\left(  \left(  xh^{-1}%
\right)  ^{-1}\right)  d\nu\left(  x\right) \\
&  =\int_{G}f\left(  x^{-1}\right)  J_{h^{-1}}\left(  x\right)  d\nu\left(
x\right)  =\int_{G}f\left(  x\right)  J_{h^{-1}}\left(  x^{-1}\right)
d\nu\left(  x\right)  ,
\end{align*}
from which it follows that $J_{h}^{\text{l}}\left(  x\right)  =J_{h^{-1}%
}\left(  x^{-1}\right)  $ for $\nu$ -- a.e. $x.$ Therefore,%
\[
\left\Vert J_{h}^{\text{l}}\right\Vert _{L^{p}\left(  \nu\right)  }=\left\Vert
J_{h^{-1}}\right\Vert _{L^{p}\left(  \nu\right)  }\leq\exp\left(
\frac{c\left(  kT\right)  \left(  p-1\right)  }{2T}d_{G_{0}}^{2}\left(
e,h^{-1}\right)  \right)
\]
which completes the proof since $d_{G_{0}}^{2}\left(  e,h^{-1}\right)
=d_{G_{0}}^{2}\left(  h,e\right)  =d_{G_{0}}^{2}\left(  e,h\right)  .$
\end{proof}

See Driver \cite{Driver1997b} for an explicit application of the projective
limit Theorem \ref{t.7.2} in the setting of loop groups and see Driver and
Gordina \cite{DG07b} for an application of the inductive limit Theorem
\ref{t.7.3} to an infinite dimensional Heisenberg group setting.

\begin{acknowledgement}
We are grateful to Alexander Grigo\'{r}yan and Laurent Saloff-Coste for their
comments and suggestions on the heat kernel bounds used in this paper. The
first author would also like to thank the Berkeley mathematics department and
the Miller Institute for Basic Research in Science for their support of this
project in its latter stages.
\end{acknowledgement}

\appendix

\section{A commutator theorem\label{s.A}}

In this section we will develop the abstract functional analytic results which
were used in the proofs of Theorems \ref{t.4.1} and \ref{t.4.3}. Results
similar to the next theorem may be found in Br\"{u}ning and Lesch
\cite{Bruning92}, Xue-Mei Li \cite{Li94,Li94b} and in Bueler \cite{Buel}.

\begin{thm}
\label{t.A.1}Let $W,X,$ and $Y$ be Hilbert spaces and $A:W\rightarrow X$ and
$B:X\rightarrow Y$ be densely defined closed (unbounded) operators such that
$\operatorname*{Ran}(A)\subset\operatorname*{Nul}(B).$ Let $Q:X\rightarrow
W\oplus Y$ be the unbounded linear operator defined by: $\mathcal{D}%
(Q)=\mathcal{D}(A^{\ast})\cap\mathcal{D}(B)$ and for $x\in\mathcal{D}(Q),$
$Qx:=(A^{\ast}x,Bx).$ Let us also define $R:W\oplus Y\rightarrow X$ by
$\mathcal{D}(R)=\mathcal{D}(A)\oplus\mathcal{D}(B^{\ast})$ and
$R(w,y):=Aw+B^{\ast}y.$ Then
\end{thm}

\begin{enumerate}
\item $\operatorname*{Ran}(A)$ and $\operatorname*{Ran}(B^{\ast})$ are orthogonal.

\item $R$ is closed.

\item $Q=R^{\ast}$ is a closed densely defined operator.

\item The operator, $L:=AA^{\ast}+B^{\ast}B,$ on $X$ is densely defined,
non-negative, and self adjoint operator. Moreover, $L:=Q^{\ast}Q.$
\end{enumerate}

\begin{proof}
We will denote all of the inner products on these Hilbert spaces by
$\langle\cdot,\cdot\rangle.$ Let $w\in\mathcal{D}(A)$ and $y\in\mathcal{D}%
(B^{\ast}).$ Since $\operatorname*{Ran}(A)\subset\operatorname*{Nul}(B),$
$0=\langle BAw,y\rangle=\langle Aw,B^{\ast}y\rangle$ which proves item 1. For
item 2., suppose that $(w_{n},y_{n})\in\mathcal{D}(R)$ are such that there
exists $(w,y)\in W\oplus Y$ and $x\in X$ such that
\begin{align*}
(w_{n},y_{n})  &  \rightarrow(w,y)\text{ as }n\rightarrow\infty\quad
\text{and}\\
R(w_{n},y_{n})  &  \rightarrow x\text{ as }n\rightarrow\infty\text{.}%
\end{align*}
We must show that $w\in\mathcal{D}(A),$ $y\in\mathcal{D}(B^{\ast})$ and that
$x=Aw+B^{\ast}y.$ We are given that $Aw_{n}+B^{\ast}y_{n}\rightarrow x$ as
$n\rightarrow\infty.$ But by the first item and the Cauchy criteria, this
implies that both $\lim_{n\rightarrow\infty}Aw_{n}$ and $\lim_{n\rightarrow
\infty}B^{\ast}y_{n}$ exist. Because both $A$ and $B^{\ast}$ are closed, this
implies that $w\in\mathcal{D}(A),$ $y\in\mathcal{D}(B^{\ast})$ and that
\[
Aw+B^{\ast}y=\lim_{n\rightarrow\infty}Aw_{n}+\lim_{n\rightarrow\infty}B^{\ast
}y_{n}=\lim_{n\rightarrow\infty}\left(  Aw_{n}+B^{\ast}y_{n}\right)  .
\]
Hence we have proved item 2.

Item 3. As $R$ is closed, it follows that $R^{\ast}$ is densely defined.
Therefore we need only show that $R^{\ast}=Q$. For this, let us recall that
$x\in\mathcal{D}(R^{\ast})$ and $R^{\ast}x=(w,y)$ iff $\langle(w,y),(w^{\prime
},y^{\prime})\rangle=\langle x,R(w^{\prime},y^{\prime})\rangle$ for all
$(w^{\prime},y^{\prime})\in\mathcal{D}(R).$ This is equivalent to the
following statements:

\begin{itemize}
\item $\langle w,w^{\prime}\rangle+\langle y,y^{\prime}\rangle=\langle
x,Aw^{\prime}+B^{\ast}y^{\prime}\rangle$ for all $w^{\prime}\in\mathcal{D}(A)$
and $y^{\prime}\in\mathcal{D}(B^{\ast}).$

\item $\langle w,w^{\prime}\rangle=\langle x,Aw^{\prime}\rangle$ and $\langle
y,y^{\prime}\rangle=\langle x,B^{\ast}y^{\prime}\rangle$ for all $w^{\prime
}\in\mathcal{D}(A)$ and $y^{\prime}\in\mathcal{D}(B^{\ast}).$

\item $x\in\mathcal{D}(A^{\ast}),$ $x\in\mathcal{D}(B^{\ast\ast}%
)=\mathcal{D}(B),$ $A^{\ast}x=w$ and $Bx=B^{\ast\ast}x=y.$

\item $x\in\mathcal{D}(Q)$ and $Qx=(w,y).$
\end{itemize}

Thus we have proved item 2. of the theorem.

Item 4. By a Theorem of Von-Neumann, \cite[Theorem X.25]{ReSi2}, $Q^{\ast}Q$
is a non-negative densely defined self adjoint operator on $X.$ So it suffices
to show that $Q^{\ast}Q=AA^{\ast}+B^{\ast}B.$

By items 2. and 3., $Q^{\ast}=R^{\ast\ast}=R.$ Therefore, $Q^{\ast}Q=RQ.$ Now
the following are equivalent:

\begin{itemize}
\item $x\in\mathcal{D}(RQ)$ and $RQx=x^{\prime}.$

\item $x\in\mathcal{D}(A^{\ast})\cap\mathcal{D}(B)$, $Qx:=(A^{\ast}%
x,Bx)\in\mathcal{D}(R),$ and $R(A^{\ast}x,Bx)=x^{\prime}.$

\item $x\in\mathcal{D}(A^{\ast})\cap\mathcal{D}(B),$ $A^{\ast}x\in
\mathcal{D}(A),$ $Bx\in\mathcal{D}(B^{\ast})$ and $AA^{\ast}x+B^{\ast
}Bx=x^{\prime}.$

\item $x\in\mathcal{D}(AA^{\ast})\cap\mathcal{D}(B^{\ast}B)$ and $AA^{\ast
}x+B^{\ast}Bx=x^{\prime}.$

\item $x\in\mathcal{D}(AA^{\ast}+B^{\ast}B)$ and $AA^{\ast}x+B^{\ast
}Bx=x^{\prime}.$
\end{itemize}

This shows $Q^{\ast}Q=AA^{\ast}+B^{\ast}B$ and thus completes the proof.
\end{proof}

\begin{thm}
[Commutator Theorem]\label{t.A.2}Let $W,$ $X,$ $Y,$ and $Z$ be Hilbert spaces
and $A:W\rightarrow X,$ $B:X\rightarrow Y,$ and $C:Y\rightarrow Z$ be densely
defined closed (unbounded) operators such that $\operatorname*{Ran}%
(A)\subset\operatorname*{Nul}(B)$ and $\operatorname*{Ran}(B)\subset
\operatorname*{Nul}(C).$ Let $L:=AA^{\ast}+B^{\ast}B$ and $S:=BB^{\ast
}+C^{\ast}C.$ Then $Be^{-tL}x=e^{-tS}Bx$ for all $x\in\mathcal{D}(B)$ and any
$t\geq0.$
\end{thm}

\begin{proof}
Let $\lambda>0.$ Observe that $BL=BB^{\ast}B$ on $\mathcal{D}(BL)=\mathcal{D}%
(AA^{\ast})\cap\mathcal{D}(BB^{\ast}B)$ and the $SB=BB^{\ast}B$ on
$\mathcal{D}(SB)=\mathcal{D}(BB^{\ast}B).$ In particular we have shown%
\begin{equation}
SB=BB^{\ast}B=BL\text{ on }\mathcal{D}(BL)=\mathcal{D}(AA^{\ast}%
)\cap\mathcal{D}(BB^{\ast}B) \label{e.A.1}%
\end{equation}
and hence,%
\begin{equation}
\left(  1+\lambda S\right)  B=B\left(  1+\lambda L\right)  \text{ on
}\mathcal{D}(BL). \label{e.A.2}%
\end{equation}
If $x\in D\left(  B\right)  $ and $x^{\prime}:=\left(  1+\lambda L\right)
^{-1}x,$ then $x^{\prime}\in D\left(  L\right)  \subset D\left(  B\right)  $
and
\[
Lx^{\prime}=\left(  1+\lambda L\right)  x^{\prime}-\lambda x^{\prime
}=x-\lambda x^{\prime}\in D\left(  B\right)  .
\]
Therefore $x^{\prime}\in D\left(  BL\right)  $ and so by Eq. (\ref{e.A.2})
applied to $x^{\prime}=\left(  1+\lambda L\right)  ^{-1}x$ we discover that,%
\[
\left(  1+\lambda S\right)  B\left(  1+\lambda L\right)  ^{-1}x=B\left(
1+\lambda L\right)  \left(  1+\lambda L\right)  ^{-1}x=Bx.
\]
Applying $\left(  1+\lambda S\right)  ^{-1}$ to both sides of this equation
shows%
\begin{equation}
B(1+\lambda L)^{-1}=(1+\lambda S)^{-1}B\text{ on }D\left(  B\right)  .
\label{e.A.3}%
\end{equation}
Multiplying Eq. (\ref{e.A.3}) on the right by $\left(  1+\lambda L\right)
^{-1}$ gives%
\[
B(1+\lambda L)^{-2}=(1+\lambda S)^{-1}B\left(  1+\lambda L\right)
^{-1}=(1+\lambda S)^{-2}B\text{ on }D\left(  B\right)  ,
\]
wherein we have used Eq. (\ref{e.A.3}) again in the second equality.
Continuing this way inductively allows us to conclude.
\begin{equation}
B(1+\lambda L)^{-n}=(1+\lambda S)^{-n}B\text{ on }D\left(  B\right)  \text{
for all }n\in\mathbb{N}. \label{e.A.4}%
\end{equation}

To complete the proof the theorem recall $e^{-tL}=s-\lim_{n\rightarrow\infty
}(1+\frac{t}{n}L)^{-n}$ and that $e^{-tS}=s-\lim_{n\rightarrow\infty}%
(1+\frac{t}{n}S)^{-n}.$ Hence, taking $\lambda=t/n$ in Eq. (\ref{e.A.4}) and
then passing to the limit allows us to conclude%
\[
\lim_{n\rightarrow\infty}B(1+\frac{t}{n}L)^{-n}x=\lim_{n\rightarrow\infty
}(1+\frac{t}{n}S)^{-n}Bx=e^{-tS}Bx~\text{ for all }~x\in D\left(  B\right)  .
\]
Since $B$ is closed, it follows that, for all $x\in D\left(  B\right)  ,$
that
\[
e^{-tL}x=\lim_{n\rightarrow\infty}(1+\frac{t}{n}L)^{-n}x\in D\left(  B\right)
\]
and
\[
Be^{-tL}x=\lim_{n\rightarrow\infty}B(1+\frac{t}{n}L)^{-n}x=e^{-tS}Bx.
\]

\end{proof}

\section{A Kato type inequality\label{s.B}}

Let $E$ be a real Euclidean vector bundle over a Riemannian manifold, $M,$
$\Gamma^{\infty}\left(  E\right)  $ $\left(  \Gamma_{c}^{\infty}\left(
E\right)  \right)  $ be the smooth (compactly supported) sections of $E,$ and
$\mathcal{H}:=L^{2}\left(  E\right)  $ be the space of square integrable
sections of $E.$ Further assume that $E$ is equipped with a metric compatible
connection, $\nabla^{E},$ and that $\square=\square^{E}$ is the associated
Bochner Laplacian on $\Gamma^{\infty}\left(  E\right)  .$ To be more explicit,
if $\left\{  e_{i}\right\}  _{i=1}^{\operatorname{rank}\left(  E\right)  }$ is
a local orthonormal frame, then%
\[
\square f=\operatorname{tr}\left(  \nabla^{T^{\ast}M\otimes E}\nabla
^{E}f\right)  =\sum_{i}\left(  \nabla_{e_{i}}^{E}\nabla_{e_{i}}^{E}%
f-\nabla_{\nabla_{e_{i}}^{TM}e_{i}}^{E}f\right)  .
\]
To simplify notation in the computations below, we will drop the superscripts,
$E$ and $TM$ from the symbols since they can be deduced from the context.

\begin{notation}
\label{n.B.1}Given a measurable section, $e:M\rightarrow E,$ and
$f\in\mathcal{H},$ let
\[
\mathrm{sgn}_{e}\left(  f\right)  :=1_{f\neq0}\frac{f}{\left\vert f\right\vert
}+1_{f=0}e=\left\{
\begin{array}
[c]{ccc}%
\frac{f}{\left\vert f\right\vert } & \text{if} & f\neq0\\
e & \text{if} & f=0
\end{array}
\right.  .
\]

\end{notation}

With this notation we have the polar decomposition, $f=\left\vert f\right\vert
\,\mathrm{sgn}_{e}\left(  f\right)  ,$ which is valid no matter what the
choice of $e.$

\begin{thm}
[Kato's Inequality]\label{t.B.2}Let $\varepsilon>0,$ $f\in\Gamma^{\infty
}\left(  E\right)  ,$ $\left\vert f\right\vert _{\varepsilon}:=\sqrt
{\left\vert f\right\vert ^{2}+\varepsilon^{2},}$ and $\hat{f}_{\varepsilon
}:=f/\left\vert f\right\vert _{\varepsilon}.$Then
\begin{align}
d\left\vert f\right\vert _{\varepsilon}  &  =\,\left\langle \hat
{f}_{\varepsilon},\nabla f\right\rangle \text{ and }\nonumber\\
\Delta_{0}\left\vert f\right\vert _{\varepsilon}  &  =\frac{1}{\left\vert
f\right\vert _{\varepsilon}}\sum_{i}\left(  \left\vert \nabla_{e_{i}%
}f\right\vert ^{2}-\left\vert \left\langle \hat{f}_{\varepsilon},\nabla
_{e_{i}}f\right\rangle \right\vert ^{2}\right)  +\left\langle \hat
{f}_{\varepsilon},\square f\right\rangle \label{e.B.1}\\
&  \geq\left\langle \hat{f}_{\varepsilon},\square f\right\rangle .
\label{e.B.2}%
\end{align}
Moreover if $\varphi\in C^{\infty}\left(  M\right)  _{+}$ and $f\in
C_{c}^{\infty}\left(  E\right)  ,$ then%
\begin{equation}
\left(  \square f,\varphi\,\mathrm{sgn}_{e}\left(  f\right)  \right)
\leq\left(  \left\vert f\right\vert ,\Delta_{0}\varphi\right)  \label{e.B.3}%
\end{equation}
where $e$ is any measurable section of $E$ such that $\left\langle \square
f\left(  x\right)  ,e\left(  x\right)  \right\rangle _{x}=0$ and $\left\vert
e\left(  x\right)  \right\vert _{x}=1$ on the set where $f=0.$
\end{thm}

\begin{proof}
This theorem is mostly a straightforward computation. (See \cite{HSU1977},
where a local coordinate version of this calculation is done.) We start by
computing the gradient of $\left\vert f\right\vert _{\varepsilon}$ as%
\[
d\left\vert f\right\vert _{\varepsilon}=\frac{1}{2\sqrt{\left\vert
f\right\vert ^{2}+\varepsilon^{2}}}d\left\vert f\right\vert ^{2}=\frac
{1}{\sqrt{\left\vert f\right\vert ^{2}+\varepsilon^{2}}}\,\left\langle
f,\nabla_{\cdot}f\right\rangle .
\]
With this in hand we have the following formula for the Hessian of $\left\vert
f\right\vert _{\varepsilon}$%
\[
\nabla d\left\vert f\right\vert _{\varepsilon}=-\left(  \left\vert
f\right\vert ^{2}+\varepsilon^{2}\right)  ^{-3/2}\left\langle f,\nabla_{\cdot
}f\right\rangle ^{2}+\frac{1}{\sqrt{\left\vert f\right\vert ^{2}%
+\varepsilon^{2}}}\left(  \,\left\langle \nabla_{\cdot}f,\nabla_{\cdot
}f\right\rangle +\left\langle f,\nabla_{\left(  \cdot,\cdot\right)  }%
^{2}f\right\rangle \right)  .
\]
Taking the trace of this result gives%
\[
\Delta_{0}\left\vert f\right\vert _{\varepsilon}=-\left(  \left\vert
f\right\vert ^{2}+\varepsilon^{2}\right)  ^{-3/2}\sum_{i}\left\vert
\left\langle f,\nabla_{e_{i}}f\right\rangle \right\vert ^{2}+\frac{1}%
{\sqrt{\left\vert f\right\vert ^{2}+\varepsilon^{2}}}\left(  \sum
_{i}\left\vert \nabla_{e_{i}}f\right\vert ^{2}+\left\langle f,\square
f\right\rangle \right)
\]
which is equivalent to Eq. (\ref{e.B.1}). Equation (\ref{e.B.2}) follows by
the Cauchy-Schwarz inequality which implies%
\[
\left\vert \nabla_{e_{i}}f\right\vert ^{2}-\left\vert \left\langle \hat
{f}_{\varepsilon},\nabla_{e_{i}}f\right\rangle \right\vert ^{2}\geq\left\vert
\nabla_{e_{i}}f\right\vert ^{2}-\left\vert \hat{f}_{\varepsilon}\right\vert
^{2}\cdot\left\vert \nabla_{e_{i}}f\right\vert ^{2}\geq0.
\]

If we now assume that $f\in\Gamma_{c}^{\infty}\left(  E\right)  $ and
$\varphi\in C^{\infty}\left(  M,[0,\infty)\right)  ,$ then multiplying Eq.
(\ref{e.B.2}) by $\varphi$ and integrating gives,
\begin{equation}
\int_{M}\left\langle \square f,\frac{f}{\left\vert f\right\vert _{\varepsilon
}}\right\rangle \varphi dV\leq\int_{M}\Delta_{0}\left\vert f\right\vert
_{\varepsilon}\varphi dV=\int_{M}\left\vert f\right\vert _{\varepsilon}%
\Delta_{0}\varphi dV \label{e.B.4}%
\end{equation}
where we have done two integrations by parts to get the last equality. Letting
$\varepsilon\downarrow0$ in Eq. (\ref{e.B.4}) then implies%
\begin{equation}
\int_{M}\left\langle \square f,\mathrm{sgn}_{0}\left(  f\right)  \right\rangle
\varphi dV\leq\int_{M}\left\vert f\right\vert \Delta_{0}\varphi dV
\label{e.B.5}%
\end{equation}
which is to say
\begin{equation}
\left\langle \square f,\mathrm{sgn}_{0}\left(  f\right)  \right\rangle
\leq\Delta_{0}\left\vert f\right\vert \text{ (in the distributional sense).}
\label{e.B.6}%
\end{equation}
If we now choose $e$ to be a measurable section of $E$ such that $\left\vert
e\right\vert =1$ and $\left\langle \square f,e\right\rangle =0,$ then
$\left\langle \square f,\mathrm{sgn}_{0}\left(  f\right)  \right\rangle
=\left\langle \square f,\mathrm{sgn}_{e}\left(  f\right)  \right\rangle $ and
we may rewrite Eqs. (\ref{e.B.5}) and (\ref{e.B.6}) as,
\[
\int_{M}\left\langle \square f,\mathrm{sgn}_{e}\left(  f\right)  \right\rangle
\varphi dV\leq\int_{M}\left\vert f\right\vert \Delta_{0}\varphi dV
\]
and%
\[
\left\langle \square f,\mathrm{sgn}_{e}\left(  f\right)  \right\rangle
\leq\Delta_{0}\left\vert f\right\vert \text{ (in the distributional sense).}%
\]
These last two equations are equivalent to Eq. (\ref{e.B.3}).
\end{proof}

\section{A local martingale\label{s.C}}

In this appendix we will continue to use the notation in Section \ref{s.5.1}
unless otherwise stated.

\begin{lem}
[Local martingale lemma]\label{l.C.1}Let $\tilde{\ell}_{t}\in\mathbb{R}$ be an
adapted continuously differentiable real valued process, $\ell_{0}\in T_{x}M,$%
\begin{equation}
\ell_{t}=Q_{t}\left[  \int_{0}^{t}Q_{\tau}^{-1}\left(  \frac{d}{d\tau}%
\tilde{\ell}_{\tau}\right)  db_{\tau}+\ell_{0}\right]  , \label{e.C.1}%
\end{equation}
$a\in\Omega_{c}^{1}\left(  M\right)  ,$ and%
\begin{equation}
Z_{t}:=\left(  a_{t}\left(  \Sigma_{t}\right)  \circ//_{t}\right)  \ell
_{t}-\delta a_{t}\left(  \Sigma_{t}\right)  \tilde{\ell}_{t}, \label{e.C.2}%
\end{equation}
be as in Eq. (\ref{e.5.10}). Then $Z_{t}$ is a local martingale whose It\^{o}
differential is given by%
\begin{equation}
dZ_{t}=\left(  \nabla_{//_{t}db_{t}}a_{t}\right)  \left(  \Sigma_{t}\right)
\circ//_{t}\ell_{t}+\left(  a_{t}\left(  \Sigma_{t}\right)  \circ
//_{t}\right)  \left(  \frac{d}{dt}\tilde{\ell}_{t}\right)  db_{t}-\left(
\nabla_{//_{t}db_{t}}a_{t}\right)  \left(  \Sigma_{t}\right)  \tilde{\ell}%
_{t}. \label{e.C.3}%
\end{equation}

\end{lem}

\begin{proof}
The proof of this lemma is purely a computation. For the sake of the reader's
understanding we will give a slightly inefficient proof designed to motivate
the form of $Z_{t}$ in Eq. (\ref{e.5.10}). Let $a_{t}$ be as in Eq.
(\ref{e.5.9}) and then set%
\[
N_{t}:=a_{t}\left(  \Sigma_{t}\right)  \circ//_{t}.
\]
Then by It\^{o}'s lemma in Eq. (\ref{e.5.2}), Theorem \ref{t.4.1}, and Bochner
identity in Eq. (\ref{e.4.3}), we find%
\begin{align}
dN_{t}  &  =\left(  \nabla_{//_{t}db_{t}}a_{t}\right)  \left(  \Sigma
_{t}\right)  \circ//_{t}+\frac{1}{2}\left(  \left(  \square-\Delta\right)
a_{t}\left(  \Sigma_{t}\right)  \right)  \circ//_{t}dt\nonumber\\
&  =\left(  \nabla_{//_{t}db_{t}}a_{t}\right)  \left(  \Sigma_{t}\right)
\circ//_{t}+\frac{1}{2}\left[  a_{t}\left(  \Sigma_{t}\right)  \circ
\operatorname{Ric}\circ//_{t}\right]  dt. \label{e.C.4}%
\end{align}
Also by It\^{o}'s lemma in Eq. (\ref{e.5.1}) and item 4. of Theorem
\ref{t.4.1},
\begin{align}
d\left[  \delta a_{t}\left(  \Sigma_{t}\right)  \right]   &  =d\left[  \left(
e^{\left(  T-t\right)  \bar{\Delta}_{0}/2}\delta a\right)  \left(  \Sigma
_{t}\right)  \right] \nonumber\\
&  =\left(  \nabla_{//_{t}db_{t}}\left[  e^{\left(  T-t\right)  \bar{\Delta
}_{0}/2}\delta a\right]  \right)  \left(  \Sigma_{t}\right)  =\left(
\nabla_{//_{t}db_{t}}\left[  \delta a_{t}\right]  \right)  \left(  \Sigma
_{t}\right)  . \label{e.C.5}%
\end{align}

Now suppose $\ell_{t}\in T_{x}M$ and $\tilde{\ell}_{t}\in\mathbb{R}$ are
arbitrary continuous Brownian semi-martingales such that%
\[
d\ell_{t}=\alpha_{t}\,db_{t}+\beta_{t}\,dt\text{ and }d\tilde{\ell}_{t}%
=\tilde{\alpha}_{t}\,db_{t}+\tilde{\beta}_{t}\,dt
\]
with $\alpha_{t},$ $\beta_{t},$ $\tilde{\alpha}_{t},$ and $\tilde{\beta}_{t}$
being continuous adapted processes with values in $\operatorname*{End}\left(
T_{x}M\right)  ,$ $T_{x}M,$ $T_{x}M^{\ast},$ and $\mathbb{R}$ respectively and
let
\begin{equation}
Z_{t}=N_{t}\ell_{t}-\left(  \delta a_{t}\right)  \left(  \Sigma_{t}\right)
\tilde{\ell}_{t}. \label{e.C.6}%
\end{equation}
Making use of Eqs. (\ref{e.C.4}) and (\ref{e.C.5}), the It\^{o} differential
of $Z$ in Eq. (\ref{e.C.6}) is,%
\begin{align}
dZ_{t}  &  =\left(  \nabla_{//_{t}db_{t}}a_{t}\right)  \left(  \Sigma
_{t}\right)  \circ//_{t}\ell_{t}+\frac{1}{2}\left[  a_{t}\left(  \Sigma
_{t}\right)  \circ\operatorname{Ric}\circ//_{t}\ell_{t}\right]  dt\nonumber\\
&  +\left(  a_{t}\left(  \Sigma_{t}\right)  \circ//_{t}\right)  \left[
\alpha_{t}db_{t}+\beta_{t}\,dt\right]  +\left(  \nabla_{//_{t}e_{i}}%
a_{t}\right)  \left(  \Sigma_{t}\right)  \circ//_{t}\alpha_{t}e_{i}%
dt\nonumber\\
&  -\left(  \nabla_{//_{t}db_{t}}\left[  \delta a_{t}\right]  \right)  \left(
\Sigma_{t}\right)  \tilde{\ell}_{t}-\delta a_{t}\left(  \Sigma_{t}\right)
\left[  \tilde{\alpha}_{t}db_{t}+\tilde{\beta}_{t}dt\right] \nonumber\\
&  -\left(  \nabla_{//_{t}e_{i}}\left[  \delta a_{t}\right]  \right)
\tilde{\alpha}_{t}e_{i}dt\nonumber\\
&  =\left(  \nabla_{//_{t}db_{t}}a_{t}\right)  \left(  \Sigma_{t}\right)
\circ//_{t}\ell_{t}+\left(  a_{t}\left(  \Sigma_{t}\right)  \circ
//_{t}\right)  \alpha_{t}db_{t}\nonumber\\
&  -\left(  \nabla_{//_{t}db_{t}}\left[  \delta a_{t}\right]  \right)  \left(
\Sigma_{t}\right)  \tilde{\ell}_{t}-\delta a_{t}\left(  \Sigma_{t}\right)
\tilde{\alpha}_{t}db_{t}\nonumber\\
&  +\left(
\begin{array}
[c]{c}%
\frac{1}{2}\left[  a_{t}\left(  \Sigma_{t}\right)  \circ\operatorname{Ric}%
\circ//_{t}\ell_{t}\right]  +\left(  a_{t}\left(  \Sigma_{t}\right)
\circ//_{t}\right)  \beta_{t}+\left(  \nabla_{//_{t}e_{i}}a_{t}\right)
\left(  \Sigma_{t}\right)  \circ//_{t}\alpha_{t}e_{i}\\
-\delta a_{t}\left(  \Sigma_{t}\right)  \tilde{\beta}_{t}-\left(
\nabla_{//_{t}e_{i}}\left[  \delta a_{t}\right]  \right)  \tilde{\alpha}%
_{t}e_{i}%
\end{array}
\right)  dt. \label{e.C.7}%
\end{align}

Our goal is to choose $\alpha_{t},$ $\beta_{t},$ $\tilde{\alpha}_{t},$ and
$\tilde{\beta}_{t}$ in such as way that $Z_{t}$ is a local martingale. To do
this we need to make the term in the parenthesis in Eq. (\ref{e.C.7}) vanish.
Grouping the terms according to the number of derivatives on $a_{t},$ the term
in parenthesis in Eq. (\ref{e.C.7}) will vanish provided%
\begin{align*}
\frac{1}{2}\left[  a_{t}\left(  \Sigma_{t}\right)  \circ\operatorname{Ric}%
\circ//_{t}\ell_{t}\right]  +\left(  a_{t}\left(  \Sigma_{t}\right)
\circ//_{t}\right)  \beta_{t}  &  =0,\\
\left(  \nabla_{//_{t}e_{i}}a_{t}\right)  \left(  \Sigma_{t}\right)
\circ//_{t}\alpha_{t}e_{i}-\delta a_{t}\left(  \Sigma_{t}\right)  \tilde
{\beta}_{t}  &  =0,\\
\text{and \qquad}\left(  \nabla_{//_{t}e_{i}}\left[  \delta a_{t}\right]
\right)  \tilde{\alpha}_{t}e_{i}  &  =0.
\end{align*}
Moreover because of Eq. (\ref{e.4.1}), these equations may be satisfied by
choosing $\tilde{\alpha}\equiv0$ (so that $\tilde{\ell}_{t}$ is differentiable
and $\frac{d\tilde{\ell}_{t}}{dt}=\tilde{\beta}_{t}),$
\[
\beta_{t}=-\frac{1}{2}//_{t}^{-1}\operatorname{Ric}//_{t}\ell_{t}=:-\frac
{1}{2}\operatorname{Ric}^{//_{t}}\ell_{t},
\]
and%
\[
\alpha_{t}=\tilde{\beta}_{t}I_{T_{x}M}=\frac{d\tilde{\ell}_{t}}{dt}I_{T_{x}%
M}.
\]
Thus we have shown,%
\[
Z_{t}:=\left(  a_{t}\left(  \Sigma_{t}\right)  \circ//_{t}\right)  \ell
_{t}-\delta a_{t}\left(  \Sigma_{t}\right)  \tilde{\ell}_{t},
\]
is a local martingale provided $\tilde{\ell}_{t}$ is an adapted $C^{1}$ --
process and $\ell$ solves%
\begin{equation}
d\ell_{t}=\frac{d\tilde{\ell}_{t}}{dt}db_{t}-\frac{1}{2}\operatorname{Ric}%
^{//_{t}}\ell_{t}\,dt. \label{e.C.8}%
\end{equation}
To solve this equation for $\ell_{t},$ let $Q_{t}$ solve the ODE in Eq.
(\ref{e.5.3}) and write $\ell_{t}=Q_{t}k_{t}$ where $k_{t}:=Q_{t}^{-1}\ell
_{t}.$ Plugging this expression for $\ell_{t}$ into Eq. (\ref{e.C.8}) using,%
\[
d\ell_{t}=-\frac{1}{2}\operatorname{Ric}^{//_{t}}Q_{t}k_{t}dt+Q_{t}dk_{t},
\]
implies,
\[
-\frac{1}{2}\operatorname{Ric}^{//_{t}}Q_{t}k_{t}dt+Q_{t}dk_{t}=\frac
{d\tilde{\ell}_{t}}{dt}db_{t}-\frac{1}{2}\operatorname{Ric}^{//_{t}}Q_{t}%
k_{t}\,dt
\]
from which we learn, $dk_{t}=Q_{t}^{-1}\frac{d\tilde{\ell}_{t}}{dt}db_{t}.$
Integrating this equation and multiplying the result on the left by $Q_{t}$
gives Eq. (\ref{e.C.1}). Equation (\ref{e.C.3}) now follows from Eq.
(\ref{e.C.7}) with $\tilde{\alpha}=0$ and $\alpha_{t}=\frac{d\tilde{\ell}_{t}%
}{dt}I_{T_{x}M}.$
\end{proof}

\section{Wang's dimension free Harnack inequality\label{s.D}}

Suppose that $p_{T}\left(  \cdot,\cdot\right)  >0$ is the heat kernel at time
$T>0$ on a complete connected Riemannian manifold $\left(  M\right)  $ and for
measurable $f:M\rightarrow\lbrack0,\infty),$ let
\[
\left(  P_{T}f\right)  \left(  x\right)  :=\int_{M}p_{T}\left(  x,y\right)
f\left(  y\right)  dV\left(  y\right)  .
\]
Hence if $f\in L^{2}\left(  V\right)  ,$ then $P_{T}f=e^{T\bar{\Delta}_{0}%
/2}f.$ The following lemma reflects the fact that $\left(  L^{p}\right)
^{\ast}$ and $L^{p^{\prime}}$ are isometrically isomorphic Banach spaces for
$1<p<\infty$ and $p^{\prime}=p/\left(  p-1\right)  $ -- the conjugate exponent
to $p.$

\begin{lem}
\label{l.D.1}Let $x,y\in M,$ $T>0,$ $p\in(1,\infty),$ and $C\in(0,\infty].$
Then%
\begin{equation}
\left[  \left(  P_{T}f\right)  \left(  x\right)  \right]  ^{p}\leq
C^{p}\left(  P_{T}f^{p}\right)  \left(  y\right)  \text{ for all }f\geq0
\label{e.D.1}%
\end{equation}
if and only if%
\begin{equation}
\left(  \int_{M}\left[  \frac{p_{T}\left(  x,z\right)  }{p_{T}\left(
y,z\right)  }\right]  ^{p^{\prime}}p_{T}\left(  y,z\right)  dV\left(
z\right)  \right)  ^{1/p^{\prime}}\leq C. \label{e.D.2}%
\end{equation}

\end{lem}

\begin{proof}
Since%
\[
\left(  P_{T}f\right)  \left(  x\right)  =\int_{M}\frac{p_{T}\left(
x,z\right)  }{p_{T}\left(  y,z\right)  }f\left(  z\right)  p_{T}\left(
y,z\right)  dV\left(  z\right)  ,
\]
if $d\mu\left(  z\right)  :=p_{T}\left(  y,z\right)  dV\left(  z\right)  $ and
$g\left(  x\right)  :=\frac{p_{T}\left(  x,\cdot\right)  }{p_{T}\left(
y,\cdot\right)  },$ then
\begin{equation}
\left(  P_{T}f\right)  \left(  x\right)  =\int_{M}f\left(  x\right)  g\left(
x\right)  d\mu\left(  x\right)  . \label{e.D.3}%
\end{equation}
Since $g\geq0$ and $L^{p}\left(  \mu\right)  ^{\ast}$ is isomorphic to
$L^{p^{\prime}}\left(  \mu\right)  ^{\ast}$ under the pairing in Eq.
(\ref{e.D.3}), it follows that
\[
\left\Vert g\right\Vert _{L^{p^{\prime}}\left(  \mu\right)  }=\sup_{f\geq
0}\frac{\int_{M}f\left(  x\right)  g\left(  x\right)  d\mu\left(  x\right)
}{\left\Vert f\right\Vert _{L^{p}\left(  \mu\right)  }}=\sup_{f\geq0}%
\frac{\left(  P_{T}f\right)  \left(  x\right)  }{\left[  \left(  P_{T}%
f^{p}\right)  \left(  y\right)  \right]  ^{1/p}}.
\]
The last equation may be written more explicitly as,%
\[
\left(  \int_{M}\left[  \frac{p_{T}\left(  x,z\right)  }{p_{T}\left(
y,z\right)  }\right]  ^{p^{\prime}}p_{T}\left(  y,z\right)  dV\left(
z\right)  \right)  ^{1/p^{\prime}}=\sup_{f\geq0}\frac{\left(  P_{T}f\right)
\left(  x\right)  }{\left[  \left(  P_{T}f^{p}\right)  \left(  y\right)
\right]  ^{1/p}},
\]
and from this equation the lemma easily follows.
\end{proof}

The following theorem appears in \cite{Wang97a,Wang2004} -- also see
also see \cite{}.

\begin{thm}
[Wang's Harnack inequality]\label{t.D.2}Suppose that $M$ is a complete
connected Riemannian manifold such that $\operatorname*{Ric}\geq kI$ for some
$k\in\mathbb{R}.$ Then for all $p>1,$ $f\geq0,$ $T>0,$ and $x,y\in M,$ we
have
\begin{equation}
\left(  P_{T}f\right)^{p}\left(  y\right)  \leq\left(
P_{T}f^{p}\right) \left(  z\right)  \exp\left(
p^{\prime}\frac{k}{e^{kT}-1}d^{2}\left(
y,z\right)  \right), \label{e.D.4}%
\end{equation}
where $p^{\prime}=p/\left(  p-1\right)  $ is the conjugate exponent to $p.$
\end{thm}

In applying Wang's results the reader should use $k=-K,$ $V\equiv0,$ and
replace $T$ by $T/2$ since Wang's generator is $\Delta$ rather than
$\Delta/2.$

\begin{cor}
\label{c.D.3}Let $\left(  M,g\right)  $ be a complete Riemannian manifold such
that $\operatorname{Ric}\geq kI$ for some $k\in\mathbb{R}.$ Then for every
$y,z\in M$ and $p\in\lbrack1,\infty),$%
\begin{equation}
\left(  \int_{M}\left[  \frac{p_{T}\left(  y,x\right)  }{p_{T}\left(
z,x\right)  }\right]  ^{p}p_{T}\left(  z,x\right)  dV\left(  x\right)
\right)  ^{1/p}\leq\exp\left(  \frac{c\left(  kT\right)  \left(  p-1\right)
}{2T}d^{2}\left(  y,z\right)  \right)  \label{e.D.5}%
\end{equation}
where $c\left(  \cdot\right)  $ is defined as in Eq. (\ref{e.1.7}),
$p_{t}\left(  x,y\right)  $ is the heat kernel on $M$ and $d\left(
y,z\right)  $ is the Riemannian distance from $x$ to $y$ for $x,y\in M.$
\end{cor}

\begin{proof}
From Lemma \ref{l.D.1} and Theorem \ref{t.D.2} with
\[
C=\exp\left(  \frac{p^{\prime}}{p}\frac{k}{e^{kT}-1}d^{2}\left(  y,z\right)
\right)  =\exp\left(  \frac{1}{p-1}\frac{k}{e^{kT}-1}d^{2}\left(  y,z\right)
\right)  ,
\]
it follows that it follows that
\[
\left(  \int_{M}\left[  \frac{p_{T}\left(  x,z\right)  }{p_{T}\left(
y,z\right)  }\right]  ^{p^{\prime}}p_{T}\left(  y,z\right)  dV\left(
z\right)  \right)  ^{1/p^{\prime}}\leq\exp\left(  \frac{1}{p-1}\frac{k}%
{e^{kT}-1}d^{2}\left(  y,z\right)  \right)  .
\]
Using $p-1=\left(  p^{\prime}-1\right)  ^{-1}$ and then interchanging the
roles of $p$ and $p^{\prime}$ gives Eq. (\ref{e.D.5}).
\end{proof}

For comparison sake, recall that the classical Li - Yau Harnack inequality
(see Li and Yau \cite{Li-Yau86} and Davies \cite[Theorem 5.3.5]{Davies89})
states if $\alpha>1,$ $s>0,$ and $\operatorname{Ric}\geq-K$ for some $K\geq0,$
then%
\begin{equation}
\frac{p_{t}\left(  y,x\right)  }{p_{t+s}\left(  z,x\right)  }\leq\left(
\frac{t+s}{t}\right)  ^{d\alpha/2}\exp\left(  \frac{\alpha d^{2}\left(
y,z\right)  }{2s}+\frac{d\cdot\alpha Ks}{8\left(  \alpha-1\right)  }\right)  ,
\label{e.D.6}%
\end{equation}
for all $x,y,z\in M^{d}$ and $t>0.$ However when $s=0,$ Eq. (\ref{e.D.6})
gives no information on $p_{t}\left(  y,x\right)  /p_{t}\left(  z,x\right)  $
when $y\neq z.$

\begin{rem}
\label{r.D.4}Since our heat equation is determined by $\Delta_{0}/2$ rather
than $\Delta_{0},$ the reader should replace $t$ and $s$ by $t/2$ and $s/2$
when applying the results in \cite{Li-Yau86,Davies89}.
\end{rem}

\section{Consequences of Hamilton's estimates\label{s.E}}

Let $T\in\left(  0,\infty\right)  ,$ $M$ $\left(  d=\dim\left(  M\right)
\right)  $ be a complete Riemannian manifold with $\operatorname{Ric}\geq-KI$
for some $K\geq0,$ and let $V\left(  x,r\right)  :=\operatorname*{Vol}\left(
B\left(  x,r\right)  \right)  $ be the volume of the ball, $B\left(
x,r\right)  ,$ centered at $x\in M$ with radius $r>0.$ Suppose, for $0\leq
t\leq t_{1},$ that $u\left(  t,x\right)  $ is a positive solution to the heat
equation, $\frac{\partial}{\partial t}u=\frac{1}{2}\Delta u.$ The Hamilton
type gradient bounds \cite{Hamilton93,Souplet-Zhang06,Kotschwar07} state if
\[
m:=\sup\left\{  u\left(  t,x\right)  :0\leq t\leq t_{1},~x\in M\right\}
\]
then
\begin{equation}
t|\nabla\log(u\left(  t,x\right)  )|^{2}\leq2(1+Kt)\log(m/u\left(  t,x\right)
)\text{ for all }\left(  t,x\right)  \in\left[  0,t_{1}\right]  \times M.
\label{e.E.1}%
\end{equation}
The standard heat kernel bounds (see for example Theorems 5.6.4, 5.6.6, and
5.4.12 in Sallof-Coste \cite{Saloff-Coste2002} and for more detailed bounds
see \cite{Li-Yau86,Davies89,Saloff-Coste92,Davies93,Grigoryan1995}) which
state there exist constants, $c=c\left(  K,d,T\right)  $ and $C=C\left(
K,d,T\right)  ,$ such that,%
\begin{equation}
\frac{c}{V\left(  x,\sqrt{t/2}\right)  }\exp\left(  -C\frac{d^{2}(x,y)}%
{t}\right)  \leq p(t,x,y)\leq\frac{C}{V\left(  x,\sqrt{t/2}\right)  }%
\exp\left(  -c\frac{d^{2}(x,y)}{t}\right)  , \label{e.E.2}%
\end{equation}
for all $x,y\in M$ and $t\in(0,T].$

Let $s\in(0,T],$ $o\in M,$ $t_{1}=s/2$ and $u\left(  t,x\right)
=p_{s/2+t}\left(  o,x\right)  .$ Combining Eqs. (\ref{e.E.1}) and
(\ref{e.E.2}) then shows,%
\begin{equation}
t|\nabla_{x}\log p_{s/2+t}\left(  o,x\right)  )|^{2}\leq2(1+Kt)\log\left(
\frac{C}{c}\frac{V\left(  0,\sqrt{s/4+t/2}\right)  }{V\left(  o,\sqrt
{s/4}\right)  }\exp\left(  C\frac{d^{2}(o,y)}{s/2+t}\right)  \right)  .
\label{e.E.3}%
\end{equation}
Taking $t=s/2$ in Eq. (\ref{e.E.3}) and then replacing $s$ by $t$ in the
resulting inequality implies,%
\begin{equation}
\frac{t}{2}|\nabla_{x}\log p_{t}\left(  o,x\right)  )|^{2}\leq2(1+K\frac{t}%
{2})\log\left(  \frac{C}{c}\frac{V\left(  0,\sqrt{t/2}\right)  }{V\left(
o,\sqrt{t/4}\right)  }\exp\left(  C\frac{d^{2}(o,y)}{t}\right)  \right)  .
\label{e.E.4}%
\end{equation}
Using the volume estimate (see \cite{CGT82} and \cite[Theorem 5.6.4]%
{Saloff-Coste2002}),%
\[
\frac{V(x,\sigma)}{V(x,s)}\leq\left(  \frac{\sigma}{s}\right)  ^{d}\exp\left(
\sqrt{\left(  d-1\right)  K}\sigma\right)  ~\forall~x\in M\text{ and }0\leq
s<\sigma,
\]
it follows that%
\begin{equation}
\frac{V(x,\sqrt{t/2})}{V(x,\sqrt{t/4})}\leq2^{d/2}\exp\left(  \sqrt{\left(
d-1\right)  Kt/2}\right)  \leq2^{d/2}\exp\left(  \sqrt{\left(  d-1\right)
KT/2}\right)  . \label{e.E.5}%
\end{equation}
Combining Eqs. (\ref{e.E.4}) and (\ref{e.E.5}) then allows us to conclude that
there exist constants, $c_{1}$ and $c_{2}$ depending on $T,K,$ and $d\ $such
that%
\begin{equation}
\left\vert \nabla_{x}\log p_{t}\left(  o,x\right)  )\right\vert \leq\left(
\frac{c_{1}}{\sqrt{t}}+c_{2}\frac{d\left(  o,x\right)  }{t}\right)  ~\text{
for all }~t\in(0,T]\text{ and }o,x\in M. \label{e.E.6}%
\end{equation}
For this estimate in the compact case with its relations to stochastic
analysis, see \cite{Driver1994b,Mal-stroock96,Stroock96,Stroock98,Hsu99a}.

\begin{prop}
\label{p.E.1}Continuing the notation and assumptions used above, there exist
constants, $C_{1}\left(  d,K\right)  $ and $C_{2}\left(  d,K,t\right)  $ such
that,%
\begin{equation}
\int_{M}\exp\left(  \lambda\left\vert \nabla_{x}\log p_{t}\left(  o,x\right)
)\right\vert \right)  p_{t}\left(  o,x\right)  dx\leq C\left(  d,K,t\right)
\exp\left(  C\left(  d,K\right)  \lambda^{2}/t\right)  \label{e.E.7}%
\end{equation}
for all $o\in M$ and $t\in(0,T].$
\end{prop}

\begin{proof}
Let $v\left(  r\right)  :=\operatorname{Vol}\left(  B\left(  o,r\right)
\right)  ,$ $\kappa:=\sqrt{K/\left(  d-1\right)  },$ $\gamma:=\left(
d-1\right)  \kappa=\sqrt{K\left(  d-1\right)  },$ and $\omega_{d-1}$ be the
volume of the standard $d-1$ sphere. Using Bishop's comparison theorem (see
\cite{CheegerEbinBook,schoen_yau}) which states,
\begin{equation}
dv\left(  r\right)  \leq\omega_{d-1}\left(  \frac{\sinh\kappa r}{\kappa
}\right)  ^{d-1}dr\leq\left(  \frac{\omega_{d-1}}{2\kappa}\right)
^{d-1}e^{\kappa\left(  d-1\right)  r}dr, \label{e.E.8}%
\end{equation}
along with the estimates in Eqs. (\ref{e.E.2}) and (\ref{e.E.6}), we have%
\begin{align}
&  \int_{M}\exp\left(  \lambda\left\vert \nabla_{x}\log p_{t}\left(
o,x\right)  )\right\vert \right)  p_{t}\left(  o,x\right)  dx\nonumber\\
&  \leq Ct^{-d/2}\int_{0}^{\infty}\exp\left(  \lambda\left(  \frac{c_{1}%
}{\sqrt{t}}+c_{2}\frac{r}{t}\right)  \right)  \exp\left(  -\frac{C}{2t}%
r^{2}\right)  dv\left(  r\right) \nonumber\\
&  \leq C\left(  \frac{\omega_{d-1}}{2\kappa}\right)  ^{d-1}t^{-d/2}\int
_{0}^{\infty}\exp\left(  \lambda\left(  \frac{c_{1}}{\sqrt{t}}+c_{2}\frac
{r}{t}\right)  \right)  \exp\left(  -\frac{C}{2t}r^{2}\right)  e^{\gamma
r}dr\label{e.E.9}\\
&  =C(d,K,T)t^{-d/2}\exp\left(  \lambda\frac{c_{1}}{\sqrt{t}}\right)  \int
_{0}^{\infty}\exp\left(  \left(  \gamma+\lambda\frac{c_{2}}{t}\right)
r\right)  \exp\left(  -\frac{C}{2t}r^{2}\right)  dr. \label{e.E.10}%
\end{align}
Equation (\ref{e.E.7}) follows easily from Eq. (\ref{e.E.10}) and the
following two estimates
\[
c_{1}\frac{\lambda}{\sqrt{t}}\leq\frac{1}{2}\left(  c_{1}^{2}+\frac
{\lambda^{2}}{2t}\right)
\]
and%
\begin{align}
\int_{0}^{\infty}\exp\left(  \left(  \gamma+\lambda\frac{c_{2}}{t}\right)
r\right)   &  \exp\left(  -\frac{C}{2t}r^{2}\right)  dr\nonumber\\
&  \leq\int_{-\infty}^{\infty}\exp\left(  \left(  \gamma+\lambda\frac{c_{2}%
}{t}\right)  r\right)  \exp\left(  -\frac{C}{2t}r^{2}\right)  dr\nonumber\\
&  \qquad\qquad=\sqrt{2\pi t/C}\exp\left(  \frac{t}{2C}\left(  \gamma
+\lambda\frac{c_{2}}{t}\right)  ^{2}\right)  . \label{e.E.11}%
\end{align}

\end{proof}

\begin{rem}
\label{r.E.2}When $M=\mathbb{R}^{d},$ using Laplace asymptotics, one may show;%
\[
\lim_{d\rightarrow\infty}e^{-\frac{\lambda}{\sqrt{t}}\sqrt{d-1}}%
\int_{\mathbb{R}^{d}}\exp\left(  \lambda\left\vert \nabla_{x}\log p_{t}\left(
o,x\right)  )\right\vert \right)  p_{t}\left(  o,x\right)  dx=e^{\lambda
^{2}/4t}~\forall~t,\lambda>0.
\]
In particular, this implies that we can not take both $C\left(  d,0,t\right)
$ and $C\left(  d,0\right)  $ in Eq. (\ref{e.E.7}) to be independent of the
dimension, $d=\dim\left(  M\right)  .$
\end{rem}

\input{end}

\def\cprime{$'$}
\providecommand{\bysame}{\leavevmode\hbox to3em{\hrulefill}\thinspace}
\providecommand{\MR}{\relax\ifhmode\unskip\space\fi MR }
\providecommand{\MRhref}[2]{%
  \href{http://www.ams.org/mathscinet-getitem?mr=#1}{#2}
}
\providecommand{\href}[2]{#2}

\end{document}

%% file: har-rev-title.tex
\def\H{\mathcal{H}}
\def\P{\mathcal{P}}
\def\WC{W_{\mathbb C}}
\def\HC{H_{\mathbb C}}

\def\hfootnote#1{}

\def\enddoc{\end{document}}

\title[Integrated Harnack Inequalities]{Integrated Harnack Inequalities on Lie
Groups}
\author[Driver]{Bruce K. Driver{$^{\dagger }$}}
\thanks{\footnotemark {$^\dagger$}This research was supported in part by NSF Grant DMS-0504608 and the Miller Institute at the University of California, at Berkeley.}
\address{Department of Mathematics, 0112\\
University of California, San Diego \\
La Jolla, CA 92093-0112 } \email{driver{@}euclid.ucsd.edu}
\author[Gordina]{Maria Gordina{$^{\ddag }$}}
\thanks{\footnotemark {$^\ddag$}Research was supported in part by NSF Grant DMS-0706784.}
\address{Department of Mathematics\\
University of Connecticut\\
Storrs, CT 06269, U.S.A. } \email{gordina@math.uconn.edu}



\date{\today \ \emph{File:\jobname{.tex}}}

\begin{abstract}
We show that the logarithmic derivatives of the convolution
heat kernels on a uni-modular Lie group are exponentially integrable. This
result is then used to prove an ``integrated'' Harnack inequality for these heat
kernels. It is shown that this integrated Harnack inequality is equivalent to
a version of Wang's Harnack inequality. (A key feature of all of these
inequalities is that they are dimension independent.) Finally, we show these inequalities imply quasi-invariance properties of heat kernel measures
for two classes of infinite dimensional ``Lie'' groups. 
\end{abstract}

\maketitle

\tableofcontents


\renewcommand{\contentsname}{Table of Contents}

%% file: end.tex
\bibliographystyle{amsplain}
\bibliography{Heis}